\RequirePackage{fix-cm}
\documentclass{svjour3}                     
\smartqed  
\usepackage{graphicx}
\usepackage{indentfirst}
\usepackage{graphicx}
\usepackage{epstopdf}
\usepackage{epsfig}
\usepackage{overpic}
\usepackage{overcite}
\usepackage{array}
\usepackage{color}
\usepackage{multirow}
\usepackage{float}
\usepackage{subfigure} 
\usepackage{amsmath,amssymb}
\usepackage[mathscr]{euscript}
\usepackage{latexsym,bm}
\usepackage{booktabs}
\usepackage{bbm}
\usepackage{diagbox}
 \usepackage[all]{xy}
 \usepackage{color}
\usepackage{diagbox}
\usepackage{algorithm}  
\usepackage{algpseudocode}  
\usepackage{hyperref}
\usepackage{graphicx}
\usepackage{mathtools}
\newcommand\D{\textup{d}}

\def\mi{\mathrm{i}}
\def\me{\mathrm{e}}

\def\bx{\bm{x}}

\def\bk{\bm{k}}

\def\by{\bm{y}}

\def\bn{\bm{n}}

\def\bzeta{\bm{\zeta}}

\def\pdo{{\rm \Psi} \textup{DO}}
\begin{document}

\title{Performance evaluations on the parallel CHAracteristic-Spectral-Mixed (CHASM) scheme}
\subtitle{As supplementary materials of ``A characteristic-spectral-mixed scheme for six-dimensional Wigner-Coulomb dynamics''}

\titlerunning{Performance evaluations on CHASM}        

\author{Yunfeng Xiong
        \and
        Yong Zhang
        \and
        Sihong Shao
}


\institute{Y. Xiong \at
              CAPT, LMAM and School of Mathematical Sciences, Peking University, Beijing, China. \\
              \email{xiongyf@math.pku.edu.cn}           
           \and
           Y. Zhang \at
              The School of Mathematics, Tianjin University, Tianjin, China. \\
               \email{ sunny5zhang@163.com}
           \and
           S. Shao \at
               CAPT, LMAM and School of Mathematical Sciences, Peking University, Beijing, China. \\
               To whom correspondence should be addressed. \\
                \email{sihong@math.pku.edu.cn}
              }
\date{\today}

\maketitle

\begin{abstract}

Performance evaluations on the deterministic algorithms for 6-D problems are rarely found in literatures except some recent advances in the Vlasov and Boltzmann community [Dimarco et al. (2018), Kormann et al. (2019)], due to the extremely high complexity. Thus a detailed comparison among various techniques shall be useful to the researchers in the related fields. We try to make a thorough evaluation on a parallel CHAracteristic-Spectral-Mixed (CHASM) scheme to support its usage. CHASM utilizes the cubic B-spline expansion in the spatial space and spectral expansion in the momentum space, which many potentially overcome the computational burden in solving classical and quantum kinetic equations in 6-D phase space. Our purpose is three-pronged. First, we would like show that by imposing some effective Hermite boundary conditions, the local cubic spline can approximate to the global one as accurately as possible. Second, we will illustrate the necessity of adopting the truncated kernel method in calculating the pseudodifferential operator with a singular symbol, since the widely used pseudo-spectral method [Ringhofer (1990)]  might fail to properly tackle the singularity. Finally, we make a comparison among non-splitting Lawson schemes and Strang operator splitting. Our numerical results demonstrate the advantage of the one-stage Lawson predictor-corrector scheme over multi-stage ones as well as the splitting scheme in both accuracy and stability.

\end{abstract}

\newsavebox{\tablebox}
\setcounter{tocdepth}{3}
\tableofcontents

\section{Introduction to the characteristic method}

The characteristic methods, especially those within the semi-Lagrangian framework, have proved very successful in solving kinetic equations and other nonlocal problems \cite{CrouseillesLatuSonnendrucker2009,Kormann2015,KormannReuterRampp2019,XiongChenShao2016,GuoLiWang2018b,DimarcoLoubereNarskiRey2018}. In order to make the materials self-contained, we will briefly review their basic settings. 

\subsection{The Lawson integrators for partial integro-differential equations}

Consider the model problem
\begin{equation}
\frac{\partial}{\partial t} y(\bx, t) = \mathcal{L} y(\bx, t) + \mathcal{N} y(\bx, t),
\end{equation}
where $\mathcal{L}$ is the linear local operator and $\mathcal{N}$ is the nonlocal operator. Under the Lawson transformation $v(\bx, t) = \me^{(t_{n-1}-t)\mathcal{L}} y(\bx, t)$, it yields that
\begin{equation}
\frac{\partial}{\partial t} v(\bx, t) = \me^{(t_{n-1} - t)\mathcal{L}} \mathcal{N}( \me^{(t_{n-1} - t)\mathcal{L}} v(\bx, t)).
\end{equation}
Applying a $q$-step Adams method and transforming back to original variable yields the Lawson-Adams method,  \begin{equation}
y^n(\bx) = \me^{\tau\mathcal{L}} y^{n-1}(\bx)  + \sum_{k=0}^{q} \beta_k \me^{k\tau\mathcal{L}} \mathcal{N} y^{n-k}(\bx),
\end{equation} 
where $\tau = t_{n} - t_{n-1}$ is the time stepsize,
and $y^n(\bx)$ denotes the solution at $n$-th step. 

Specifically, for the partial integro-differential equation with a nonlocal operator $\Theta_V[f]$, e.g., the Boltzmann equation and the Wigner equation, of the form:
\begin{equation}\label{eq.Wigner}
\frac{\partial }{\partial t}f(\bm{x}, \bm{k}, t)+ \frac{\hbar \bm{k}}{m} \cdot \nabla_{\bm{x}} f(\bm{x},\bm{k}, t)  = \Theta_V[f](\bx, \bk, t).
\end{equation}
The commonly used Lawson schemes are collected as follows.
\begin{itemize}

\item[(1)] One-stage Lawson predictor-corrector scheme (LPC-1)
\begin{equation*}
\boxed{
\begin{split}
\textup{P}:\widetilde{f}^{n+1}(\bx, \bk) &= f^{n}(\mathcal{A}_\tau(\bx, \bk)) + \tau \Theta_{V}[f^{n}](\mathcal{A}_\tau(\bx, \bk)), \\
\textup{C}: f^{n+1}(\bx, \bk) &=  f^{n}(\mathcal{A}_\tau(\bx, \bk)) + \frac{\tau}{2} \Theta_V[\widetilde{f}^{n+1}](\bx, \bk) + \frac{\tau}{2} \Theta_{V}[f^{n}](\mathcal{A}_\tau(\bx, \bk)).
\end{split}
}
\end{equation*} 

\item[(2)] Two-stage Lawson-Adams predictor-corrector scheme (LAPC-2): 
\begin{equation*}
\boxed{
\begin{split}
\textup{P}: \widetilde{f}^{n+1}(\bx, \bk) =& f^{n}(\mathcal{A}_\tau(\bx, \bk)) + \frac{3\tau}{2} \Theta_{V}[f^{n}](\mathcal{A}_\tau(\bx, \bk)) \\
&- \frac{\tau}{2} \Theta_{V}[f^{n-1}](\mathcal{A}_{2\tau}(\bx, \bk)), \\
\textup{C}: f^{n+1}(\bx, \bk) = & f^{n}(\mathcal{A}_\tau(\bx, \bk)) + \frac{5\tau}{12} \Theta_V[\widetilde{f}^{n+1}](\bx, \bk) + \frac{8\tau}{12} \Theta_{V}[f^{n}](\mathcal{A}_\tau(\bx, \bk)) \\
&- \frac{\tau}{12} \Theta_{V}[f^{n-1}](\mathcal{A}_{2\tau}(\bx, \bk)).
\end{split}
}
\end{equation*} 

\item[(3)]
Three-stage Lawson-Adams predictor-corrector scheme (LAPC-3): 
\begin{equation*}
\boxed{
\begin{split}
\textup{P}: \widetilde{f}^{n+1}(\bx, \bk) = &f^{n}\mathcal{A}_\tau(\bx, \bk)) + \frac{23\tau}{12}  \Theta_V[f^{n}](\mathcal{A}_{\tau}(\bx, \bk)) \\
&- \frac{16\tau}{12}  \Theta_V[f^{n-1}](\mathcal{A}_{2\tau}(\bx, \bk)) + \frac{5\tau}{12}  \Theta_V[f^{n-2}](\mathcal{A}_{3\tau}(\bx, \bk)), \\
\textup{C}:f^{n+1}(\bx, \bk) = &f^{n}\mathcal{A}_\tau(\bx, \bk)) + \frac{9\tau}{24} \Theta_V[\widetilde{f}^{n+1}](\bx, \bk) + \frac{19\tau}{24}  \Theta_V[f^{n}](\mathcal{A}_{\tau}(\bx, \bk)) \\
&- \frac{5\tau}{24}  \Theta_V[f^{n-1}](\mathcal{A}_{2\tau}(\bx, \bk)) + \frac{\tau}{24}  \Theta_V[f^{n-2}](\mathcal{A}_{3\tau}(\bx, \bk)).
\end{split}
}
\end{equation*}
\end{itemize}
Here we use the notation $\mathcal{A}_\tau(\bx, \bk) = (\bx - \frac{\hbar \bk}{m} \tau, \bk)$. The Lawson scheme exploits the exact advection along the characteristic line, i.e.,  the semigroup $\me^{-\frac{\hbar \tau}{m} \bk \cdot \nabla_{\bx}} f(\bx, \bk, t) =  f(\mathcal{A}_\tau(\bx, \bk), t - \tau)$. The convergence order of $q$-stage Lawson predictor-corrector scheme is between $q$ and $q+1$ as it can be regarded as an implicit integrator with incomplete iteration. In practice, the one-step predictor-corrector scheme LPC-1 is used to obtain missing starting points for multistep schemes LAPC-2 and LAPC-3.

Apart from the non-splitting scheme, another commonly used scheme is the operator splitting (OS). Take the Strang splitting as an example.
\begin{equation*}
\boxed{
\begin{split}
&\textup{Half-step advection}: f^{n+1/2}(\bx, \bk) = f^{n}(\mathcal{A}_{\tau/2}(\bx, \bk)), \\
&\textup{Full-step of $\pdo$}:\widetilde{f}^{n+1/2}(\bx, \bk) = f^{n+1/2}(\bx, \bk) + \tau \Theta_V[f^{n+1/2}](\bx, \bk), \\
&\textup{Half-step advection}: f^{n+1}(\bx, \bk) = \widetilde{f}^{n+1/2}(\mathcal{A}_{\tau/2}(\bx, \bk)). \\
\end{split}
}
\end{equation*}
The Strang splitting adopted here is a first-order scheme overall as one of the subproblems is integrated by the backward Euler method.

\subsection{Cubic spline interpolation}

The standard way  to evaluate $ f^{n}(\mathcal{A}_{\tau}(\bx, \bk))$ and  $ \Theta_V[f^{n}](\mathcal{A}_{\tau}(\bx, \bk))$ on the shifted grid is realized by interpolation via a specified basis expansion of $f^n$. Typical choices include the spline wavelets \cite{CrouseillesLatuSonnendrucker2006,CrouseillesLatuSonnendrucker2009}, the Fourier pseudo-spectral basis and the Chebyshev polynomials \cite{ChenShaoCai2019}.  Regarding the fact that the spatial advection is essentially local, we only consider the cubic B-spline as it is a local wavelet basis with low numerical dissipation, and the cost scales as $\mathcal{O}(N_x^d)$ with $d$ the dimensionality  \cite{CrouseillesLatuSonnendrucker2009}.

Now we focus on the unidimensional uniform setting as the multidimensional spline can be constructed by its tensor product. Suppose the computational domain is $[x_0, x_{N}]$ containing  $N+1$ grid points with uniform spacing $h = \frac{x_{N} - x_0}{N}$. The projection of $\varphi(x)$ onto the cubic spline basis is given by
\begin{equation}\label{interpolation}
\varphi(x) \approx s(x) = \sum_{\nu = -1}^{N+1}  \eta_{\nu} B_{\nu}(x) \quad \textup{subject to} \quad \varphi(x_i) = s(x_i), \quad i = 0, \dots, N.
\end{equation}
$B_\nu$ is the cubic B-spline with compact support over four grid points, 
\begin{equation}
B_{\nu}(x) = 
\left\{
\begin{split}
&\frac{(x - x_{\nu-2})^3}{6h^3}, \quad  x \in [x_{\nu-2}, x_{\nu-1}],\\
&-\frac{(x - x_{\nu-1})^3}{2h^3} + \frac{(x - x_{\nu-1})^2}{2h^2} + \frac{(x - x_{\nu-1})}{2h} + \frac{1}{6}, \quad x \in [x_{\nu-1}, x_{\nu}],\\
&-\frac{(x_{\nu+1} - x)^3}{2h^3} +\frac{(x_{\nu+1} - x)^2}{2h^2} + \frac{(x_{\nu+1} - x)}{2h} + \frac{1}{6}, \quad x \in [x_{\nu}, x_{\nu+1}],\\
&\frac{(x_{\nu+2} - x)^3}{6h^3}, \quad x \in [x_{\nu+1}, x_{\nu+2}],\\
&0, \quad \textup{otherwise},
\end{split}
\right.
\end{equation}
implying that $B_{\nu - 1}, B_{\nu}, B_{\nu+1}, B_{\nu+2}$ overlap a grid interval $(x_{\nu}, x_{\nu+1})$ \cite{MalevskyThomas1997}.

Now it requires to solve the $N+3$ coefficients $\bm{\eta} = (\eta_{-1}, \dots, \eta_{N+1})$. Since only $B_{i \pm 1}(x_i) = \frac{1}{6}$ and $B_i(x_i) = \frac{2}{3}$, substituting it into Eq.~\eqref{interpolation} yields $N+1$ equations for $N+3$ variables,
\begin{equation}\label{three_term_relation}
\varphi(x_i) = \frac{1}{6} \eta_{i-1} + \frac{2}{3} \eta_{i} + \frac{1}{6} \eta_{i+1}, \quad 0 \le i \le N.
\end{equation}

Two additional equations are needed to determine the unique solution of $\bm{\eta}$, which are given by a specified boundary condition at both ends of the interval.  For instance, consider the Hermite boundary condition (also termed the clamped spline) \cite{CrouseillesLatuSonnendrucker2006},
\begin{equation}
s^{\prime}(x_0) = \phi_L, \quad s^{\prime}(x_{N}) = \phi_R,
\end{equation}
where $\phi_L$ and $\phi_R$ are parameters to be determined.  In particular, when $\phi_L = \phi_R = 0$, it reduces to the Neumann boundary condition on both sides. Since 
\begin{equation}
s^{\prime}(x_i) = - \frac{1}{2h} \eta_{i-1} +  \frac{1}{2h} \eta_{i+1}, \quad i = 0, \dots, N,
\end{equation}
it is equivalent to add two constraints,
\begin{equation}
\phi_L =  -\frac{1}{2h} \eta_{-1} + \frac{1}{2h} \eta_1,\quad \phi_R = -\frac{1}{2h} \eta_{N-1} + \frac{1}{2h} \eta_{N+1}.
\end{equation}
%
%
%
Thus all the coefficients can be obtained straightforwardly by solving the equation 
\begin{equation}\label{cubic_spline}
A (\eta_{-1}, \dots, \eta_{N+1})^T = (\phi_L, \varphi(x_0), \dots, \varphi(x_{N}), \phi_R)^T.
\end{equation}
Note that $(N+3)\times (N+3)$ coefficient matrix $A$ has an explicit LU decomposition,
\begin{equation}\label{coefficient_matrix}
A = \frac{1}{6}
\begin{pmatrix}
-3/h & 0 & 3/h & 0 & \cdots & 0 \\
1     & 4 & 1     & 0 &           & \vdots \\
0     & 1 & 4     & 1 &  & \vdots \\
\vdots & \vdots & \vdots & \vdots & \ddots & \vdots \\ 
\vdots &           & 0 & 1 & 4 & 1 \\
0         &  0       & 0 & -3/h & 0 & 3/h \\
\end{pmatrix}.
\end{equation}
where
\begin{equation}\label{explicit_L}
L = 
\begin{pmatrix}
1 & 0 & 0 &\cdots & \cdots & 0 \\
-h/3 & 1 & 0  &\ddots & & \vdots \\
0 & l_1 & 1  &\ddots & & \vdots \\
0 & 0    & l_2  & \ddots & & \vdots \\
\vdots & \vdots   & & l_{N} & 1 & 0 \\
0 & 0 & \cdots  & -\frac{3l_{N}}{h} & \frac{3l_{N+1}}{h} & 1 \\
\end{pmatrix}
\end{equation}
and
\begin{equation}\label{explicit_U}
U = \frac{1}{6}
\begin{pmatrix}
-3/h & 0 & 3/h & 0 &\cdots & \cdots & 0 \\
0 & d_1 & 2 & 0 &\ddots & & \vdots \\
0 & 0 & d_2 & 1 &\ddots & & \vdots \\
0 & 0    & 0  & d_3 & \ddots & & \vdots \\
\vdots & \vdots  & & & 0 & d_{N+1} & 0 \\
0 & 0 & \cdots  & & 0 & 0 & 3 d_{N+2}/h \\
\end{pmatrix},
\end{equation}
with
\begin{equation}
\begin{split}
& d_1 = 4, \quad l_1 = 1/4, \quad d_2 = 4 - 2l_1 = 7/2, \\
& l_i = 1/d_i, \quad d_{i+1} = 4 - l_i, \quad i = 2, \dots, N+1, \\
&  l_{N+1} = 1/(d_{N} d_{N+1}), \quad d_{N+2} = 1 - l_{N+1}.
\end{split}
\end{equation}

The above scheme can achieve fourth order convergence in spatial spacing $h$ and conserves the total mass. Besides, the time step in the semi-Lagrangian method is usually not restricted by the CFL condition, that is, $C = \hbar \max |k| \tau /h > 1 $ is allowed.

\section{Parallel characteristic method}

For a  6-D problem, the foremost problem is the storage of huge 6-D tensors as the memory to store a $101^3 \times 64^3$ grid is $1.08$TB in single precision, which is still prohibitive for modern computers. 

Fortunately, the characteristic method can be realized in a distributed manner as pointed out in several pioneering works \cite{MalevskyThomas1997,CrouseillesLatuSonnendrucker2009}. Without loss of generality, we divide $N+1$ grid points on a line into $p$ uniform parts, with $M = N/p$,
\begin{align*}
\underbracket{x_0 < x_1 < \cdots < x_{M-1}}_{\textup{the 1st processor}} < \underbracket{x_{M}}_{\textup{shared}} < \cdots < \underbracket{ x_{(p-1)M}}_{\textup{shared}} < \underbracket{x_{(p-1)M+1} < \cdots < x_{pM}}_{\textup{$p$-th processor}},
\end{align*}
where the $l$-th processor only manipulates $\mathcal{X}_l$ with $\mathcal{X}_l = (x_{(l-1)M+1}, \dots, x_{l M})$, $l = 1, \dots, p$.
The grid points $x_{M}, x_{2M}, \dots, x_{(p-1)M}$ are shared by the adjacent patches. Our target is to make 
\begin{equation}
\bm{\eta}^{(l)}= (\eta_{-1}^{(l)}, \dots, \eta_{M+1}^{(l)}) \approx (\eta_{-1 +(l-1)M}, \dots, \eta_{(l-1)M+M+1}), \quad l = 1, \dots, p,
\end{equation}
say, the local spline coefficients $\bm{\eta}^{(l)}$ for $l$-th piece should approximate to those in global B-spline as accurately as possible.

\subsection{Effective Hermite boundary condition based on finite difference stencils}

In order to solve  $\bm{\eta}^{(l)}$ efficiently, Crouseilles, Latu and Sonnendr{\"u}cker suggested to impose an effective Hermite boundary condition on the shared grid points (CLS-HBC for short) \cite{CrouseillesLatuSonnendrucker2006,CrouseillesLatuSonnendrucker2009} . 
\begin{equation}\label{Hermite_boundary_condition}
 \varphi^{\prime}(x_{lM}) =  s^{\prime}(x_{lM}), \quad l = 1, \dots, p,
\end{equation}
so that it needs to solve 
\begin{align}
\varphi^{\prime}(x_{lM}) &= -\frac{1}{2h} \eta_{M-1}^{(l)} + \frac{1}{2h} \eta_{M+1}^{(l)} = -\frac{1}{2h} \eta_{-1}^{(l+1)} + \frac{1}{2h} \eta_{1}^{(l+1)}.
\end{align}

The problem is that the derivates $\varphi^{\prime}(x_{lM})$ on the adjacent points are actually unknown, so that they have to be interpolated by a finite difference stencil. The authors suggest to use the recursive relation from the spline transform matrix \eqref{coefficient_matrix} and three-term relation $\varphi(x_i) = \frac{1}{6} \eta_{i-1} + \frac{2}{3} \eta_{i} + \frac{1}{6} \eta_{i+1}$, $0 \le i \le N$.

Following \cite{CrouseillesLatuSonnendrucker2006} and taking $i = lM$, it starts from
\begin{equation}
\begin{split}
s^{\prime}(x_{i}) = & -\frac{1}{2h}\eta_{i-1} + \frac{1}{2h}\eta_{i+1} \\
= & -\frac{1}{2h} \left( \frac{3}{2} \varphi(x_{i-1}) - \frac{1}{4} \eta_{i-2} - \frac{1}{4} \eta_{i}\right) \\
& + \frac{1}{2h} \left( \frac{3}{2} \varphi(x_{i+1}) - \frac{1}{4} \eta_{i} - \frac{1}{4} \eta_{i+2}\right) \\
= & \frac{3}{4h}(\varphi(x_{i+1}) - \varphi(x_{i-1})) + \frac{1}{8h}(\eta_{i- 2} - \eta_{i+2}),
\end{split}
\end{equation}
so that it arrives at the recursive relation,
\begin{equation}\label{recursive_relation}
s^{\prime}(x_i) = \frac{3}{4h}(\varphi(x_{i-1}) + \varphi(x_{i+1})) - \frac{1}{4}(s^{\prime}(x_{i - 1}) - s^{\prime}(x_{i+1})),
\end{equation}

By further replacing $s^{\prime}(x_{i-1})$ and $s^{\prime}(x_{i+1})$ by Eq.~\eqref{recursive_relation}, it arrives at a longer expansion
\begin{equation}
\begin{split}
s^{\prime}(x_i) = & \frac{6}{7h}(\varphi(x_{i+1}) + \varphi(x_{i-1})) - \frac{3}{14h}(\varphi(x_{i+2}) + \varphi(x_{i-2})) \\
&+ \frac{1}{14}(s^{\prime}(x_{i+2}) - s^{\prime}(x_{i-2})).
\end{split}
\end{equation}

It obtains the final expansion with $\alpha = 1 - 2/14^2$,
\begin{equation}
\left(\alpha - \frac{2}{\alpha 14^2}\right) s^{\prime}(x_i) = \sum_{j = -8}^{8} \omega_j \varphi(x_{i+j}) + \frac{1}{\alpha 14^2} ( s^{\prime}(x_{i+8})  +  s^{\prime}(x_{i-8}) ),
\end{equation}
associated with the fourth-order finite difference approximation
\begin{equation}
s^{\prime}(x_{i+8}) \approx \frac{-\varphi(x_{i+10}) + 8\varphi(x_{i+9}) - 8 \varphi(x_{i+7}) + \varphi(x_{i+6})}{12h}.
\end{equation}
To sum up, it arrives at the formula
\begin{equation}
s^{\prime}(x_i) = \sum_{j = -10}^{-1} \tilde{\omega}_j^- \varphi(x_{i+j}) + \sum_{j=1}^{10}\tilde{\omega}_j^+ \varphi(x_{i+j}),
\end{equation}
where the coefficients $\tilde{\omega}_j^-$ are collected in Table \ref{Coefficients} and $\tilde{\omega}_j^+ = -\tilde{\omega}_j^-$. 
\begin{table}[!h]
  \centering
  \caption{\small Coefficients for the approximation of the derivatives \cite{CrouseillesLatuSonnendrucker2009}.
}\label{Coefficients}
 \begin{lrbox}{\tablebox}
 \begin{tabular}{cccccccc}
 \hline\hline
 $j$ & $-10$ &  $-9$  & $-8$ &  $-7$ & $-6$ &  \\
 \hline
$\tilde{\omega}_j$	&	$0.2214309755$E-5	&	$-1.771447804$E-5	&	$7.971515119$E-5	&	$-3.011461267$E-4	&	$1.113797807$E-3	\\	
\hline
 $j$ & $-5$ &  $-4$  & $-3$ &  $-2$ & $-1$ &  \\
 \hline
$\tilde{\omega}_j$	&	$-4.145187862$E-3	&	$0.01546473933$	&	$-0.05771376946$	&	$0.2153903385$	&	$-0.8038475846$	\\	
\hline\hline
 \end{tabular}
\end{lrbox}
\scalebox{0.80}{\usebox{\tablebox}}
\end{table}

At each step, $\sum_{j = -10}^{-1} \tilde{\omega}_j^- \varphi(x_{i+j})$ and $ \sum_{j=1}^{10}\tilde{\omega}_j^+ \varphi(x_{i+j})$ can be assembled by left and right processor independently, and then  data is exchanged only in adjacent processors to merge the effective boundary condition. The remaining task to solve algebraic equations in each processor independently 
\begin{equation}
A^{(l)} \bm{\eta}^{(l)} = (\phi_L^{(l)}, \varphi(x_{(l-1)M}), \dots, \varphi(x_{lM}), \phi_R^{(l)})^T,
\end{equation}
where $A^{(l)}$ is a $(M+3)\times(M+3)$ matrix with the form like Eq.~\eqref{coefficient_matrix}, and
\begin{equation}
\phi_{R}^{(l)} =  \phi_{L}^{(l+1)} \approx \sum_{j = -10}^{-1} \tilde{\omega}_j^- \varphi(x_{lM+j}) + \sum_{j=1}^{10}\tilde{\omega}_j^+ \varphi(x_{lM+j}).
\end{equation}

\subsection{Perfectly  matched boundary condition (PMBC)}

We suggest to adopt  another effective boundary condition, termed the perfectly matched boundary condition (PMBC), based on a key observation made in \cite{MalevskyThomas1997}. Specifically, we start from the exact solution of $\bm{\eta}$, with $(b_{ij}) = A^{-1}$ of size $(N+3)\times (N+3)$.

For the sake of convenience, the subindex of $(b_{ij})$ starts from $-1$ and ends at $N+1 = pM+1$. The solution of Eq.~\eqref{cubic_spline} reads that
\begin{equation}
\eta_i = b_{ii} \varphi(x_i) + \sum_{j=-1}^{i-1} b_{ij} \varphi(x_j) + \sum_{j = i+1}^{pM+1} b_{i j} \varphi(x_{j}), \quad i = -1, \dots, pM+1.
\end{equation}
We can make a truncation for $|i-j| \ge n_{nb}$ as the off-diagonal elements exhibit exponential decay away from the main diagonal,
\begin{equation}\label{truncation}
\eta_i \approx b_{ii} \varphi(x_i) + \sum_{j= i - n_{nb}+1}^{i-1} b_{ij} \varphi(x_j) + \sum_{j = i+1}^{i + n_{nb}-1} b_{i j} \varphi(x_{j}), \quad i = -1, \dots, pM+1.
\end{equation} 
Using the truncated stencils \eqref{truncation},
\begin{equation*}
\begin{split}
& \eta_{lM-1} \approx \sum_{j=(lM-1)-n_{nb}+1}^{(lM-1)+n_{nb}-1} b_{lM-1, j} \varphi(x_j) =  \sum_{j=-n_{nb}}^{n_{nb}-2} b_{lM-1, lM+j} \varphi(x_{lM+j}), \\
& \eta_{lM+1} \approx \sum_{j=(lM+1)-n_{nb}+1}^{(lM+1)+n_{nb}-1} b_{lM+1, j} \varphi(x_j) = \sum_{j=-n_{nb}+2}^{n_{nb}} b_{lM+1, lM+ j} \varphi(x_{lM+j}). \\ 
\end{split}
\end{equation*}
By further adding four more terms to complete the summations from $-n_{nb}$ to $n_{nb}$, it yields that 
\begin{equation*}
\begin{split}
-\frac{1}{2h} \eta_{lM-1} + \frac{1}{2h} \eta_{lM+1} \approx &\sum_{j=-n_{nb}}^{n_{nb}} \left(-\frac{1}{2h} b_{lM-1, lM+j} + \frac{1}{2h} b_{lM+1, lM+j}\right) \varphi(x_{lM+j})\\ 
= &\underbracket{\sum_{j=-n_{nb}}^{-1} \left(-\frac{1}{2h} b_{lM-1, lM+j} + \frac{1}{2h} b_{lM+1, lM+j}\right) \varphi(x_{lM+j})}_{\textup{stored in left processor}} \\
&~+\underbracket{\sum_{j=1}^{n_{nb}} \left(-\frac{1}{2h} b_{lM-1, lM+j} + \frac{1}{2h} b_{lM+1, lM+j}\right) \varphi(x_{lM+j})}_{\textup{stored in right processor}}  \\
&~+ \underbracket{\left(-\frac{1}{2h} b_{lM-1, lM} + \frac{1}{2h} b_{lM+1, lM}\right) \varphi(x_{lM}).}_{\textup{shared by adjacent two processors}} 
\end{split}
\end{equation*}
Thus it arrives at the formulation of PMBC
\begin{equation*}
\begin{split}
\phi_{R}^{(l)} =  \phi_{L}^{(l+1)} \approx & \underbracket{\frac{1}{2}c_{0,l} \varphi(x_{lM}) +  \sum_{j = 1}^{n_{nb}} c_{j, l}^{-} \varphi(x_{lM-j})}_{\textup{stored in left processor}} + \underbracket{\frac{1}{2}c_{0,l} \varphi(x_{lM}) + \sum_{j = 1}^{n_{nb}} c_{j, l}^{+} \varphi(x_{lM+j})}_{\textup{stored in right processor}},
\end{split}
\end{equation*}
where $c_{0, l} = -\frac{b_{lM-1, lM}}{2h} + \frac{b_{lM+1, lM}}{2h}$ and 
\begin{equation}\label{PMBC_coeffcients}
\begin{split}
&c_{j, l}^+ =  -\frac{b_{lM-1, lM+j}}{2h} + \frac{b_{lM+1, lM+j}}{2h}, \quad c_{j, l}^- =  -\frac{b_{lM-1, lM-j}}{2h} + \frac{b_{lM+1, lM-j}}{2h}.
\end{split}
\end{equation}

It deserves to mention that other spline boundary conditions can also be represented by PMBC, following the same idea in Eq.~\eqref{truncation}. When the natural boundary conditions are adopted,
\begin{equation}\label{natural_boundary}
 \frac{1}{h^2} \eta_{-1} - \frac{2}{h^2} \eta_0 +  \frac{1}{h^2} \eta_{1} = 0, \quad  \frac{1}{h^2} \eta_{N-1} - \frac{2}{h^2} \eta_N +  \frac{1}{h^2} \eta_{N+1} = 0,
\end{equation}
the coefficient matrix is
\begin{equation}\label{coefficient_matrix_2}
\widetilde A = \frac{1}{6}
\begin{pmatrix}
1/h^2 & -2/h^2 & 1/h^2 & 0 & \cdots & 0 \\
1     & 4 & 1     & 0 &           & \vdots \\
0     & 1 & 4     & 1 &  & \vdots \\
\vdots & \vdots & \vdots & \vdots & \ddots & \vdots \\ 
\vdots &           & 0 & 1 & 4 & 1 \\
0         &  0       & 0 & 1/h^2 & -2/h^2 & 1/h^2 \\
\end{pmatrix}.
\end{equation}

Denote by $(\widetilde{b}_{ij}) = \widetilde{A}^{-1}, -1\le i, j \le N+1$. Equivalently,  the equations $\widetilde A\bm{\eta}^T= (0, \varphi(x_0), \dots, \varphi(x_N),  0)^T$ can be  cast into $A\bm{\eta}^T= (\phi_{L}^{(1)}, \varphi(x_0), \dots, \varphi(x_N),  \phi_{R}^{(p)})^T$ since
\begin{equation}
\begin{split}
& \eta_{-1} \approx \sum_{j = -1}^{n_{nb} - 2} \widetilde b_{-1, j} \varphi(x_{j}), \quad \eta_{1} \approx  \sum_{j= -1 }^{n_{nb}} \widetilde b_{1, j} \varphi(x_j).
\end{split}
\end{equation} 
By adding two terms and noting that $\varphi(x_{-1}) = 0$, it yields that
\begin{equation}
\phi_{L}^{(1)} =\frac{\eta_1-\eta_{-1}}{2h} \approx \sum_{j=0}^{n_{nb}} c_{j, 0}^{-} \varphi(x_j), \quad c_{j, 0}^{-} =  \frac{1}{2h}(-\widetilde{b}_{-1, j}  + \widetilde{b}_{1, j}).
\end{equation}
Similarly, for the other end, noting that $\varphi(x_{N+1}) = 0$,
\begin{equation}
\begin{split}
& \eta_{N-1} \approx \sum_{j = -1}^{n_{nb} - 2} \widetilde b_{N-1, N-j} \varphi(x_{N-j}), \quad \eta_{N+1} \approx  \sum_{j= -1 }^{n_{nb}} \widetilde b_{N+1, N-j} \varphi(x_{N-j}),
\end{split}
\end{equation} 
so that
\begin{equation}
\phi_{R}^{(p)} = \frac{\eta_{N+1}-\eta_{N-1}}{2h}  \approx \sum_{j=0}^{n_{nb}} c_{j, p}^{+}\varphi(x_{N-j}), \quad c_{j, p}^{+} = \frac{1}{2h}(-\widetilde{b}_{N-1, N-j}  + \widetilde{b}_{N+1, N-j}).
\end{equation}

\subsection{Comparison between two effective Hermite boundary conditions}

It is shown that PMBC is more preferable than CLS-HBC in consideration of numerical accuracy. 
\begin{example}[1-D spline]
The test problem is 
\begin{equation}
\varphi(x) = \sin(x), \quad x \in [0, 8],
\end{equation}
subject to 
\begin{equation}
\varphi^{\prime}(0) = 0, \quad \varphi^{\prime}(8) = 0.
\end{equation}
\end{example}

\begin{figure}[!!h]
\centering
\includegraphics[width=0.48\textwidth,height=0.27\textwidth]{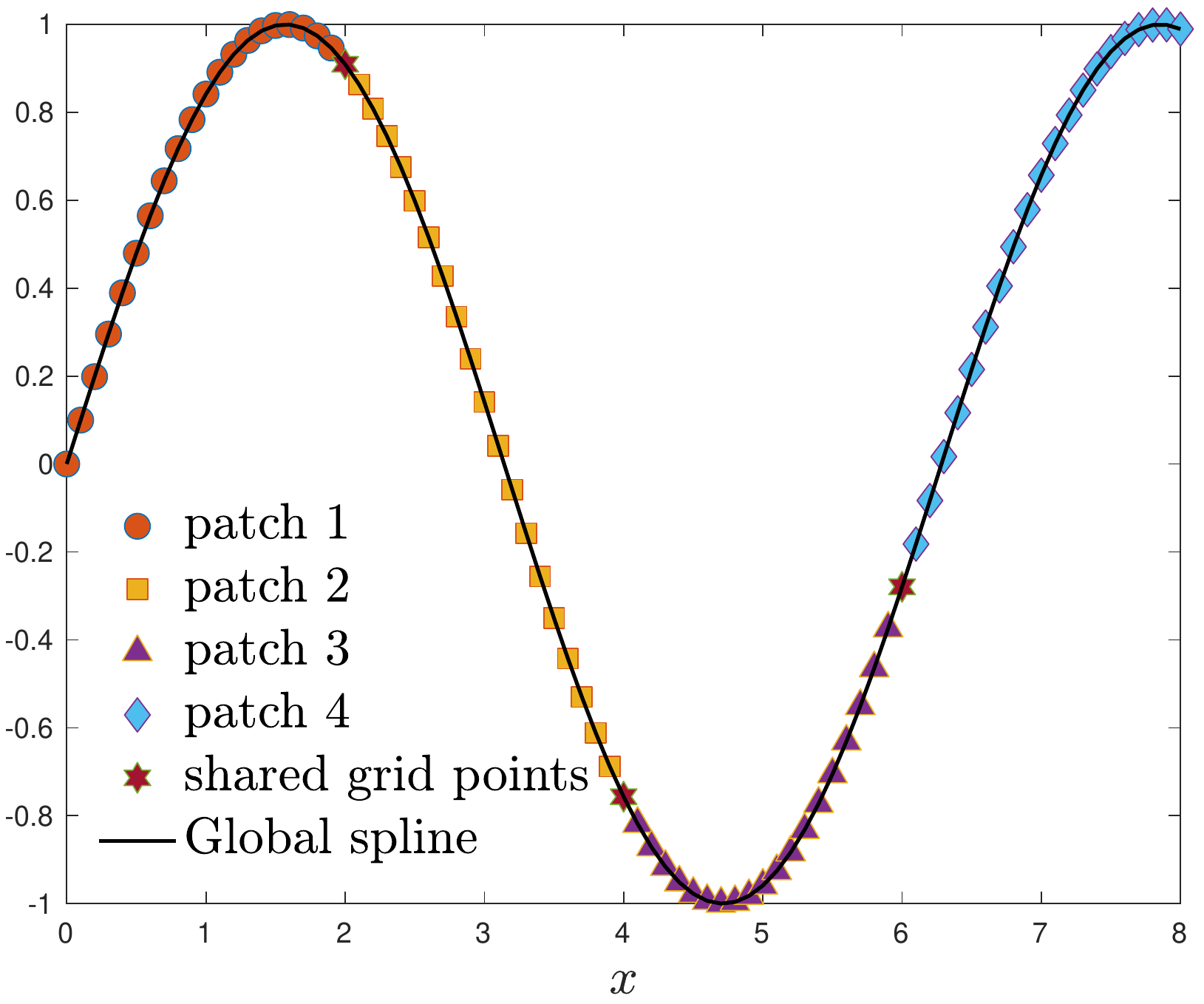}
\caption{\small An illustration of cubic B-spline on four patches. The grid points $x_{lM}$ $(l = 1, \dots, p) $ are shared by adjacent processors. In principle,  the splines on patches should approximate to the global one as accurately as possible. \label{plot_patch_spline}}
\end{figure}

For parallel implementation, the spline is decomposed into $4$ patches as given in Figure \ref{plot_patch_spline}, and each patch contains $(N-1)/4+1$ grid points .

\begin{figure}[!h]
\centering
\subfigure[$N = 81$.]{
{\includegraphics[width=0.48\textwidth,height=0.27\textwidth]{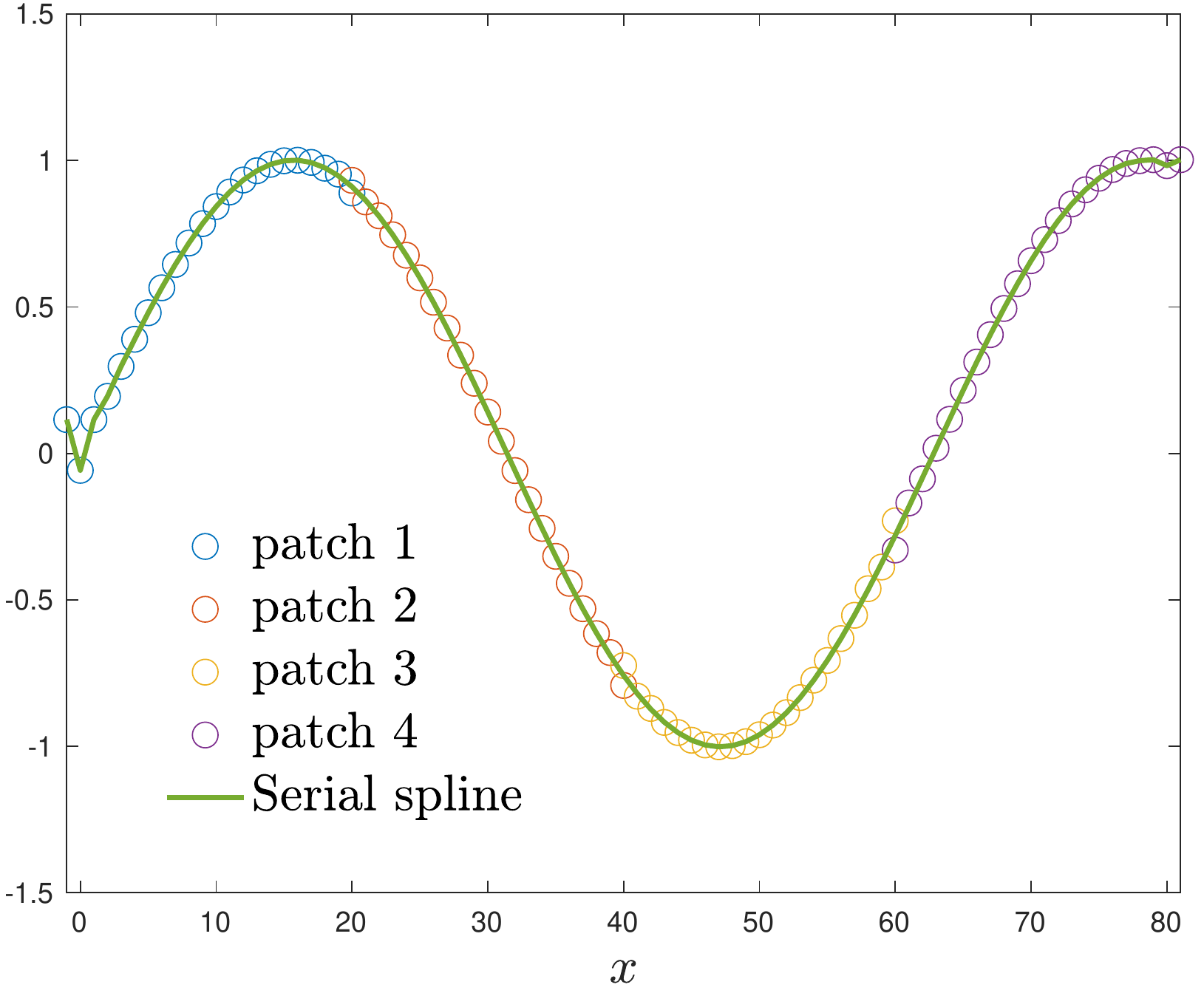}}
{\includegraphics[width=0.48\textwidth,height=0.27\textwidth]{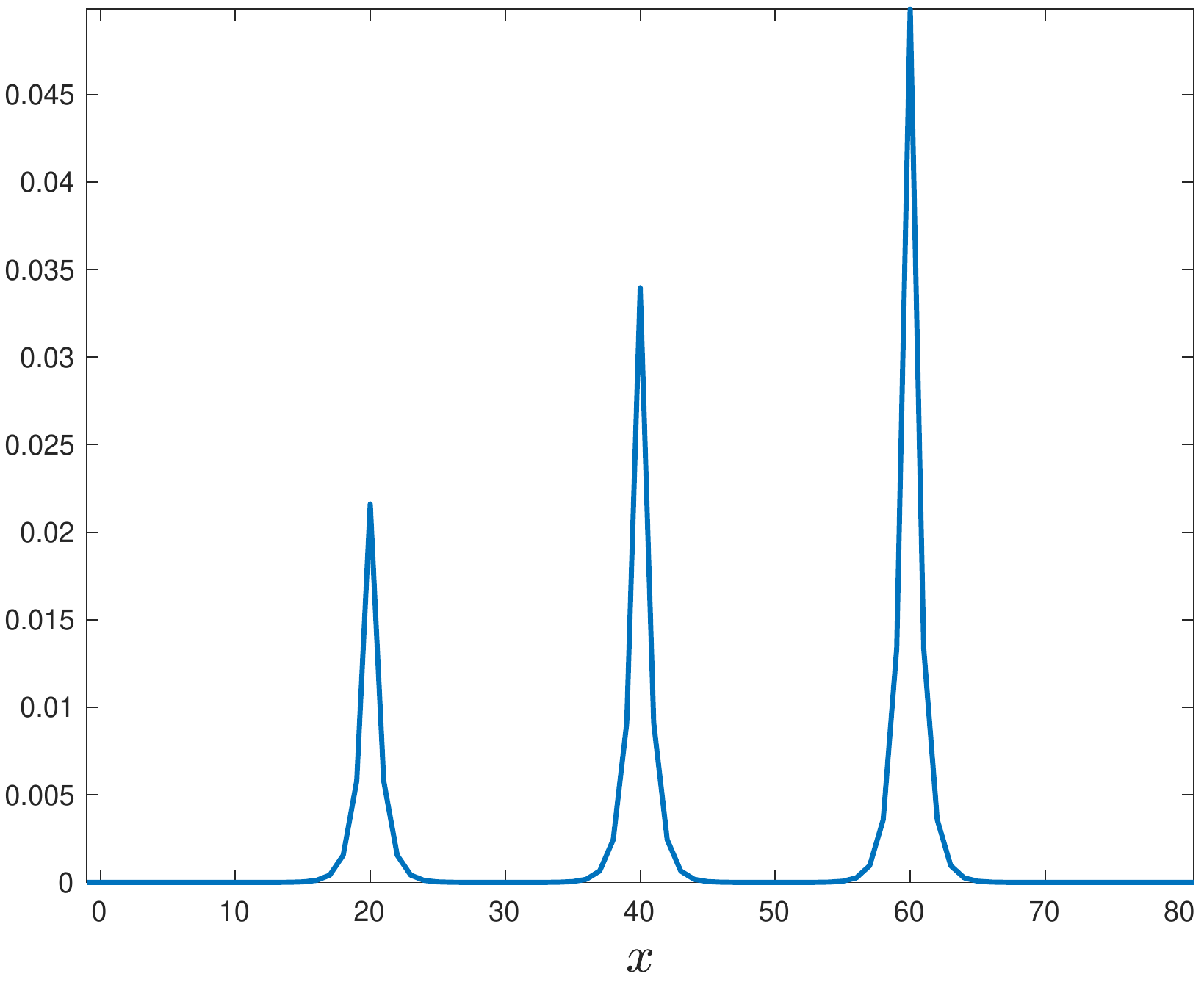}}}
\\
\centering
\subfigure[$N = 161$.]{
{\includegraphics[width=0.48\textwidth,height=0.27\textwidth]{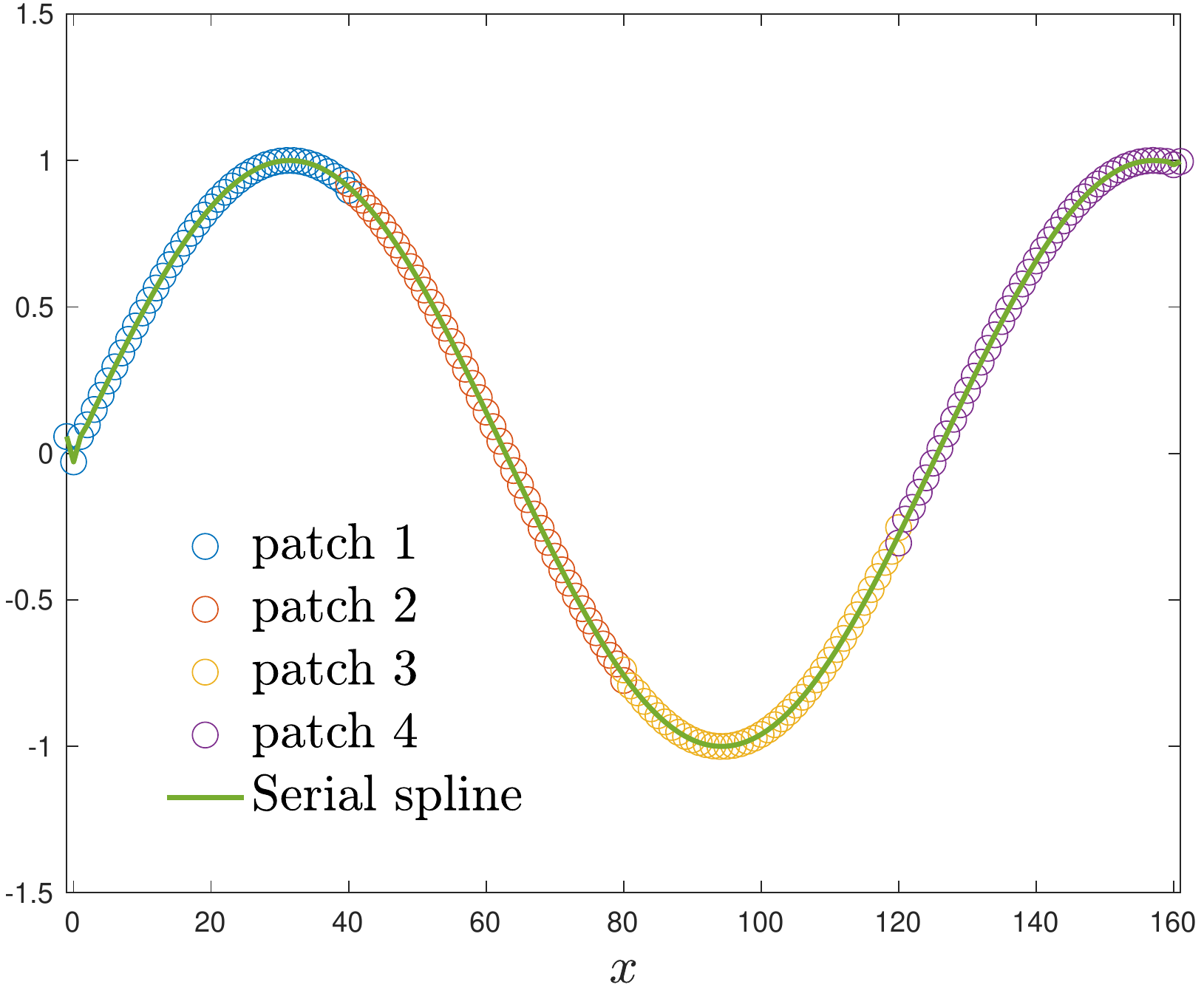}}
{\includegraphics[width=0.48\textwidth,height=0.27\textwidth]{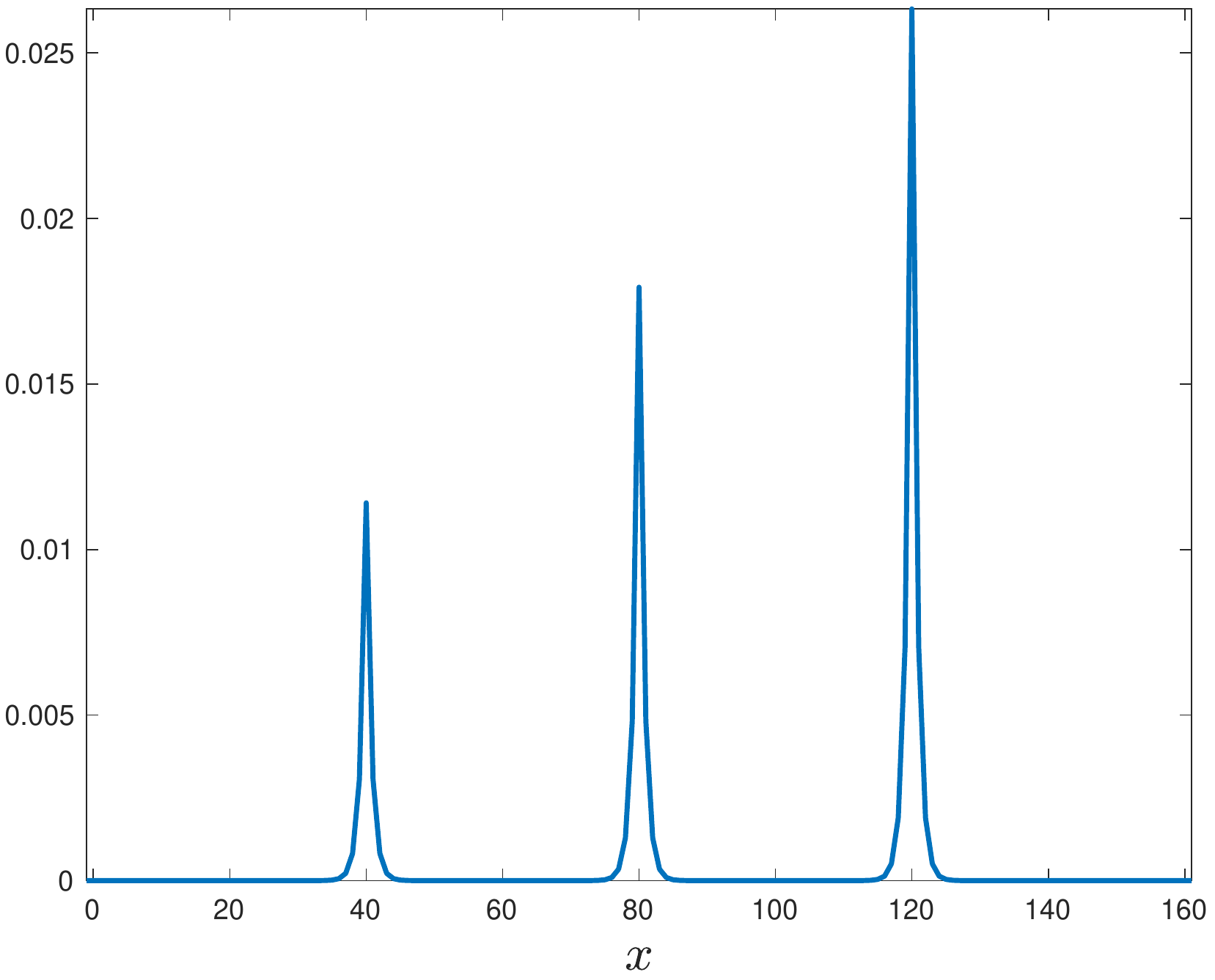}}}
\\
\centering
\subfigure[$N = 321$.]{
{\includegraphics[width=0.48\textwidth,height=0.27\textwidth]{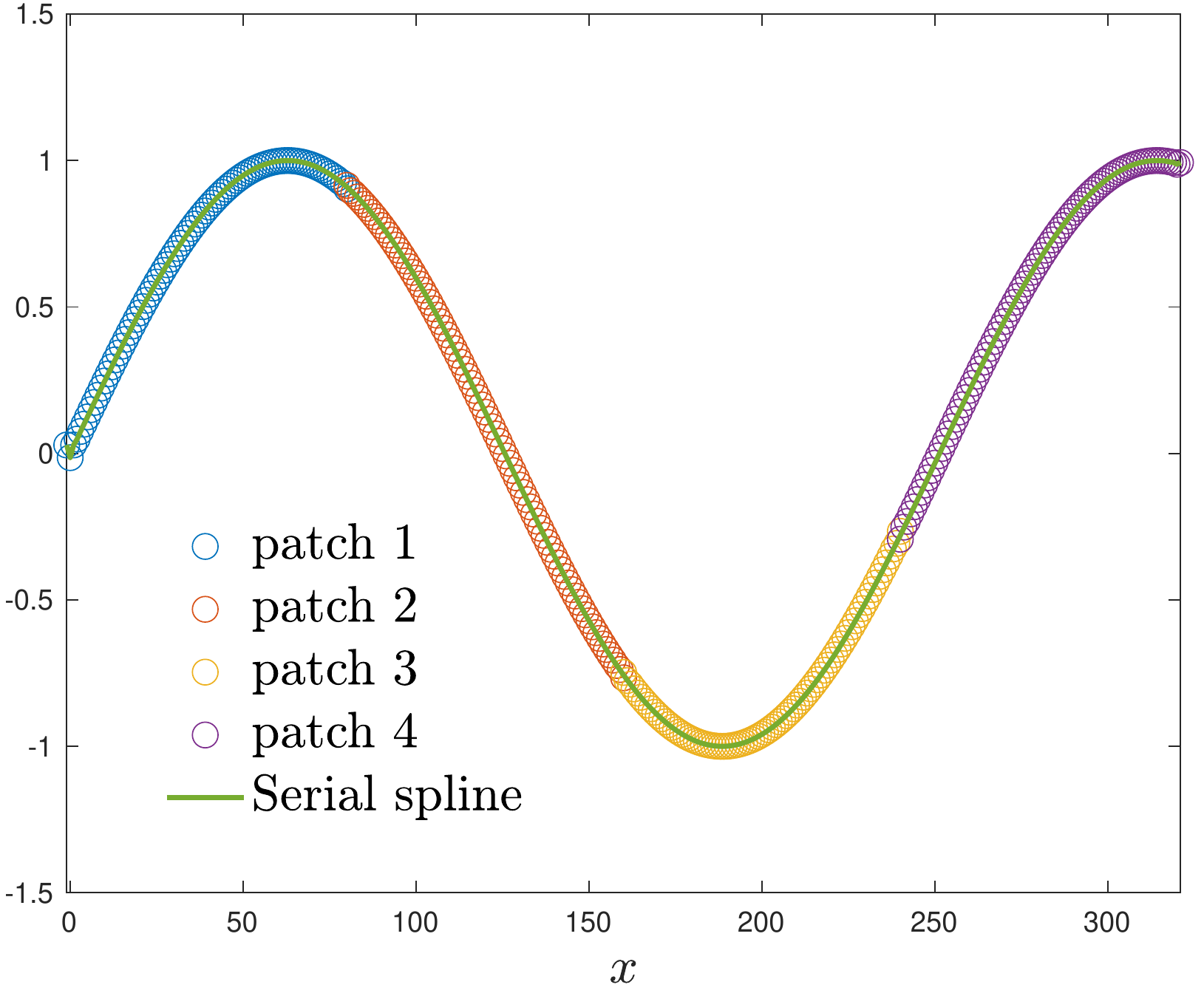}}
{\includegraphics[width=0.48\textwidth,height=0.27\textwidth]{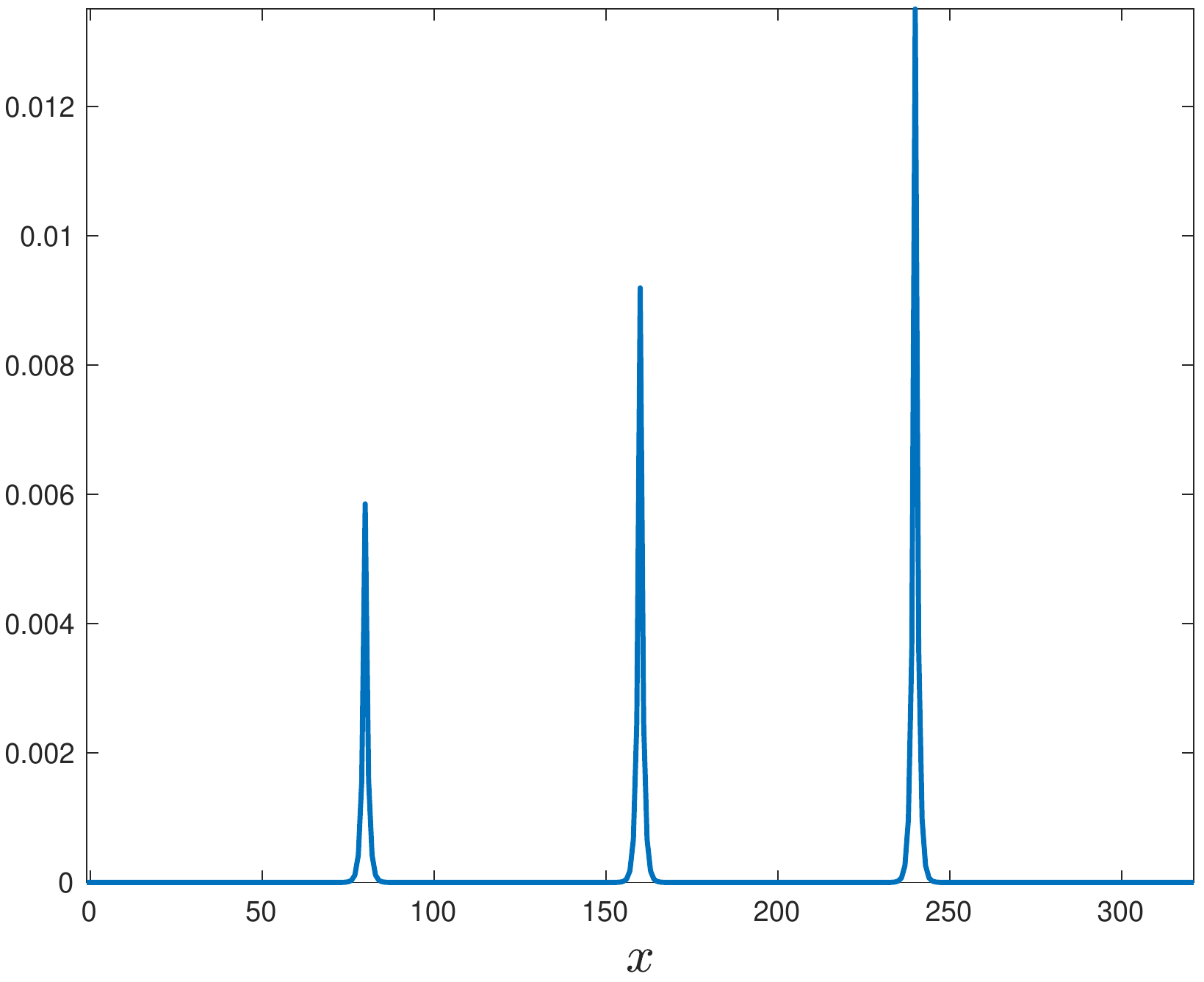}}}
\caption{\small Spline coefficients (left) and absolute errors (right) under CLS-HBC. Large errors are observed at the junction points.   \label{spline_CLS_HBC}}
\end{figure}

First, we adopt CLS-HBC and the results are shown in Figure \ref{spline_CLS_HBC}. It is observed that the errors of are concentrated at the junction points of adjacent patches. Indeed, the accuracy is improved under more collocation points (or equivalently, using smaller step size). The relative errors are less than $5\%$ when $N = 81$.

\begin{figure}[!h]
\centering
\subfigure[$n_{nb} = 5$.]{
{\includegraphics[width=0.48\textwidth,height=0.27\textwidth]{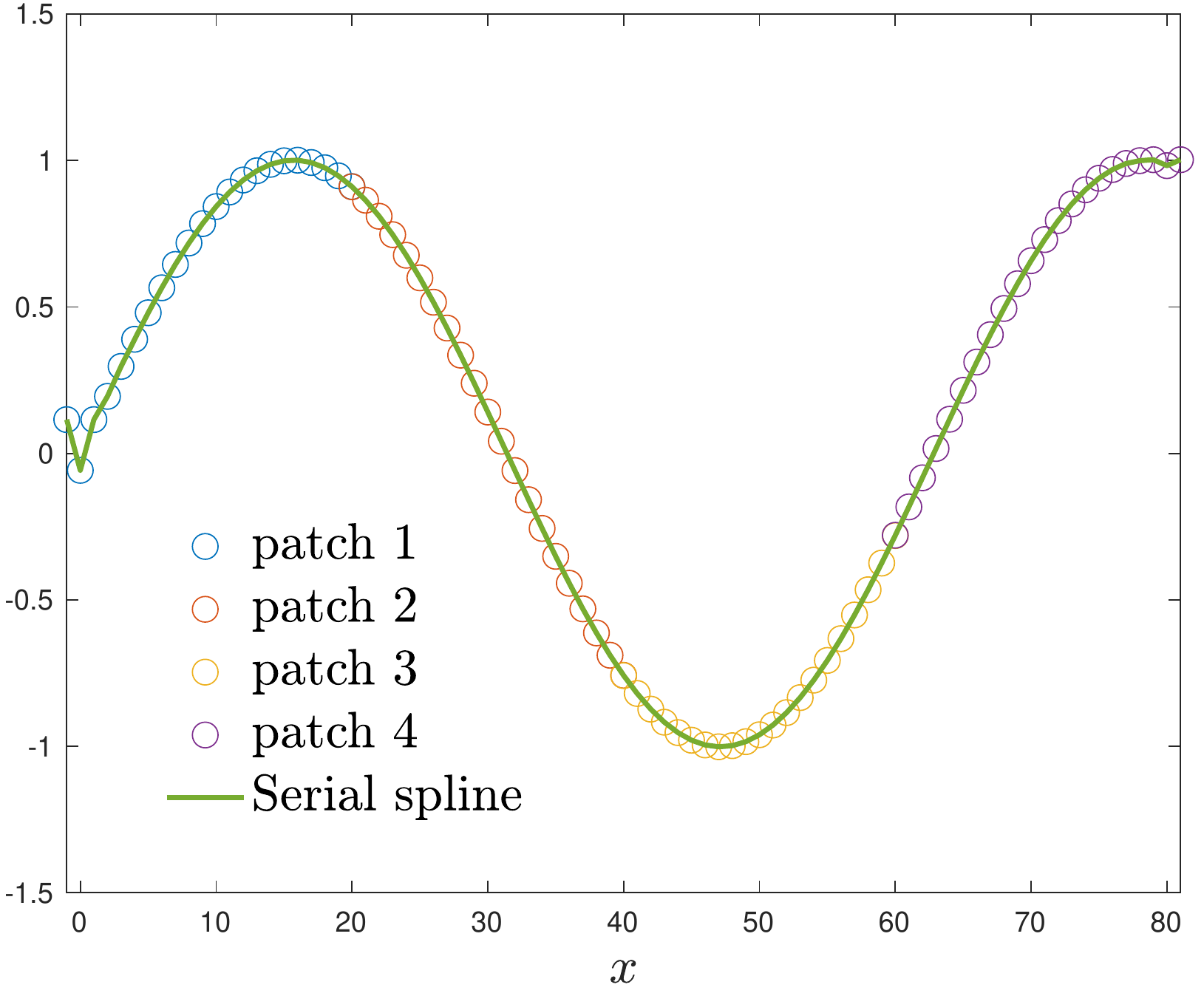}}
{\includegraphics[width=0.48\textwidth,height=0.27\textwidth]{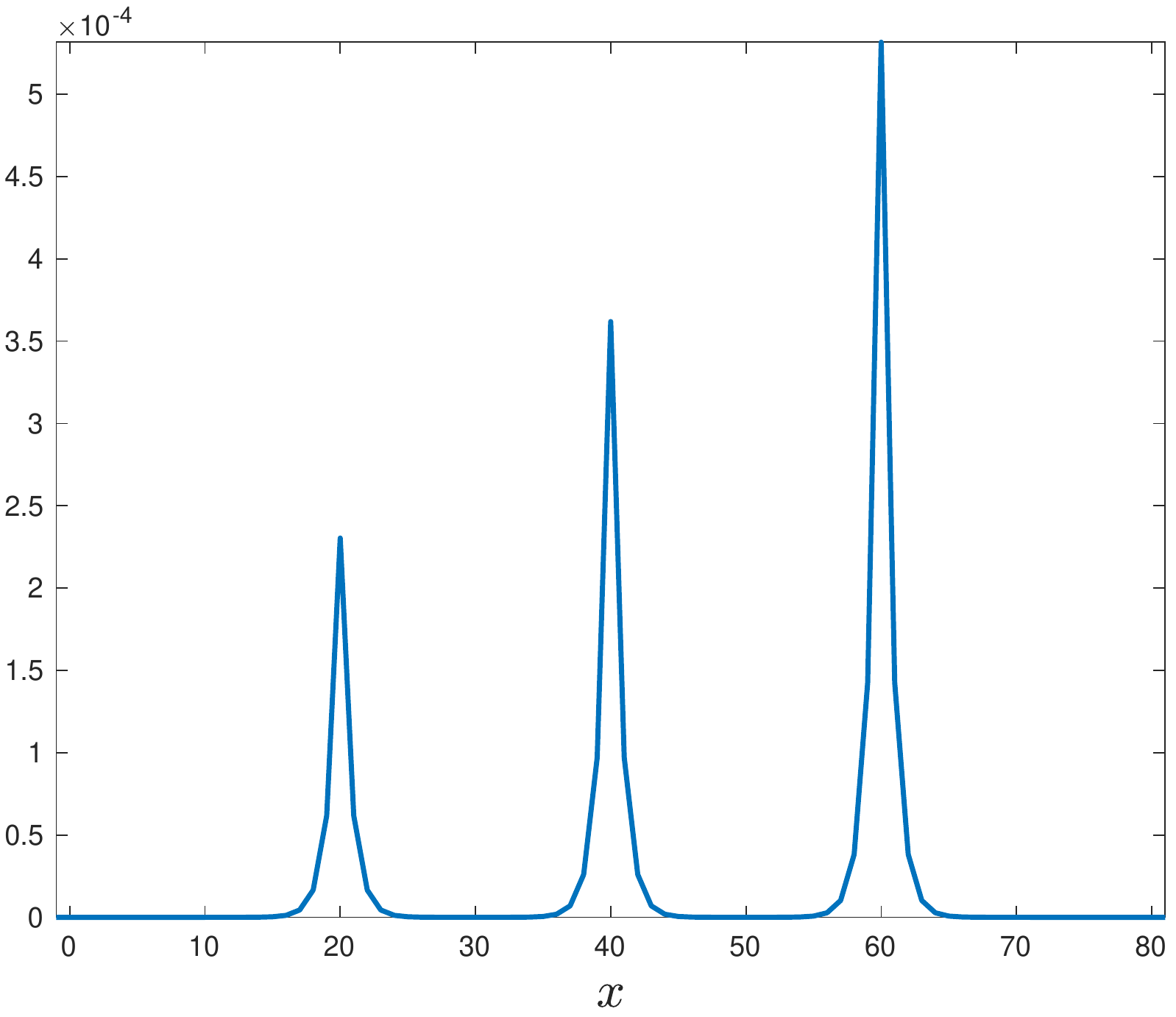}}}
\\
\centering
\subfigure[$n_{nb} = 10$.]{
{\includegraphics[width=0.48\textwidth,height=0.27\textwidth]{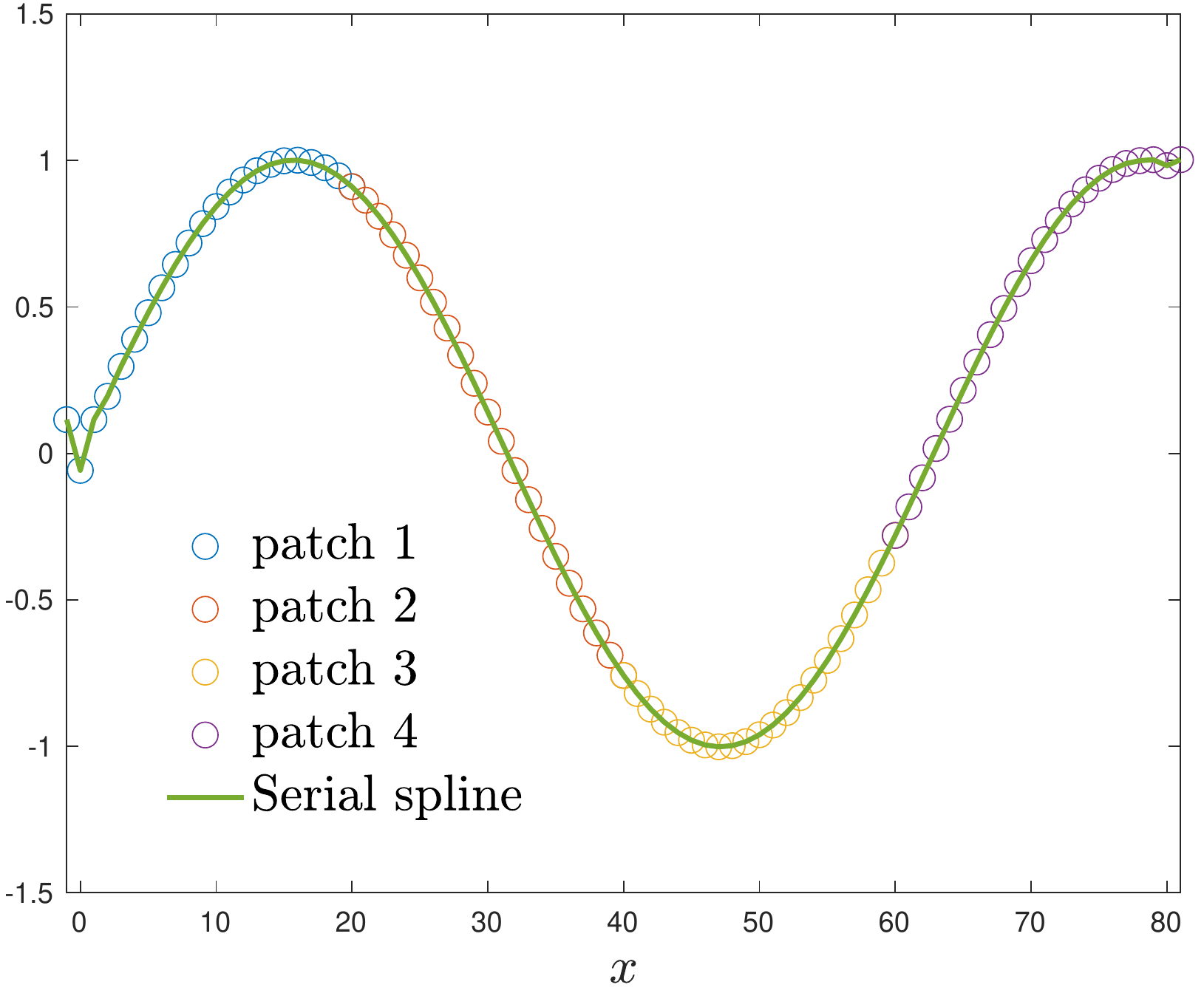}}
{\includegraphics[width=0.48\textwidth,height=0.27\textwidth]{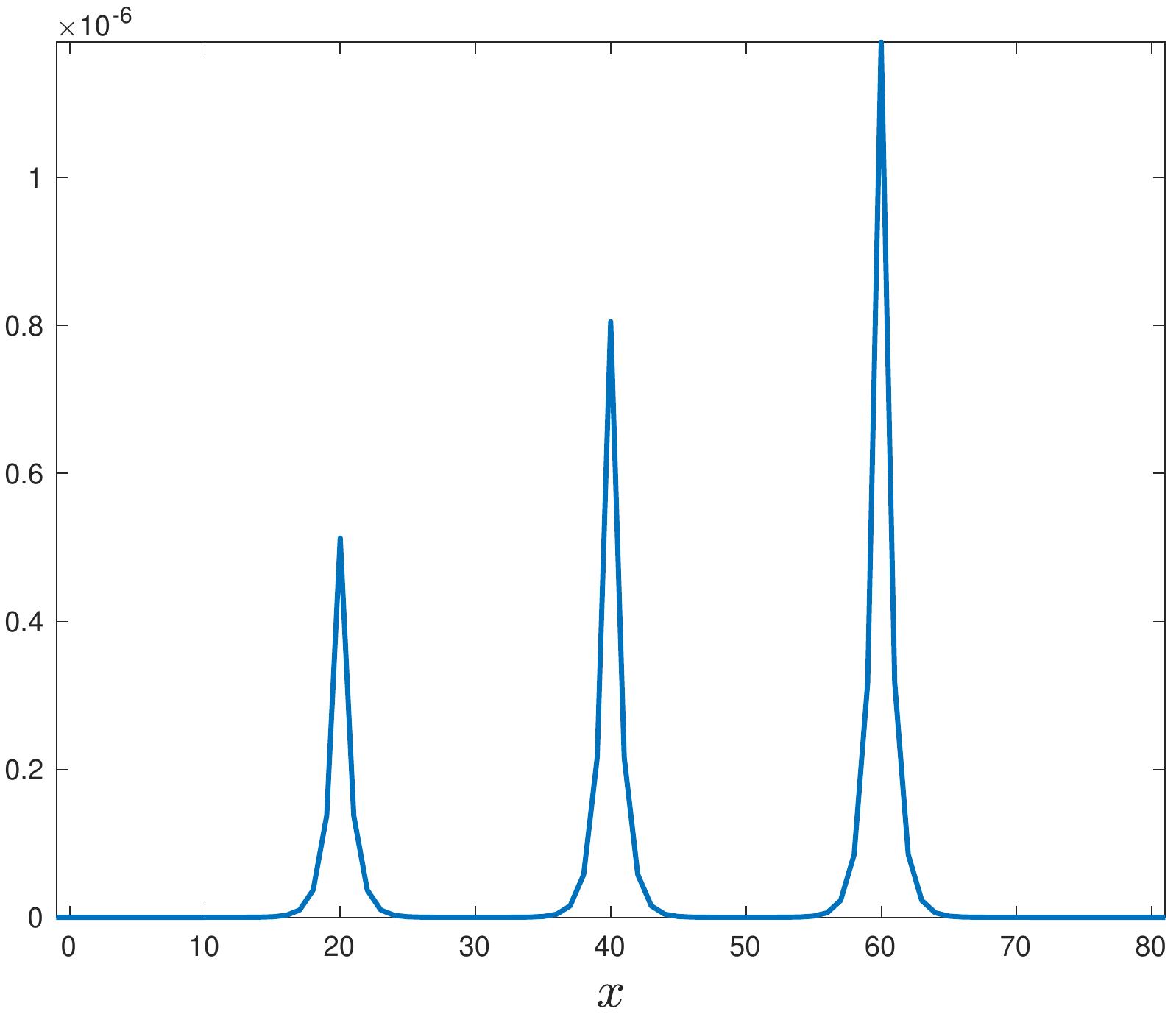}}}
\\
\centering
\subfigure[$n_{nb} = 20$.]{
{\includegraphics[width=0.48\textwidth,height=0.27\textwidth]{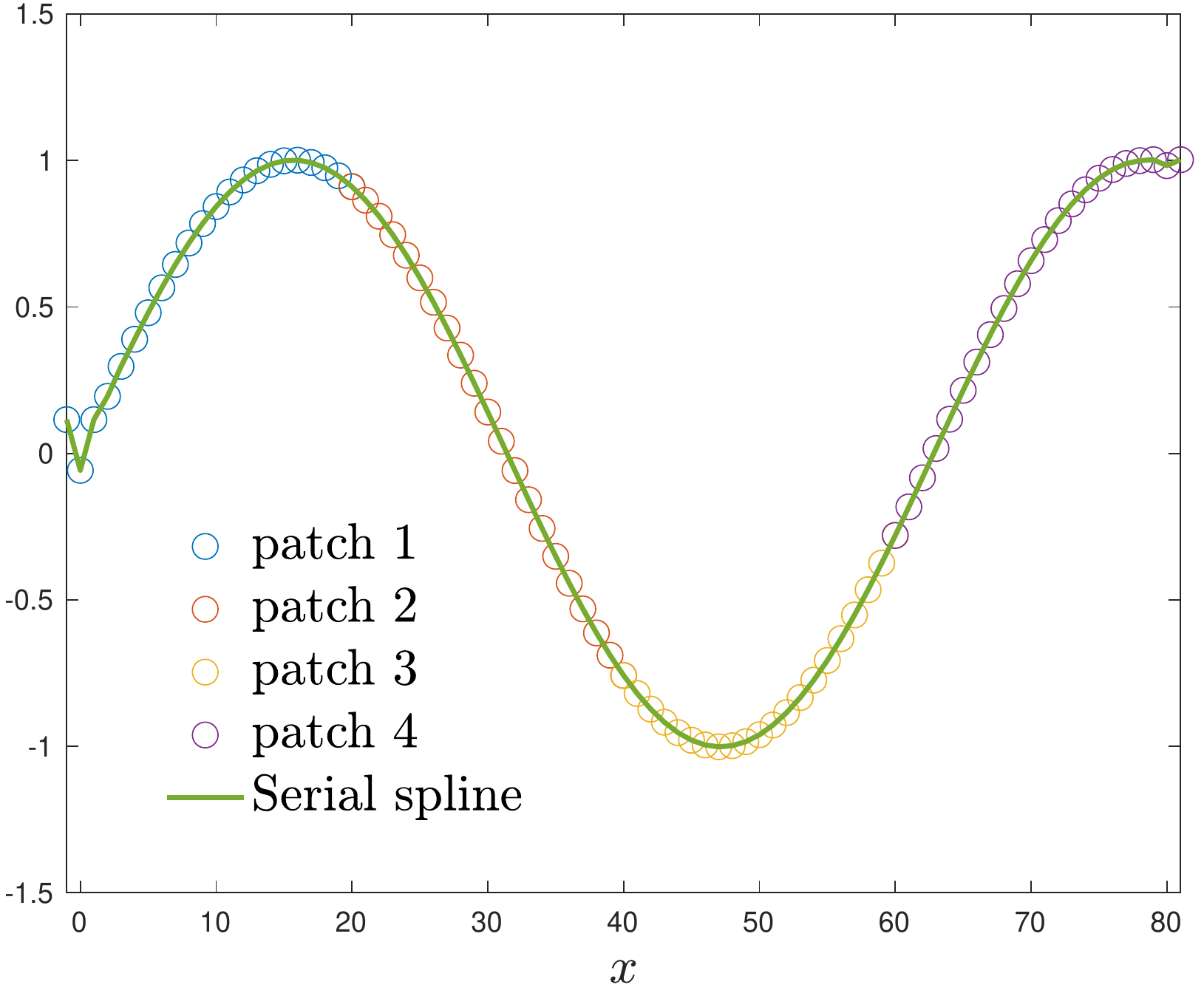}}
{\includegraphics[width=0.48\textwidth,height=0.27\textwidth]{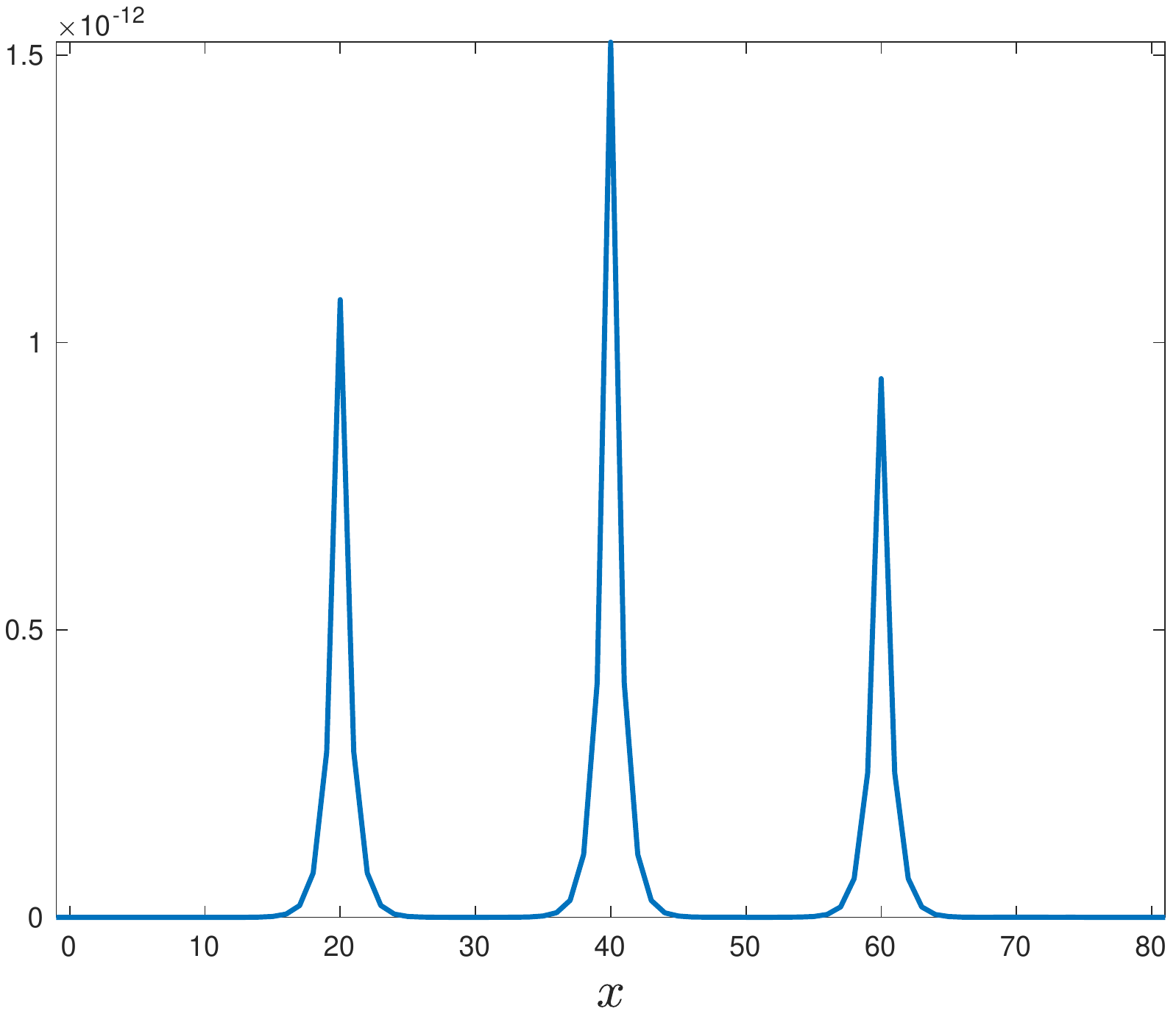}}}
\caption{\small Spline coefficients (left) and absolute errors (right) under PMBC ($n_{nb}=10$). The errors at the junction points are dramatically suppressed.    \label{spline_matched_HBC}}
\end{figure}

By contrast, the results under PMBC  are given in Figure \ref{spline_matched_HBC}. One can see that the errors are significantly smaller. When $N$ is fixed to be $161$, we find that $n_{nb} = 12$ can achieve relative error about $10^{-8}$ and $n_{nb} = 26$ can achieve that about $10^{-16}$. 

\begin{example}[Free advection of a 2-D Gaussian wavepacket\label{ex3}]

\begin{figure}[!h]
\centering
{\includegraphics[width=0.48\textwidth,height=0.27\textwidth]{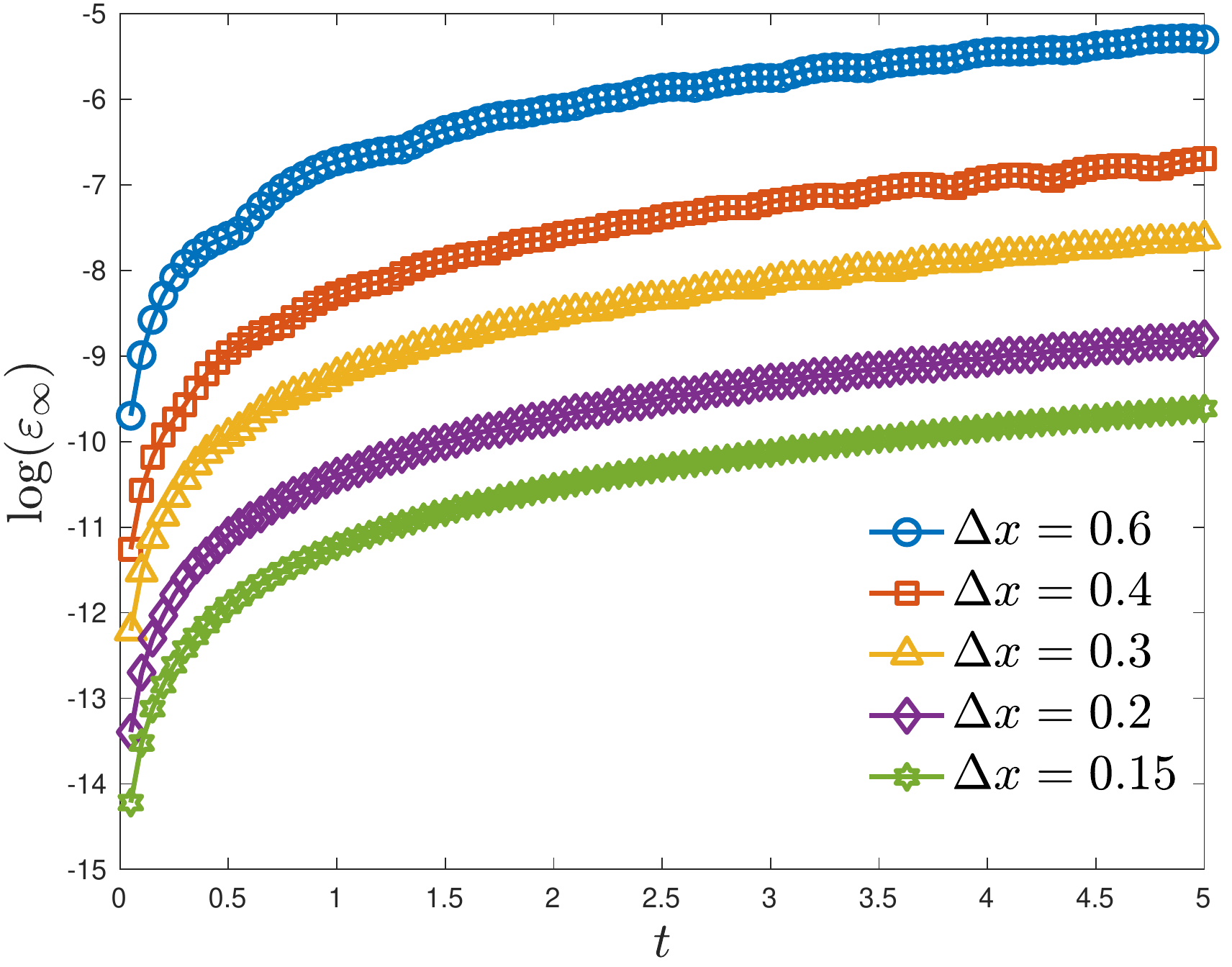}}
{\includegraphics[width=0.48\textwidth,height=0.27\textwidth]{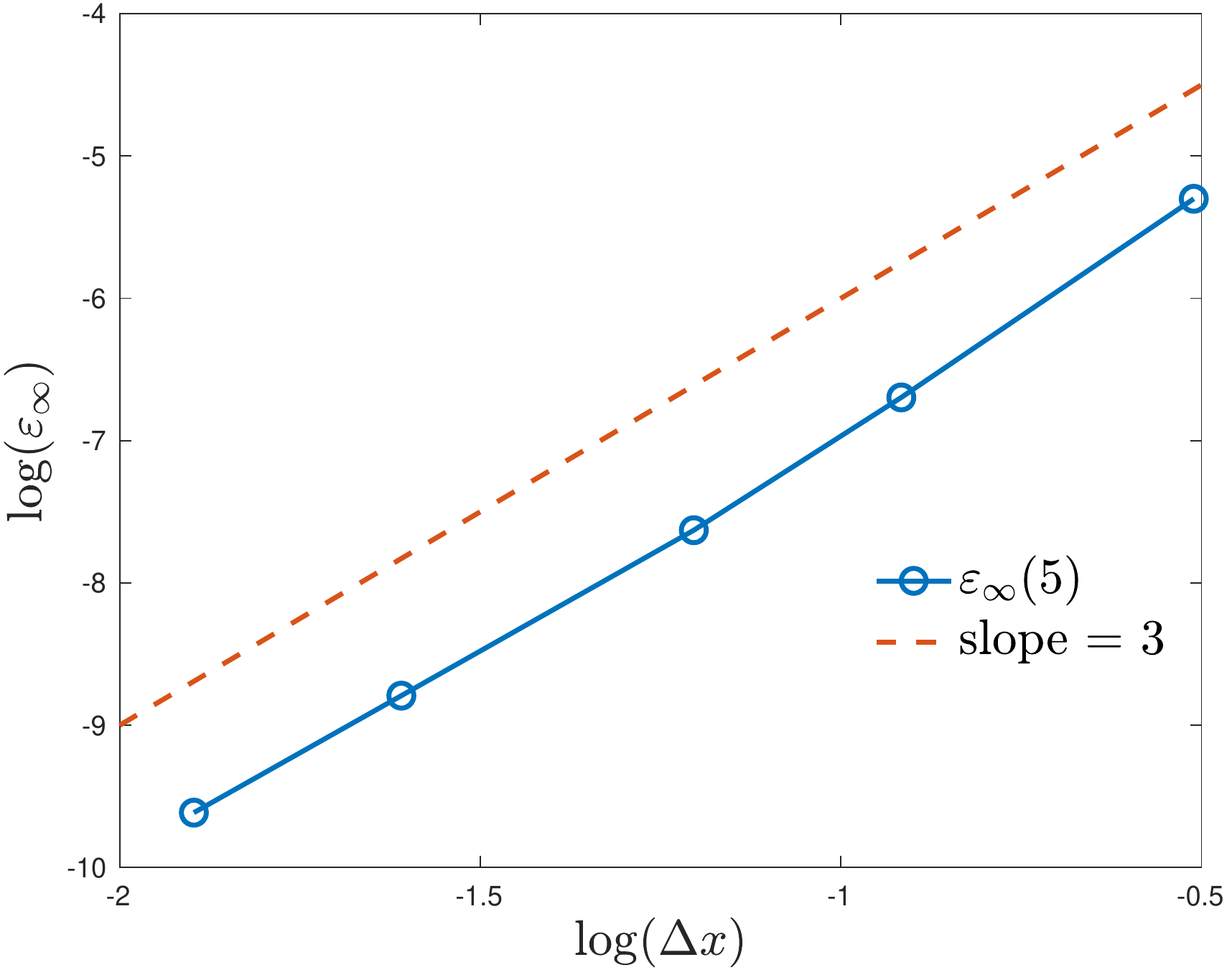}}
\caption{\small The time evolution of $\varepsilon_{\infty}(t)$ under the parallel spline reconstruction and different spatial stepsizes.  It perfectly matches the theoretical global convergence order $3$. \label{free_spline_convergence}}
\end{figure}

\begin{figure}[!h]
\centering
\subfigure[Serial cubic B-spline interpolation.]{
{\includegraphics[width=0.48\textwidth,height=0.27\textwidth]{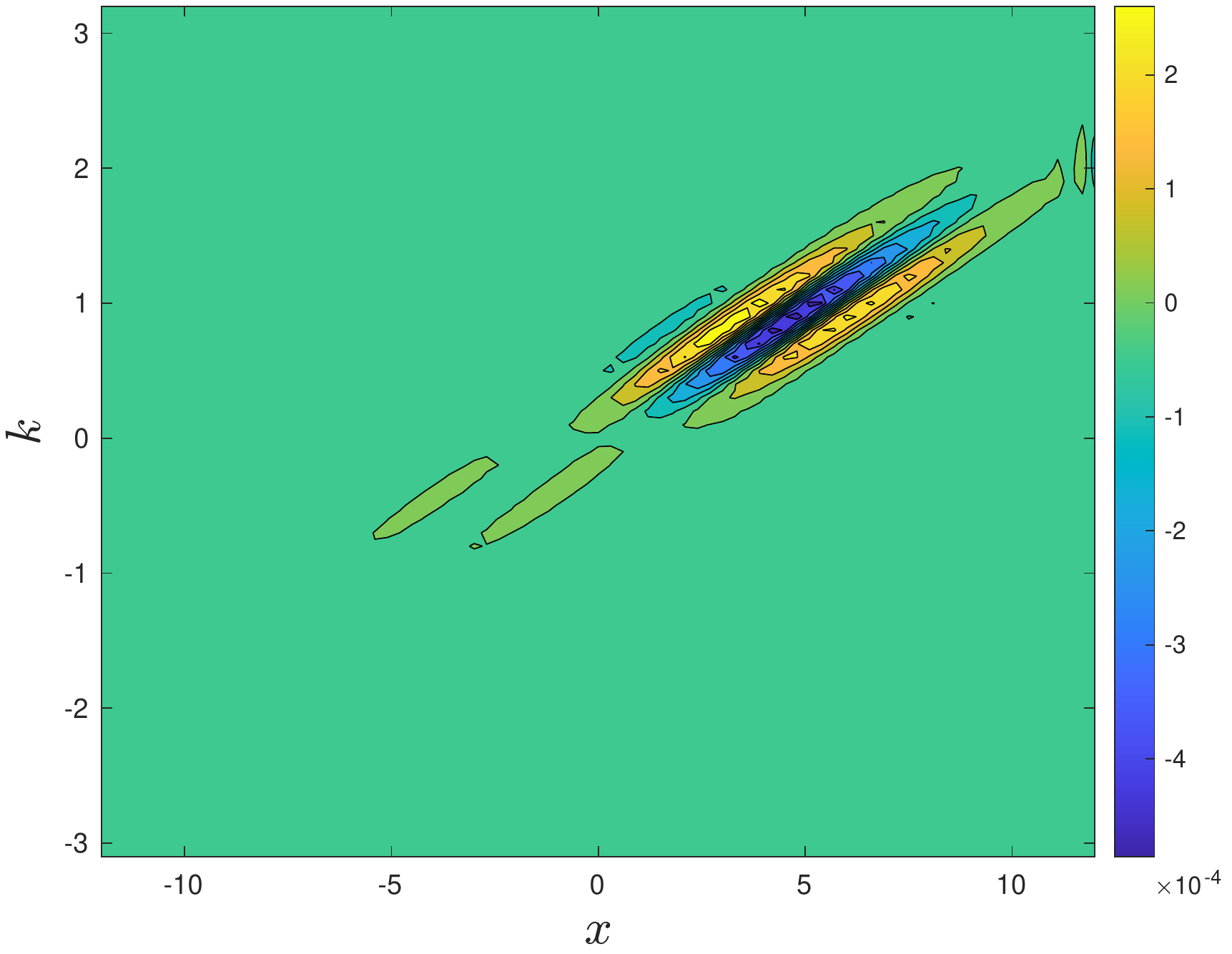}}
{\includegraphics[width=0.48\textwidth,height=0.27\textwidth]{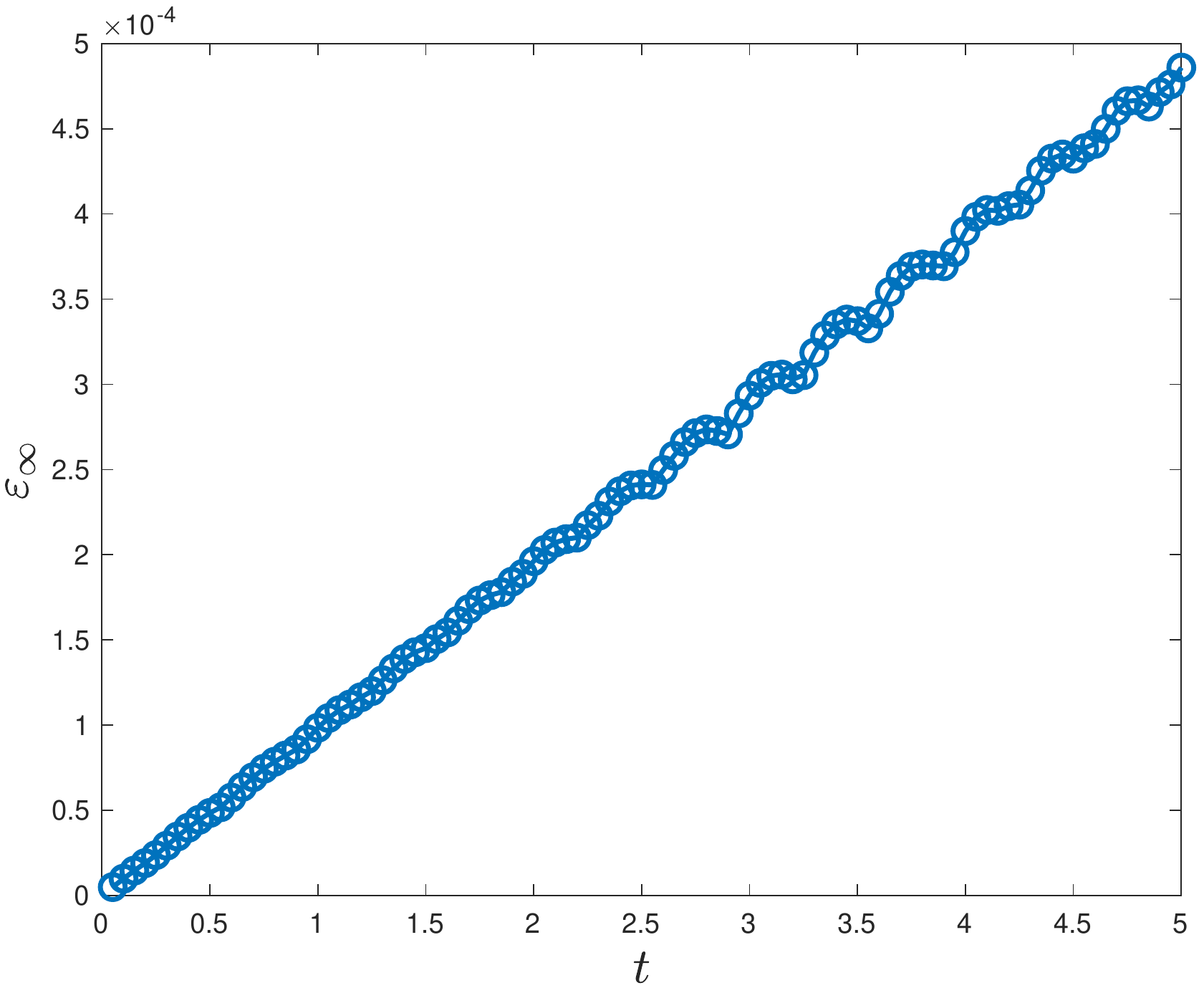}}}
\\
\centering
\subfigure[Parallel cubic B-spline interpolation with PMBC ($n_{nb} = 10$).]{
{\includegraphics[width=0.48\textwidth,height=0.27\textwidth]{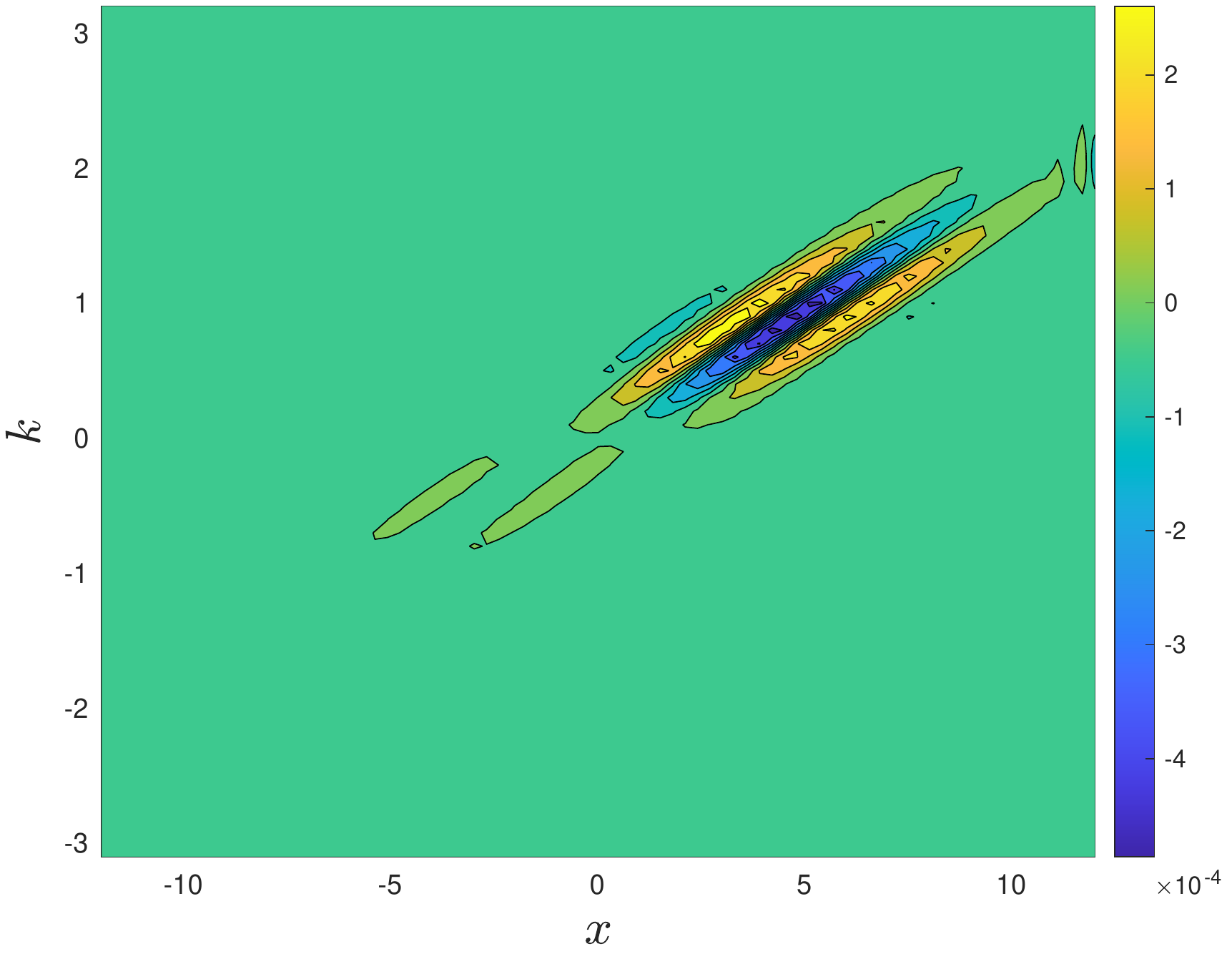}}
{\includegraphics[width=0.48\textwidth,height=0.27\textwidth]{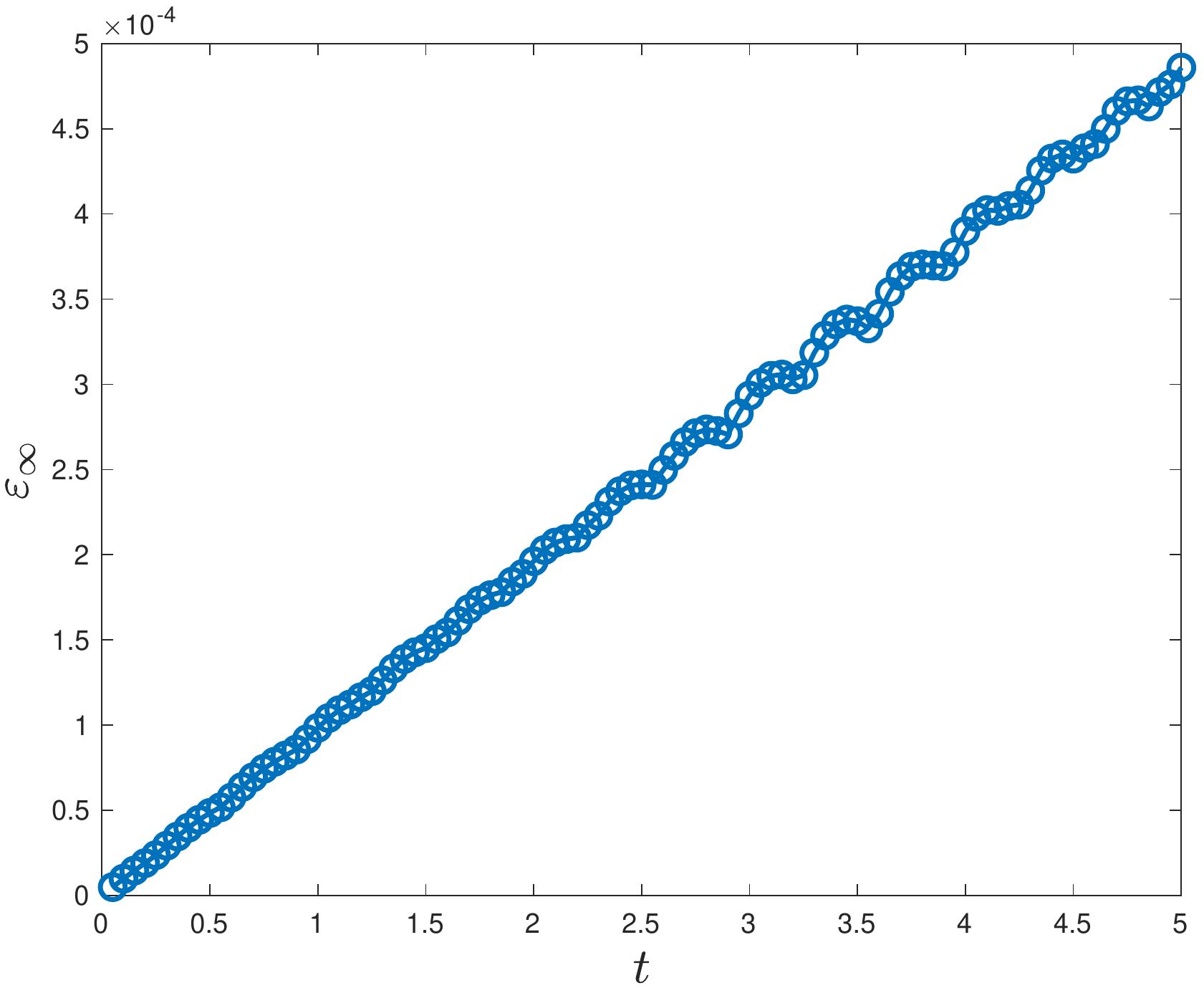}}}
\\
\centering
\subfigure[Parallel cubic B-spline interpolation with CLS-HBC.]{
{\includegraphics[width=0.48\textwidth,height=0.27\textwidth]{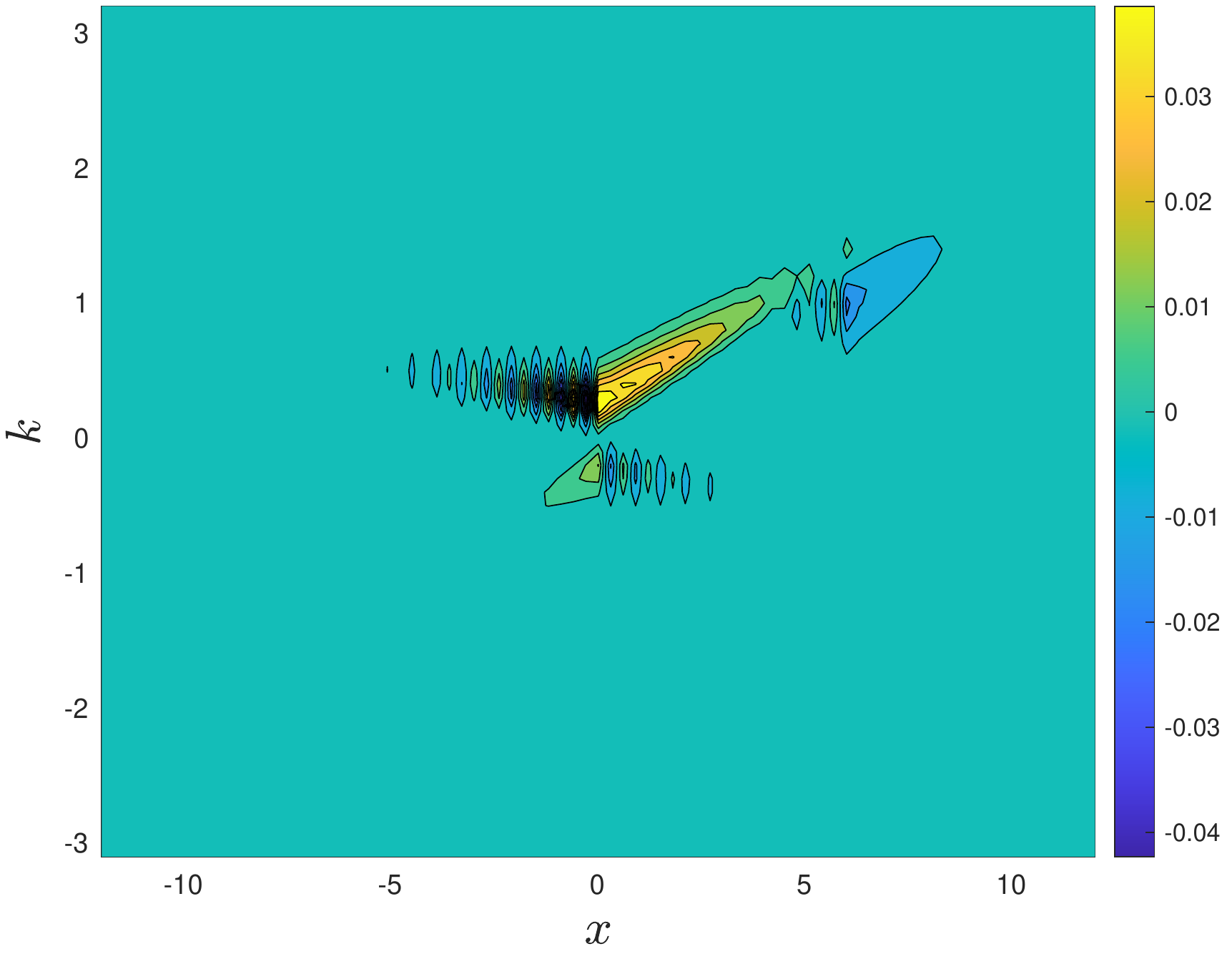}}
{\includegraphics[width=0.48\textwidth,height=0.27\textwidth]{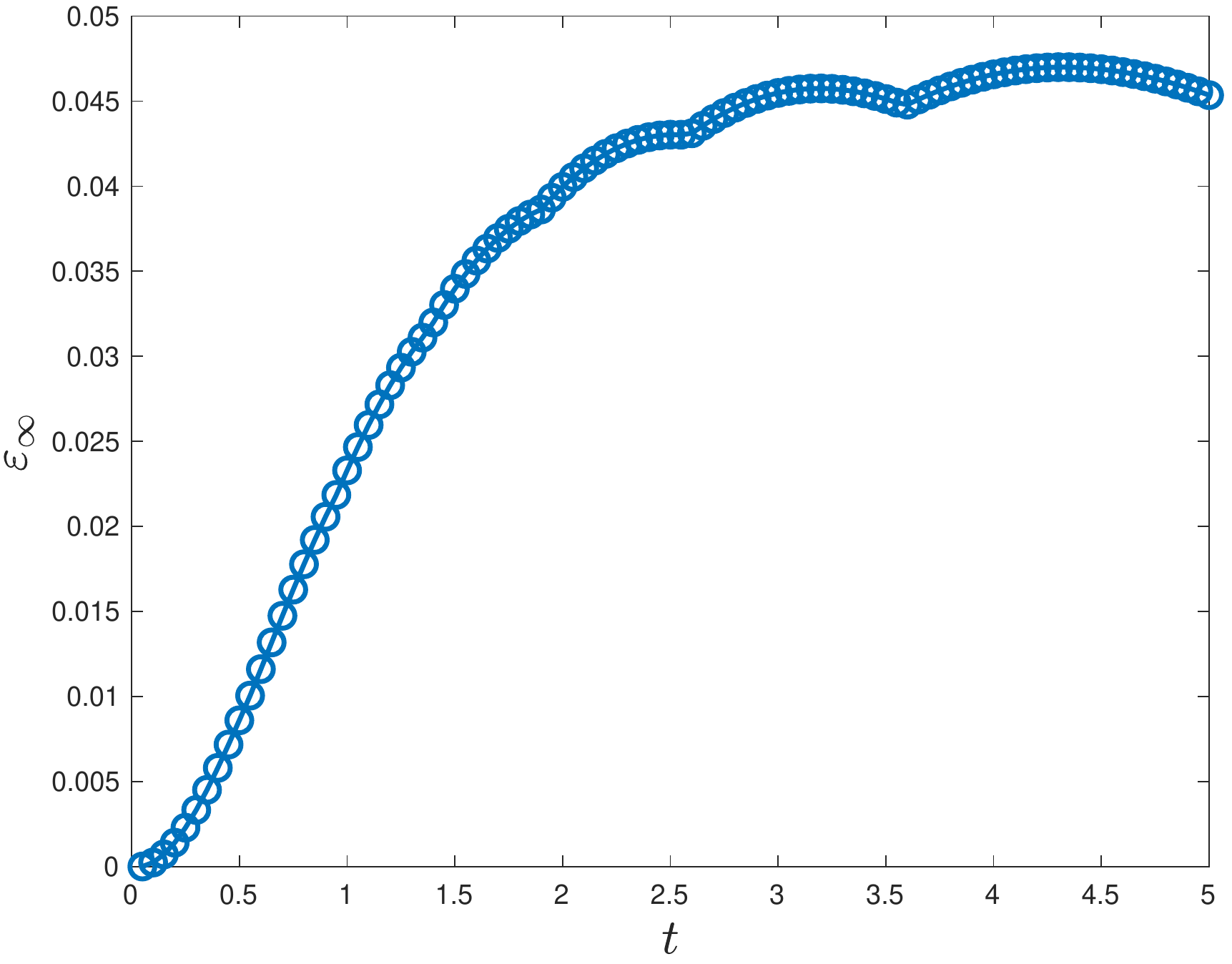}}}
\caption{\small A comparison of 2-D free advection: $f^{\textup{num}}(x, k, t) - f^{\textup{exact}}(x, k, t)$ at $t = 5$ (left) and the time evolution of $\varepsilon_{\infty}(t)$ (right) under $N = 81$. When PMBC is adopted, it produces almost the same results as that of the serial implementation. By contrast, when the CLS-HBC is adopted, small oscillations are observed at the junction points.   \label{Wigner2d_N81}}
\end{figure}

\begin{figure}[!h]
\centering
\subfigure[Serial cubic B-spline interpolation.]{
{\includegraphics[width=0.48\textwidth,height=0.27\textwidth]{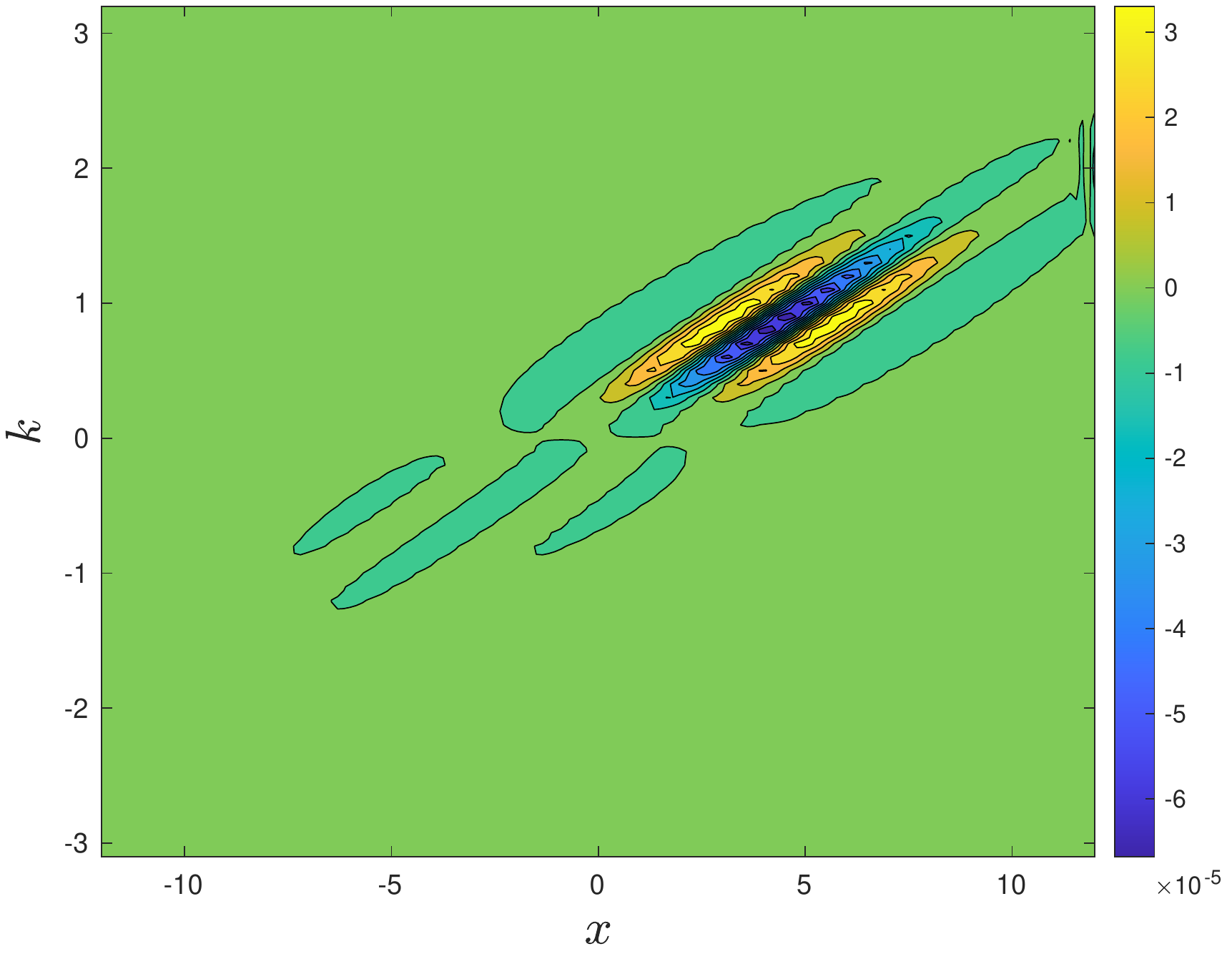}}
{\includegraphics[width=0.48\textwidth,height=0.27\textwidth]{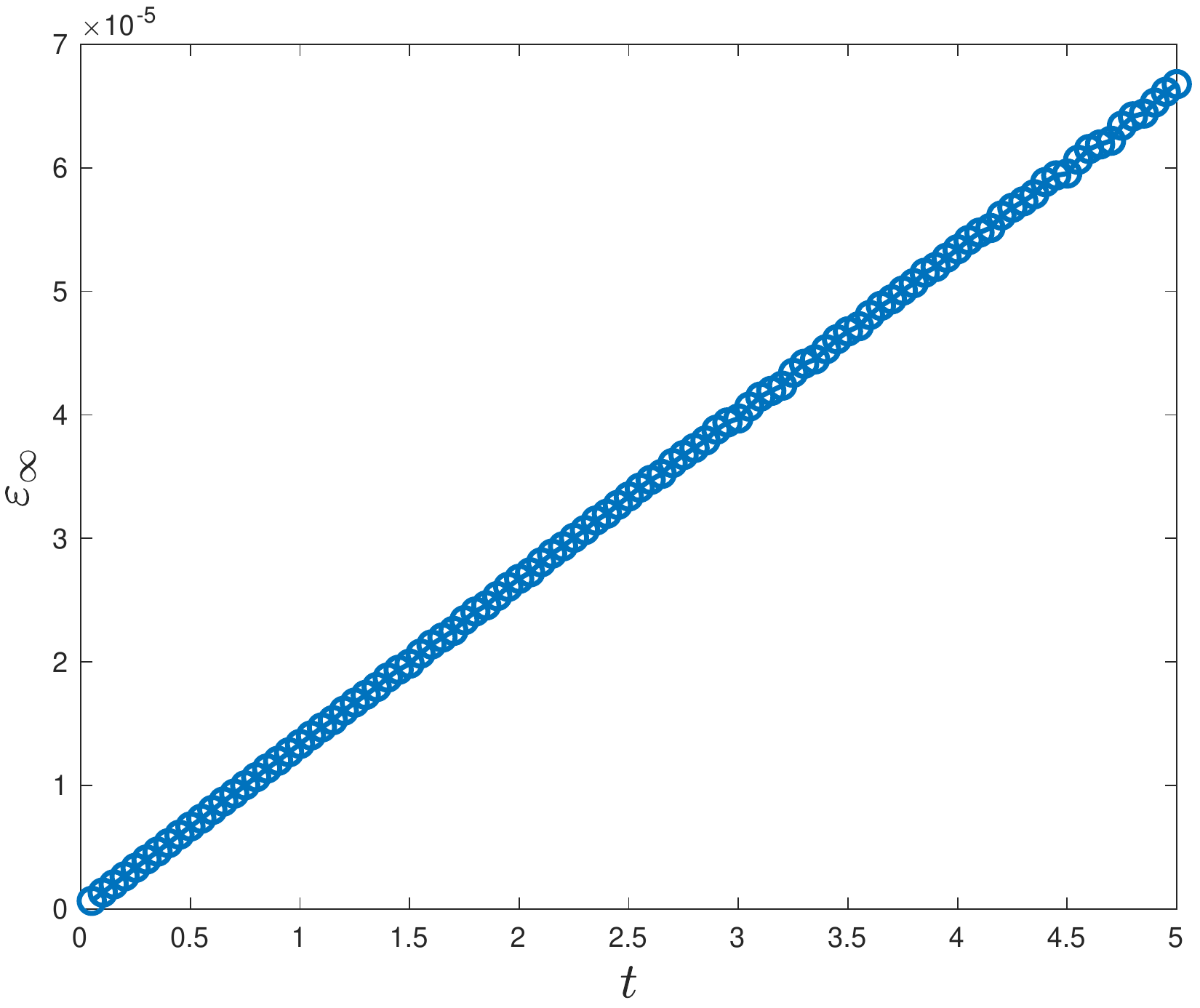}}}
\\
\centering
\subfigure[Parallel cubic B-spline interpolation with PMBC ($n_{nb} = 10$).]{
{\includegraphics[width=0.48\textwidth,height=0.27\textwidth]{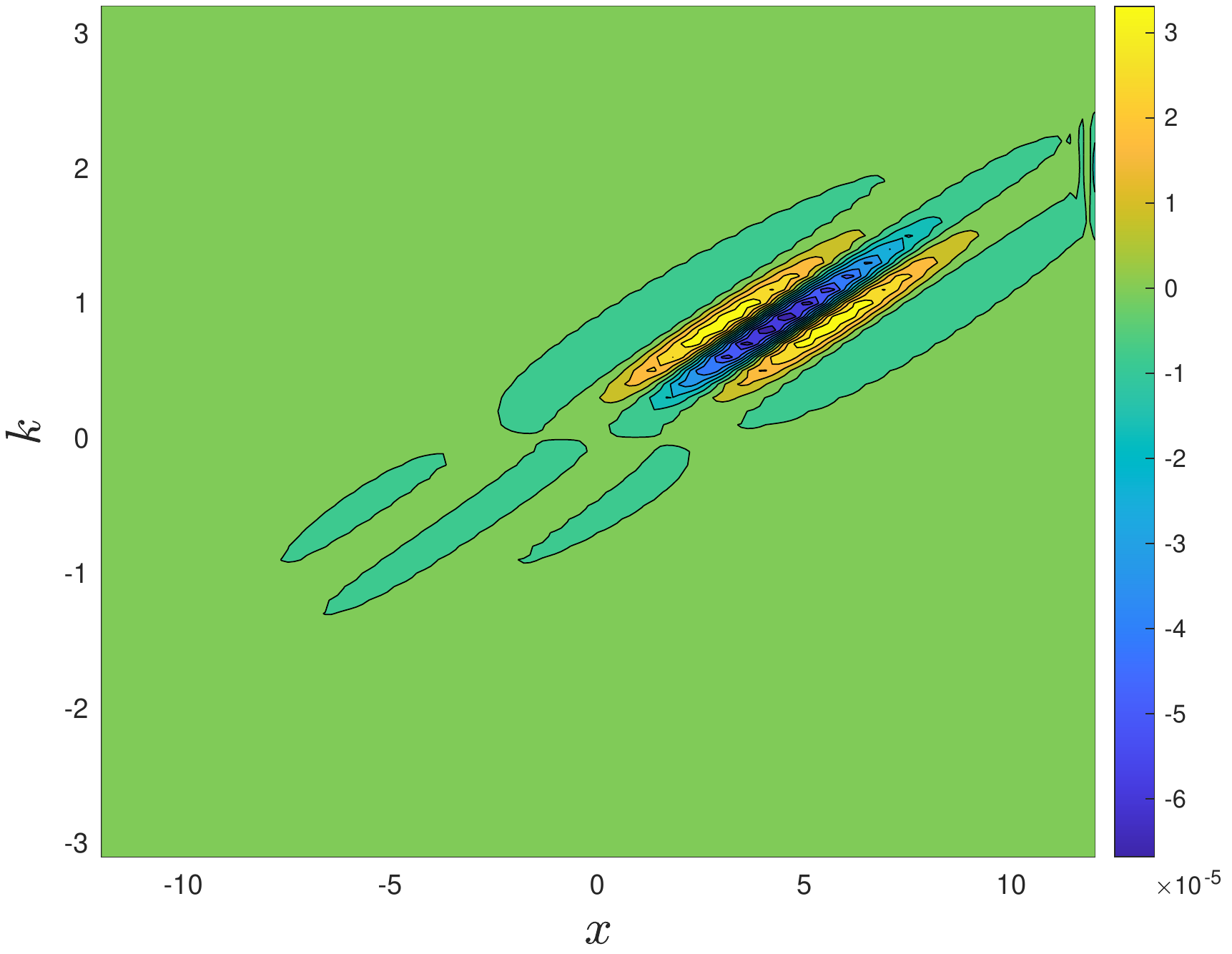}}
{\includegraphics[width=0.48\textwidth,height=0.27\textwidth]{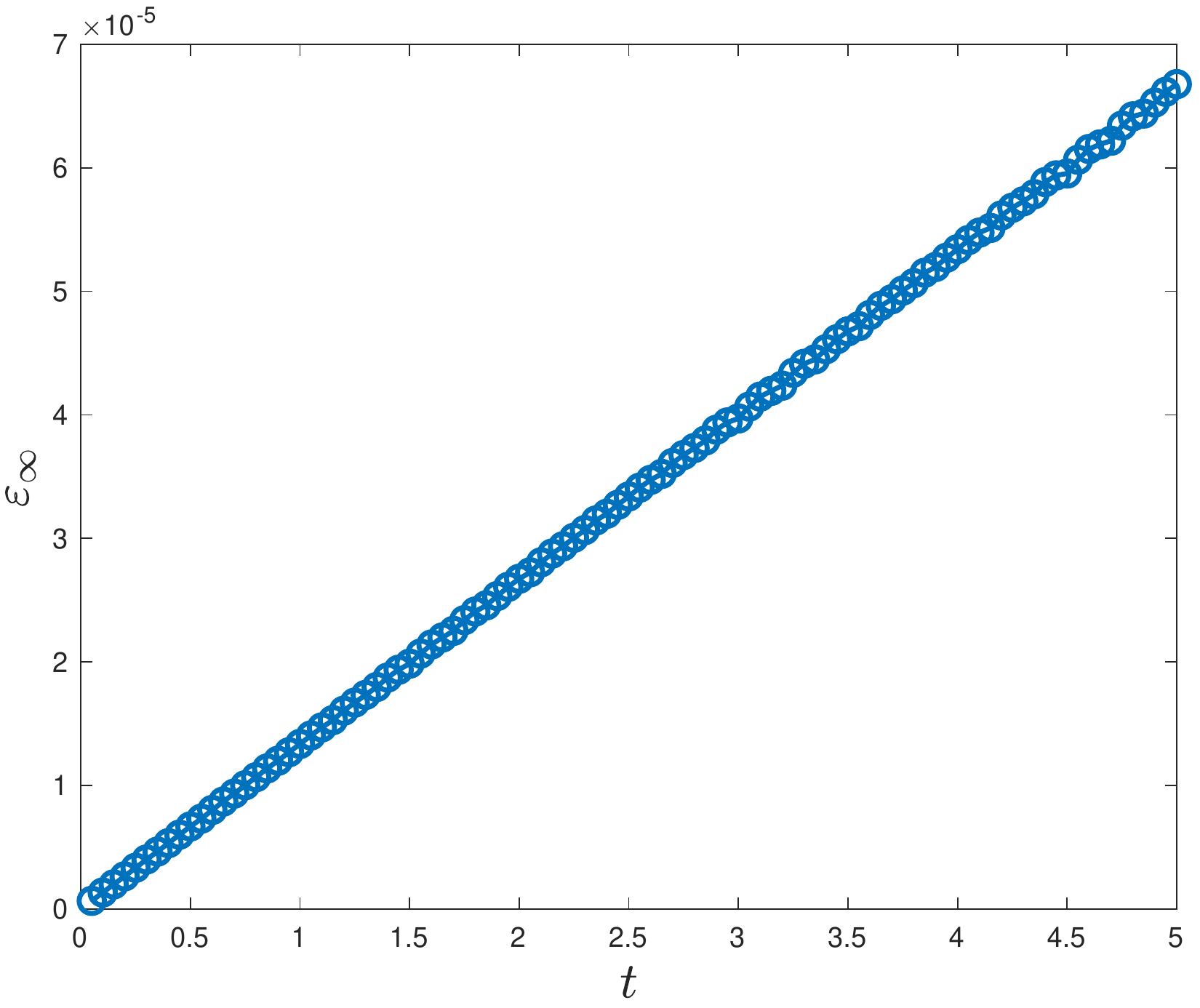}}}
\\
\centering
\subfigure[Parallel cubic B-spline interpolation with CLS-HBC.]{
{\includegraphics[width=0.48\textwidth,height=0.27\textwidth]{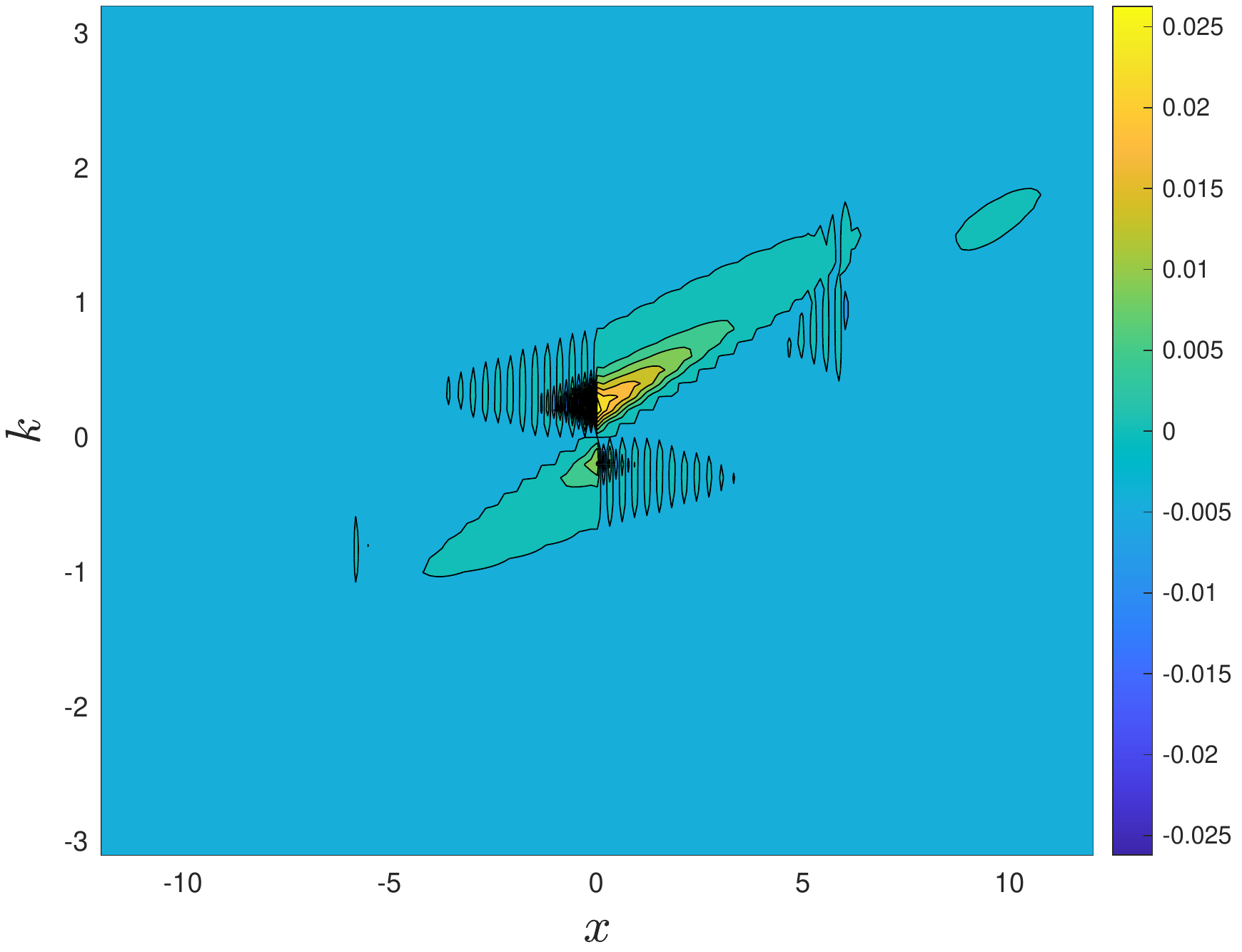}}
{\includegraphics[width=0.48\textwidth,height=0.27\textwidth]{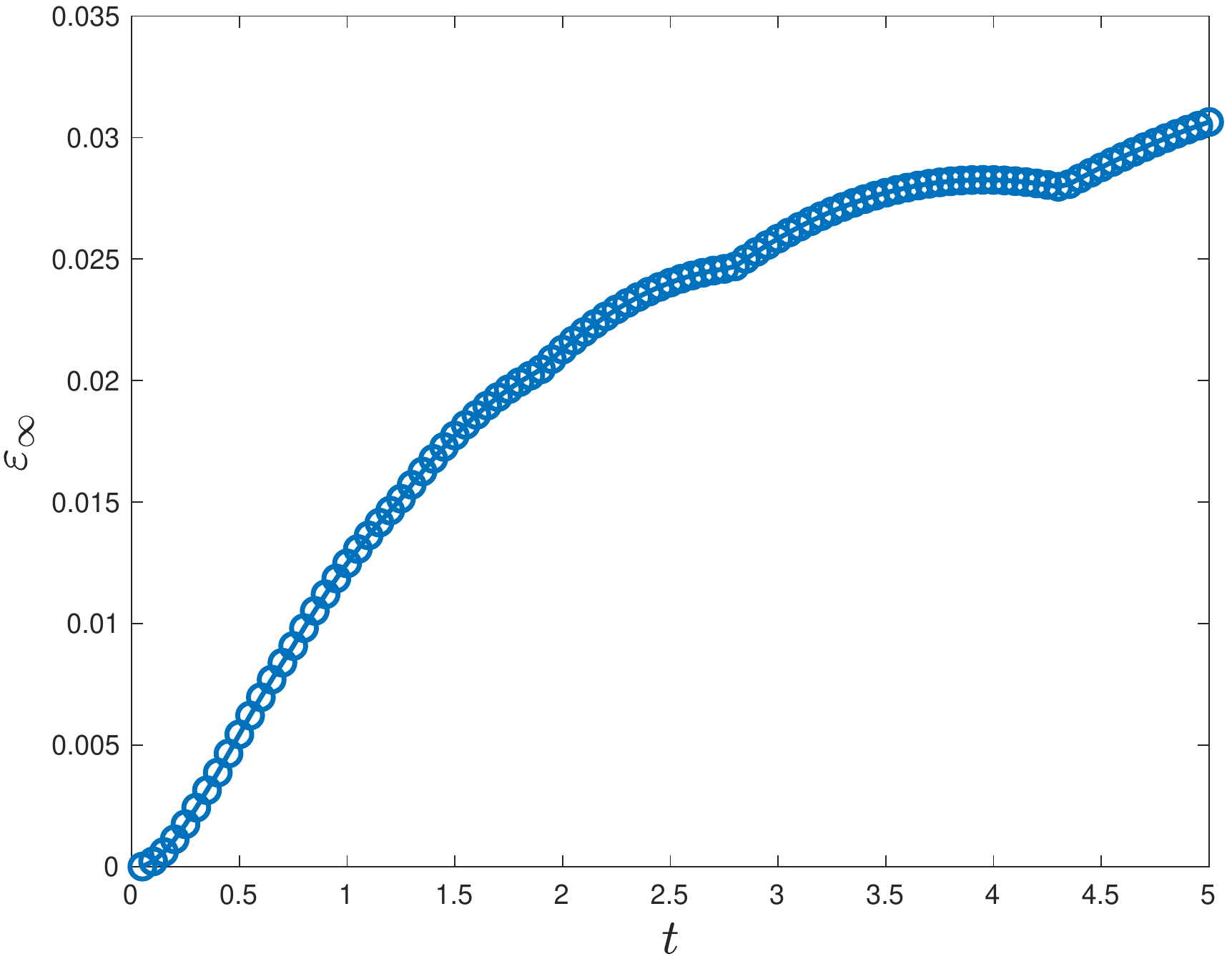}}}
\caption{\small A comparison of 2-D free advection: $f^{\textup{num}}(x, k, t) - f^{\textup{exact}}(x, k, t)$ at $t = 5$ (left) and the time evolution of $\varepsilon_{\infty}(t)$ (right) under $N = 161$. Numerical results are further improved under smaller spatial spacing.   \label{Wigner2d_N161}}
\end{figure}

The second test problem is the free-advection of the Wigner function in 2-D phase space: 
 \begin{equation}
\frac{\partial}{\partial t} f(x, k, t) = \frac{\hbar k}{m} \frac{\partial }{\partial x} f(x, k, t),
\end{equation}
with 
\begin{equation}
f(x, k, 0) = \frac{1}{\pi} \exp\left(-\frac{x^2}{2a^2} - 2a^2 (k-k_0)^2\right).
\end{equation}
The exact solution reads that
\begin{equation}
 f(x, k, t) =  \frac{1}{\pi} \exp\left(-\frac{(x - \frac{\hbar k t}{m})^2}{2a^2} - 2a^2 (k-k_0)^2\right).
\end{equation}
Here we take $a = 1$, $\hbar = m = 1$, $k_0 = 0.5$. The final time is $t_{fin} = 5$ with time step $\tau = 0.05$.
\end{example}

 To measure the numerical error, we adopt the $l^\infty$-error as the metric
\begin{equation}
\varepsilon_{\infty}(t) = \max_{(x, k) \in \mathcal{X} \times \mathcal{K}}| f^{\textup{num}}(x, k, t) - f^{\textup{exact}}(x, k, t) |,
\end{equation}
where $f^{\textup{num}}$ and $f^{\textup{exact}}$ denote the solutions produced by the spline interpolation and exact one, respectively. We make a comparison of two kinds of effective Hermite boundary conditions. When $N = 81$ and $n_{nb} = 10$ are fixed, one can see in Figure \ref{Wigner2d_N81} that the performance of the local splines under PMBC are almost the same as that of the serial spline, while the solutions under CLS-HBC exhibit small oscillations around the junction regions.  Such trend is also observed when further increasing $N$ to 161. As presented in Figure \ref{Wigner2d_N161}, the $l^\infty$-error $\varepsilon_{\infty}(5)$ decreases from $4.86\times 10^{-4}$ to $6.68 \times 10^{-5}$ and the convergence order is $3$ (see Figure \ref{free_spline_convergence}).

\subsection{The influence of different spline boundary conditions}

\begin{figure}[!h]
\centering
\subfigure[Evolution of errors. (left: Neumann boundary, right: natural boundary)]{
{\includegraphics[width=0.48\textwidth,height=0.27\textwidth]{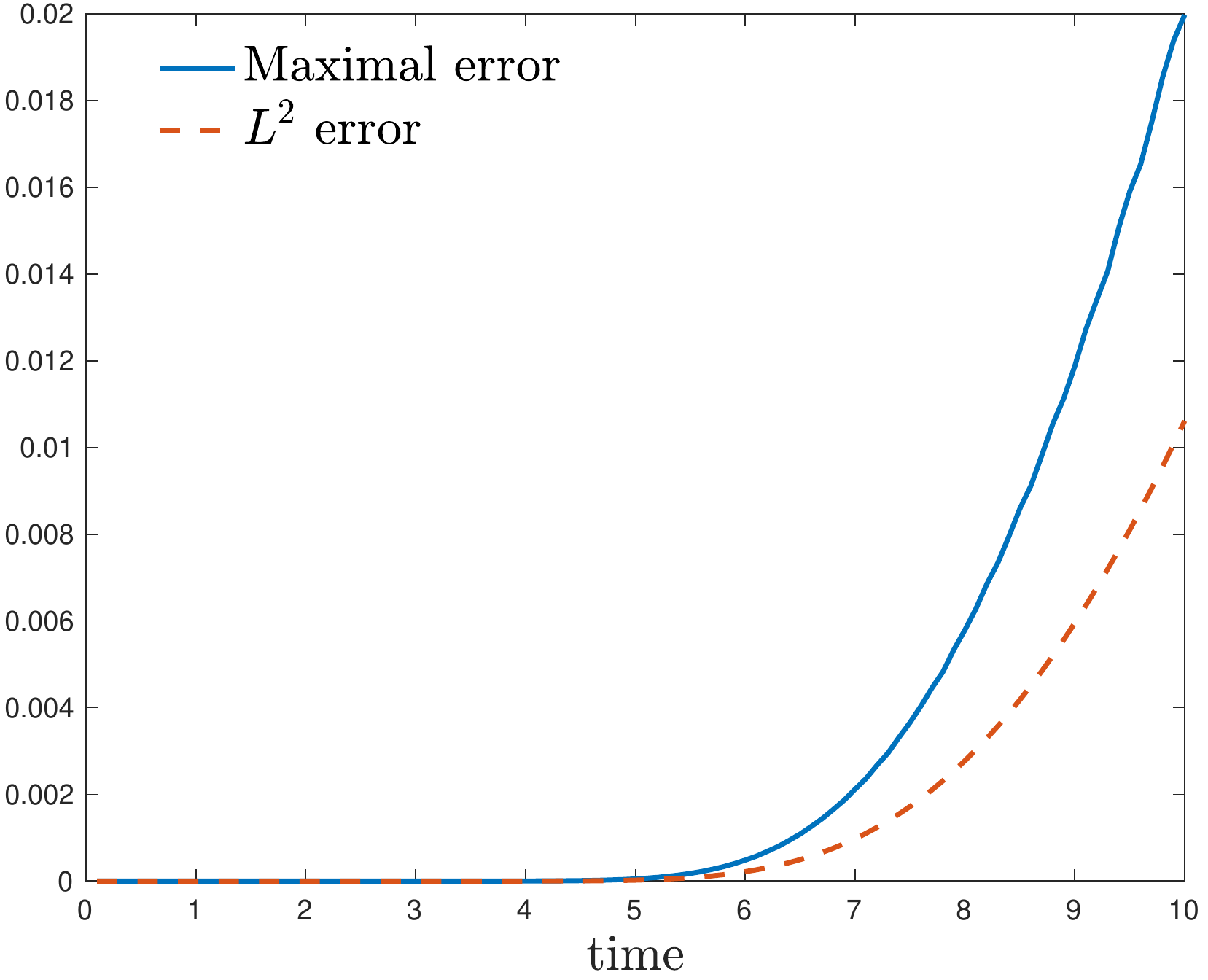}}
{\includegraphics[width=0.48\textwidth,height=0.27\textwidth]{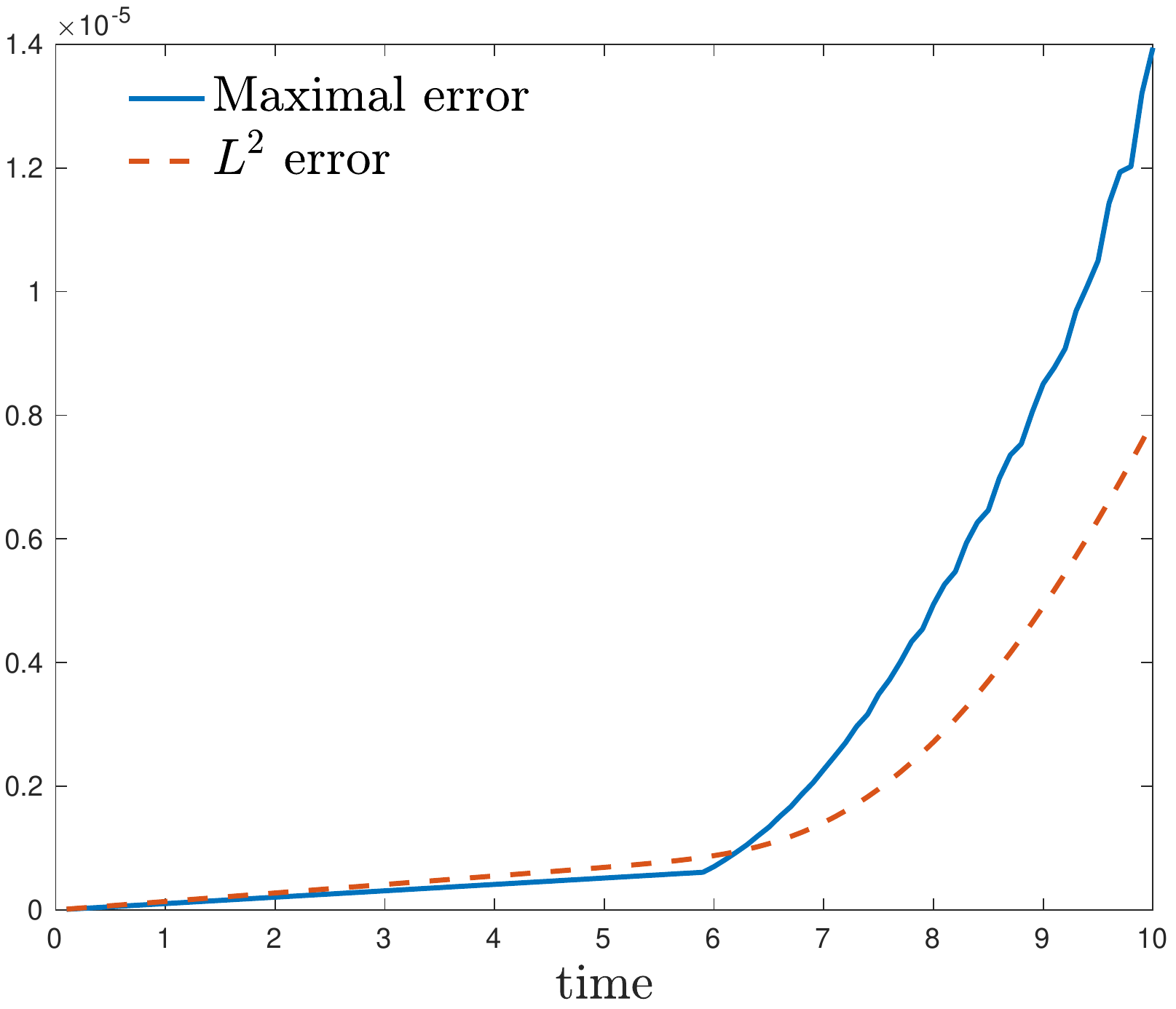}}}
\\
\centering
\subfigure[Difference at $t = 3$.  (left: Neumann boundary, right: natural boundary)]{
{\includegraphics[width=0.48\textwidth,height=0.27\textwidth]{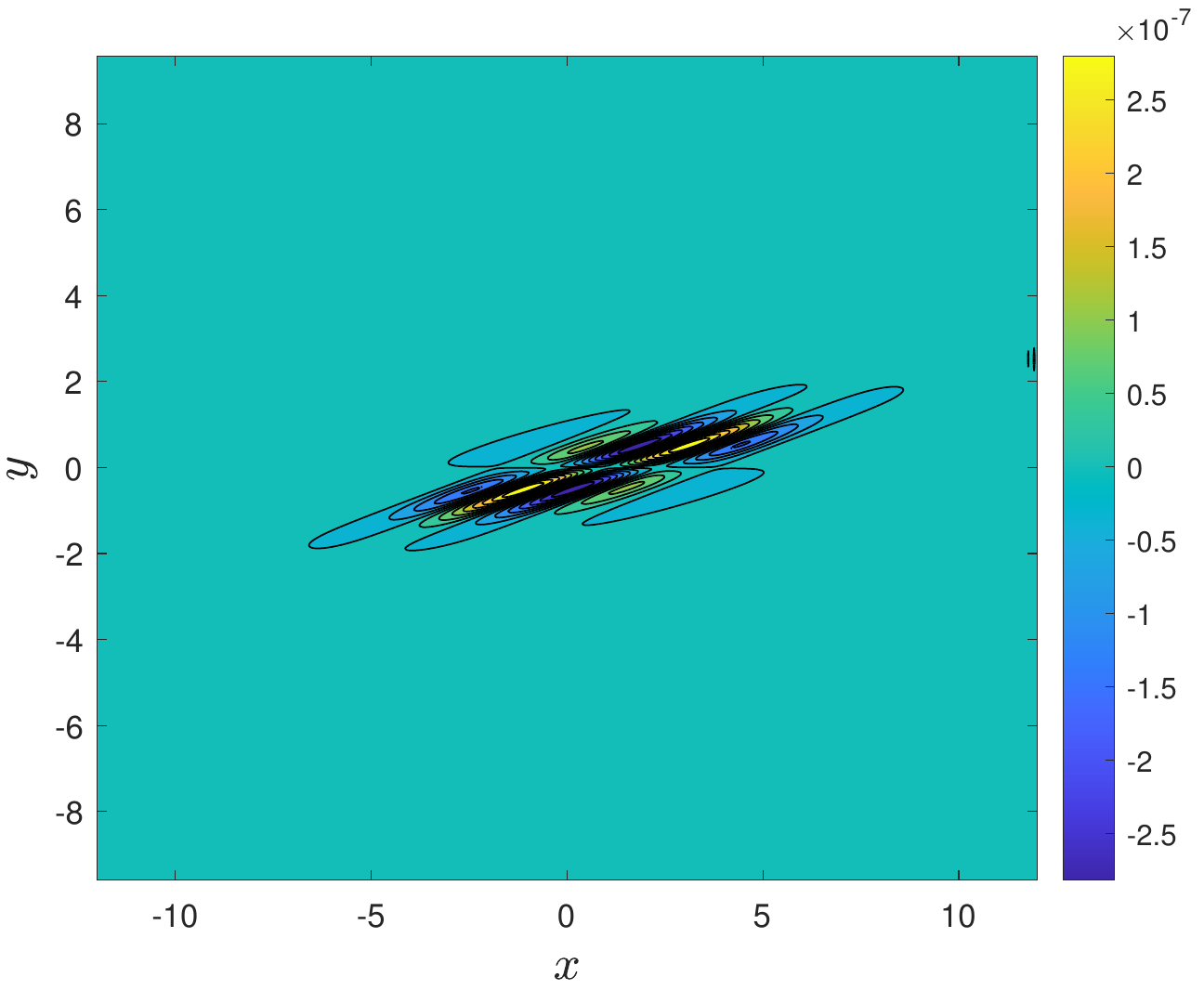}}
{\includegraphics[width=0.48\textwidth,height=0.27\textwidth]{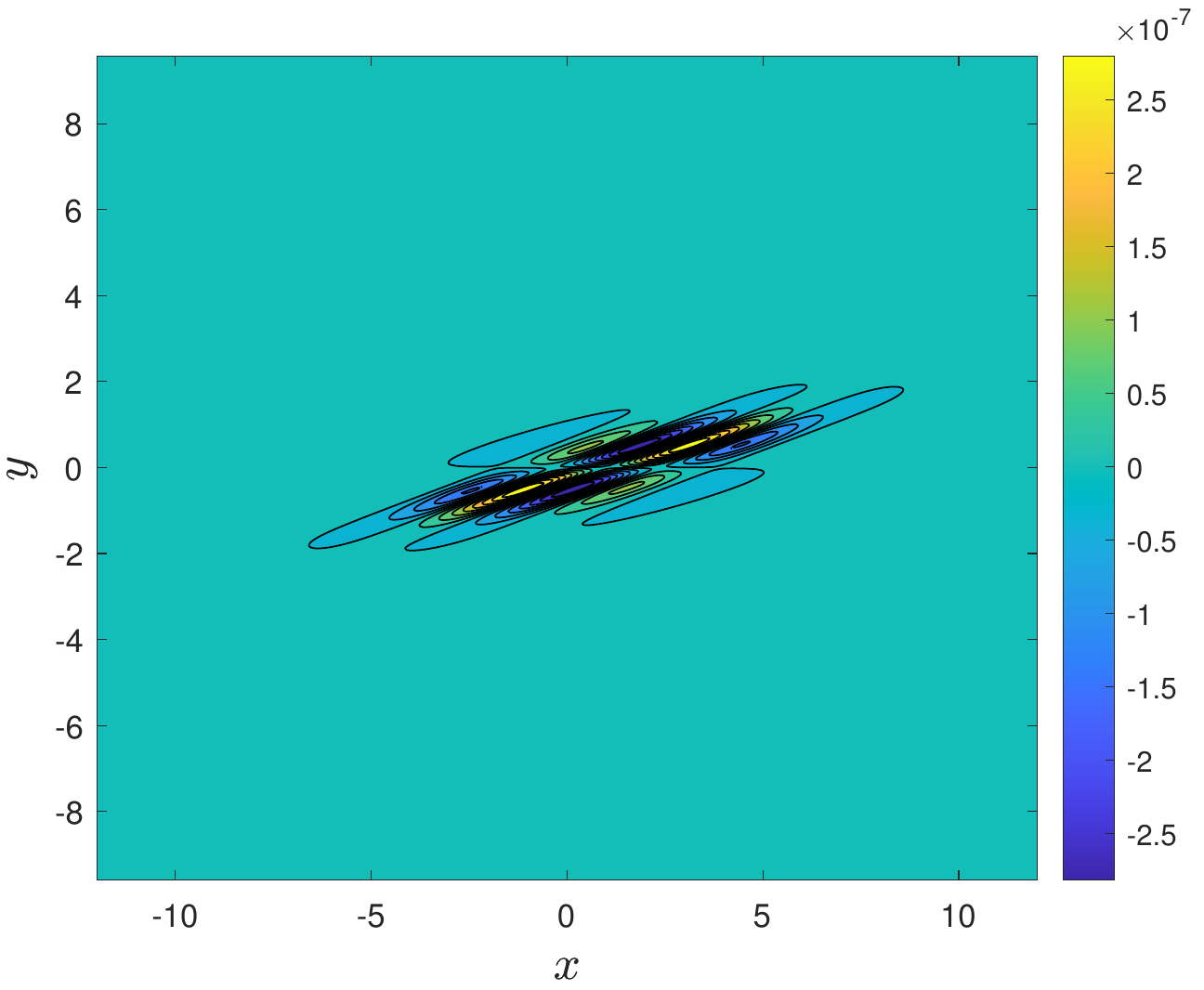}}}
\\
\centering
\subfigure[Difference at $t = 4$.  (left: Neumann boundary, right: natural boundary)]{
{\includegraphics[width=0.48\textwidth,height=0.27\textwidth]{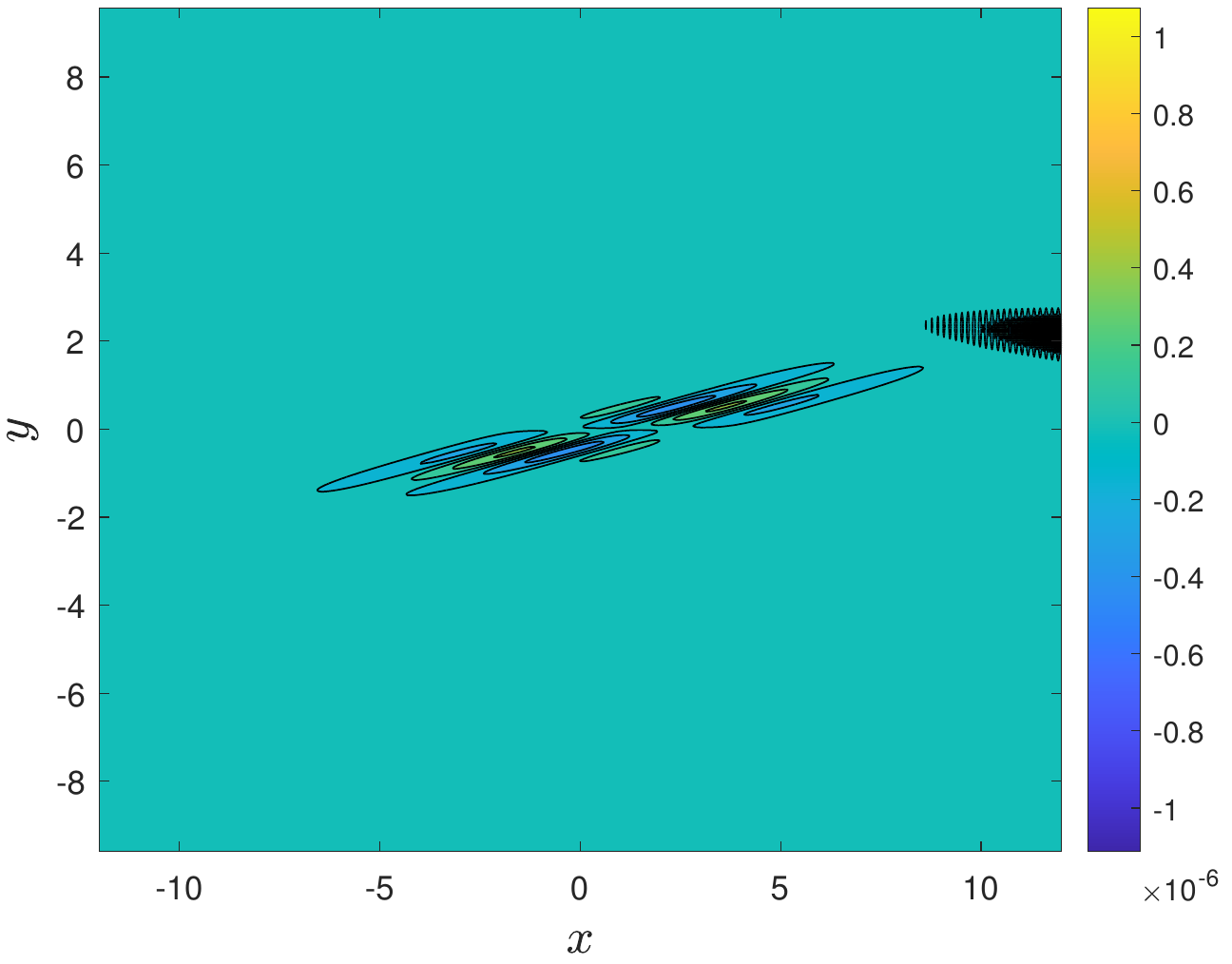}}
{\includegraphics[width=0.48\textwidth,height=0.27\textwidth]{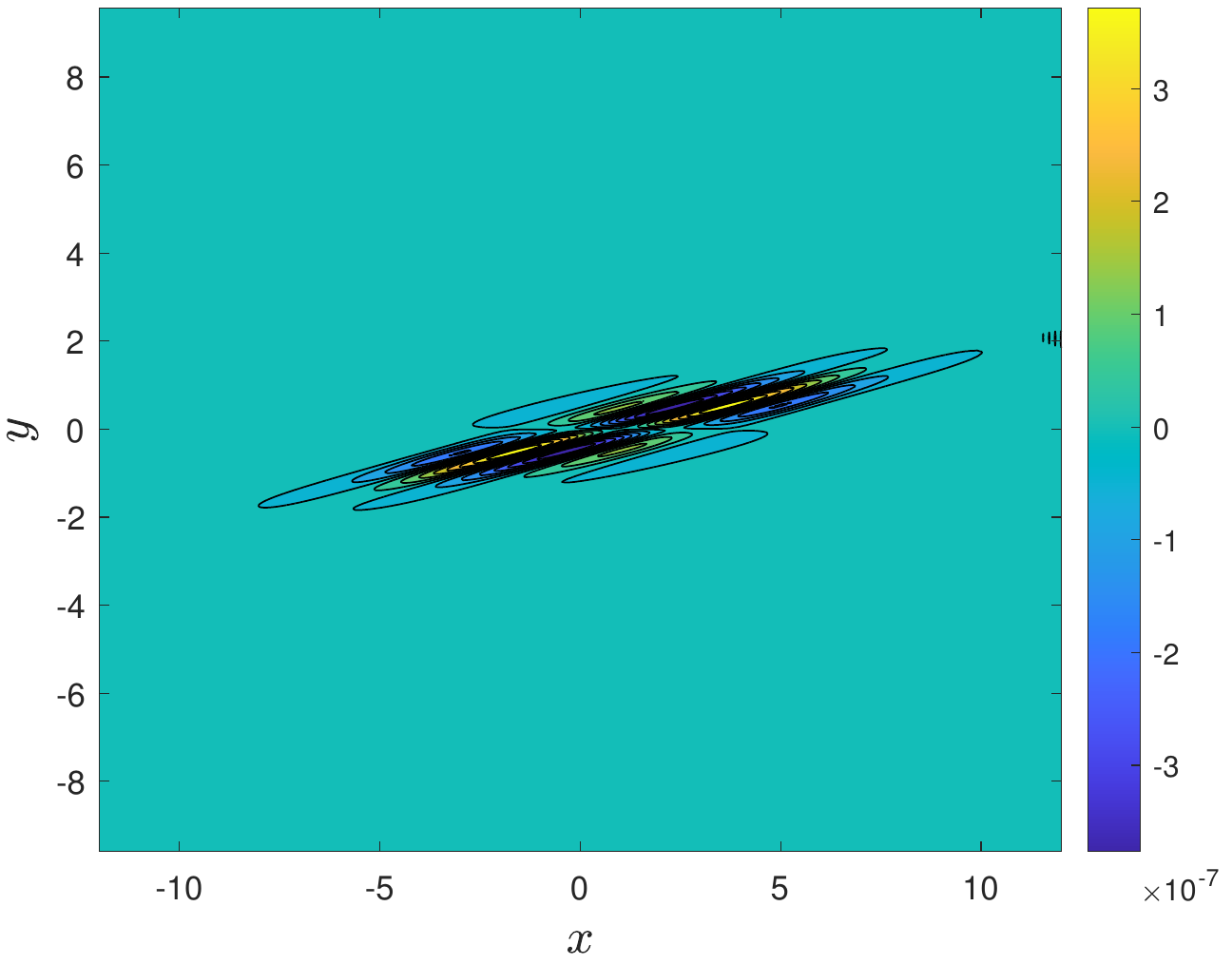}}}
\\
\centering
\subfigure[Difference at $t = 10$. (left: Neumann boundary, right: natural boundary)]{
{\includegraphics[width=0.48\textwidth,height=0.27\textwidth]{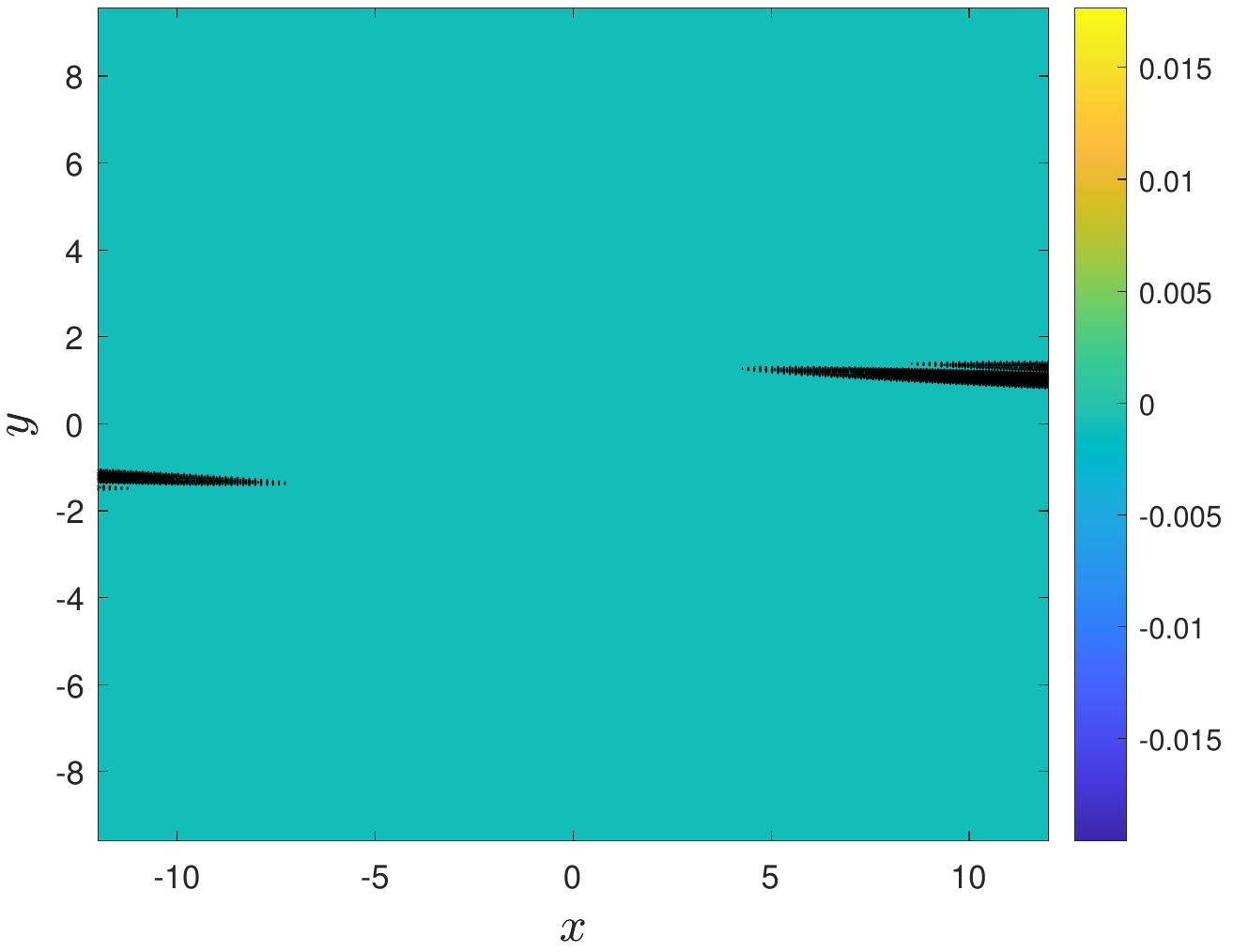}}
{\includegraphics[width=0.48\textwidth,height=0.27\textwidth]{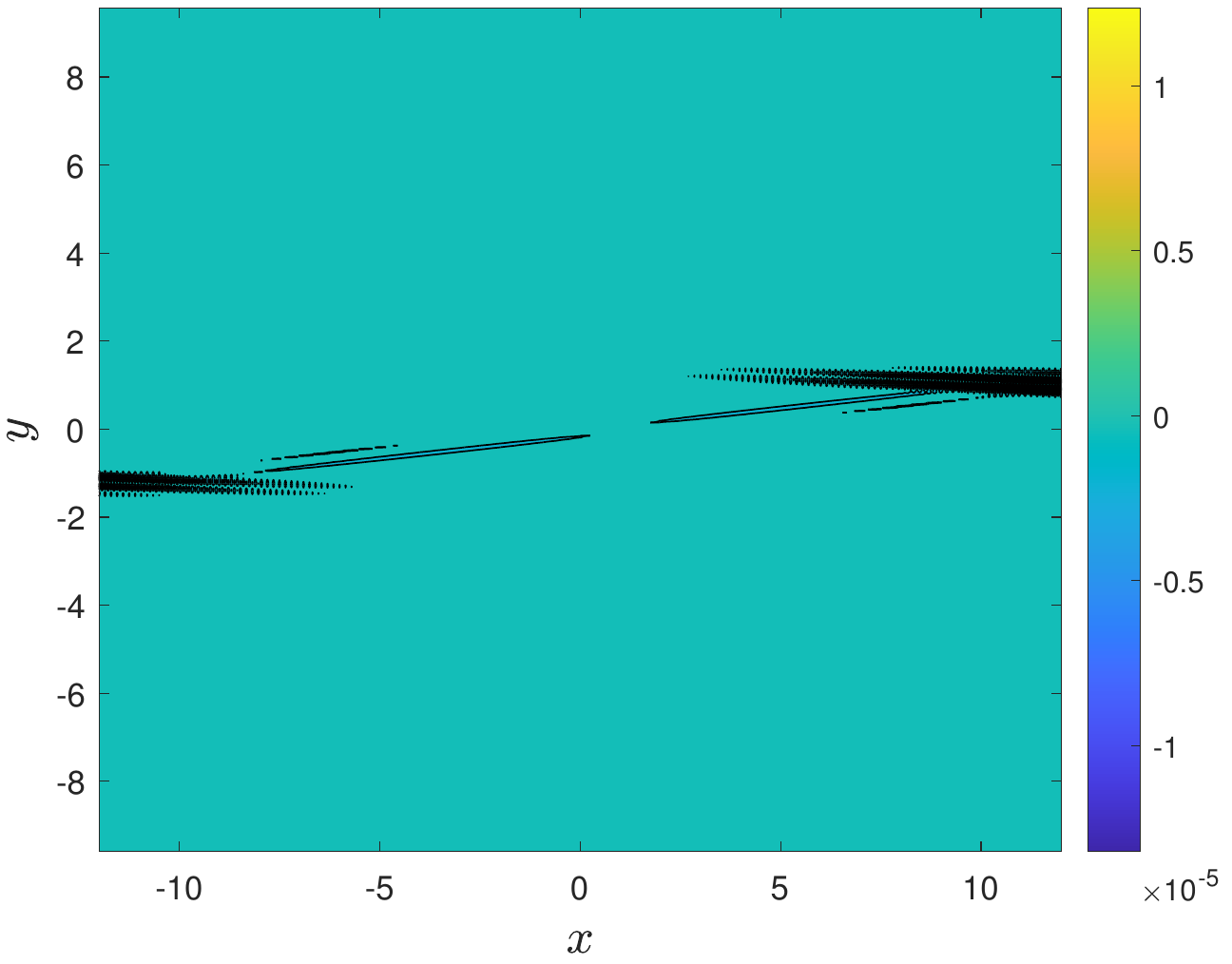}}}
\caption{\small The free advection until $t = 10$ under the Neumann boundary condition (left) or  the natural boundary condition (right) for global cubic spline. Under the Neumann boundary condition, the wave packet tends to be reflected back and leads to an evident accumulation of errors near the boundary. By contrast, the reflection of wave packet can be significantly suppressed under the natural boundary condition. \label{free_compare_Neumann_natural}}
\end{figure}

Now it turns to investigate the influence of different boundary conditions on the global spline. Again, we simulate the free advection in Example \ref{ex3} until $t_{fin} = 10$ under either the natural boundary condition \eqref{natural_boundary} or the Neumann boundary condition $f^{\prime}(x_0) = 0$ and $f^{\prime}(x_N) = 0$ imposed on the global spline.

 As seen in Figure \ref{free_compare_Neumann_natural}, when the Neumann boundary condition is adopted, the wavepacket will be reflected back when it touches the boundary and leads to a rapid accumulation of errors. By contrast, under the natural boundary condition, the reflection of wavepacket is evidently  suppressed and growth rate of errors is dramatically smaller. 
 
 The numerical evidence indicates that it is more appropriate to impose the natural cubic spline to let wavepackets leave the domain without reflecting back.

\section{Comparison between TKM and pseudo-spectral method}

The pseudo-spectral method (PSM for brevity) is a typical way to approximate the $\pdo$ \cite{Ringhofer1990,Goudon2002} 
\begin{equation}\label{def.pdo_convolution}
\Theta_V[f](\bx, \bk, t) = \frac{1}{\mi \hbar (2\pi)^3} \iint_{\mathbb{R}^{6}} \me^{-\mi (\bk - \bk^{\prime}) \cdot \by }D_V(\bx, \by, t) f(\bx, \bk^{\prime}, t) \D \by \D \bk^{\prime}
\end{equation}
with $D_V(\bx, \by, t) = V(\bx + \frac{\by}{2}) - V(\bx - \frac{\by}{2})$. 

Suppose the Wigner function $f(\bx, \bk, t)$ decays outside the finite domain $\mathcal{X} \times [-L_k, L_k]^3$, then one can impose artificial periodic boundary condition in $\bk$-space and use PSM (or the Poisson summation formula)
\begin{equation}
f(\bx, \bk, t) \approx \sum_{\bn \in \mathbb{Z}^3} \widehat{f}_{\bn}(\bx, t) \me^{\frac{2\pi \mi \bn \cdot \bk}{2L_k}}.
\end{equation} 
In addition, starting from the convolution representation of $\pdo$, it yields that
\begin{equation*}
\begin{split}
\Theta_V[f](\bx, \bk, t) & \approx   \frac{1}{\mi \hbar (2\pi)^3} \int_{\mathbb{R}^3} \me^{-\mi \bk \cdot \by} D_V(\bx, \by) \sum_{\bn \in \mathbb{Z}^3} \left(\int_{\mathbb{R}^3} \widehat{f}_{\bn}(\bx, t) \me^{\mi (\frac{\pi}{L_k} \bn -  \by) \cdot \bk^{\prime}} \D \bk^{\prime}\right) \D \by \\
&=  \frac{1}{\mi \hbar (2\pi)^3} \sum_{\bn \in \mathbb{Z}^3} \widehat{f}_{\bn}(\bx, t) \int_{\mathbb{R}^3} \me^{-\mi \bk \cdot \by} D_V(\bx, \by)  \left(\int_{\mathbb{R}^3} \me^{\mi (\frac{\pi }{L_k} \bn - \by) \cdot \bk^{\prime}} \D \bk^{\prime}\right) \D \by \\
&= \frac{1}{\mi \hbar } \sum_{\bn \in \mathbb{Z}^3} \widehat{f}_{\bn}(\bx, t) \int_{\mathbb{R}^3} \me^{-\mi \bk \cdot \by} D_V(\bx, \by) \delta(\frac{\pi \bn}{L_k}  - \by) \D \by \\
&= \frac{1}{\mi \hbar } \sum_{\bn \in \mathbb{Z}^3} \widehat{f}_{\bn}(\bx, t) D_V(\bx, \frac{\pi \bn}{L_k})   \me^{-\frac{\pi\mi}{L_k} \bk \cdot \bn}, 
\end{split}
\end{equation*}
where the third equality uses the Fourier completeness relation. By further truncating $\bn$, we arrive at the approximation formula
\begin{equation}\label{PSM_approximation}
\Theta_V[f](\bx, \bk, t) \approx \frac{1}{\mi \hbar} \sum_{\bn \in \mathcal{I} } \widehat{f}_{\bn}(\bx, t) (V(\bx - \frac{\pi \bn}{2L_k}) - V(\bx + \frac{\pi \bn}{L_k}))  \me^{-\frac{\pi\mi}{2L_k} \bk \cdot \bn},
\end{equation}
where the dual index set $\mathcal{I}$ is that
\begin{equation}
\mathcal{I}\coloneqq \{(n_1, n_2, n_3)\in \mathbb{Z}^3 | n_j = -N_k/2, \dots, N_k/2-1\}.
\end{equation}

However, we would like to report that PSM might fail to produce proper results when $V(\bx)$ has singularities and the formula \eqref{PSM_approximation} is actually not well-defined for $\bx = \pm \frac{\pi \bm{n}}{2 L_k}$. For the sake of comparison, we consider a 6-D problem under the attractive Coulomb potential.
\begin{example}\label{ex_qcs}
Consider a Quantum harmonic oscillator  $V(x) =  -{1}/{|\bx|}$ and a Gaussian wavepacket  adopted as the initial condition.
\begin{equation}
f_0(\bx, \bk) = \pi^{-3} \me^{-\frac{(x_1-1)^2 + x_2^2 + x_3^2}{2} - 2k_1^2 - 2k_2^2 - 2k_3^2}.
\end{equation}
\end{example}

\begin{figure}[!h]
\centering
\subfigure[TKM.]{

{\includegraphics[width=0.32\textwidth,height=0.18\textwidth]{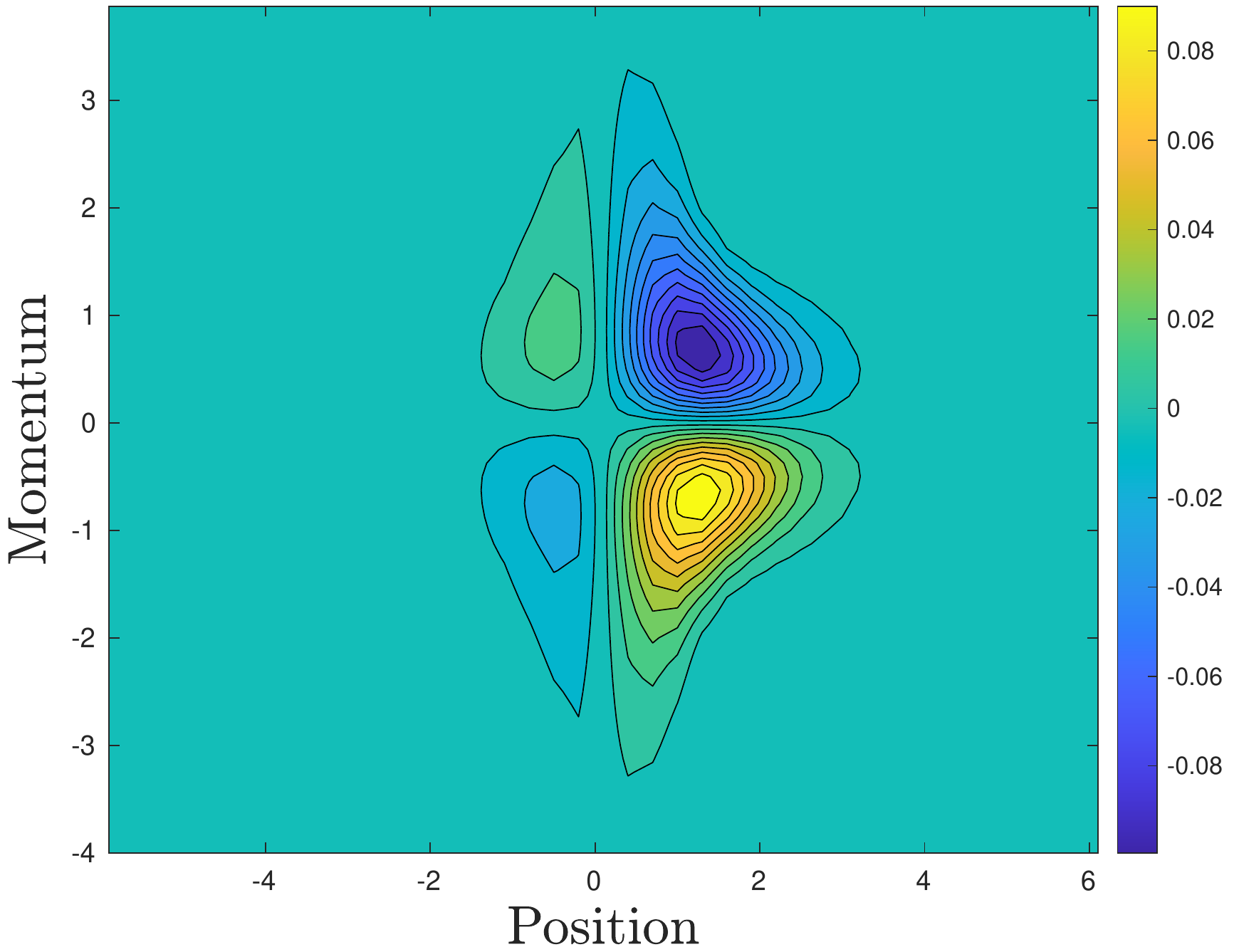}}}
\subfigure[PSM ($\delta x= 0$).]{
{\includegraphics[width=0.32\textwidth,height=0.18\textwidth]{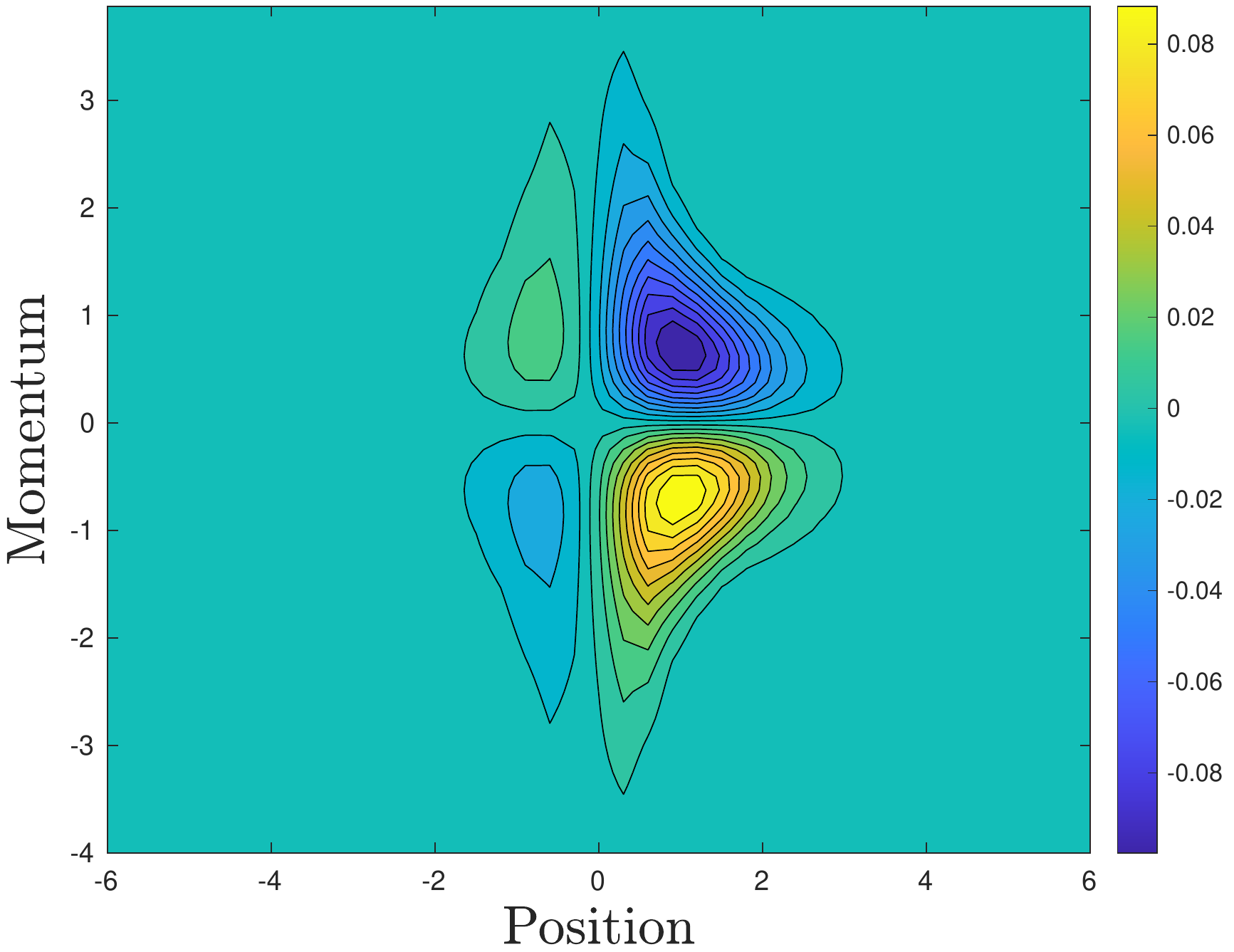}}}
\subfigure[PSM ($\delta x = 0.01$).]{
{\includegraphics[width=0.32\textwidth,height=0.18\textwidth]{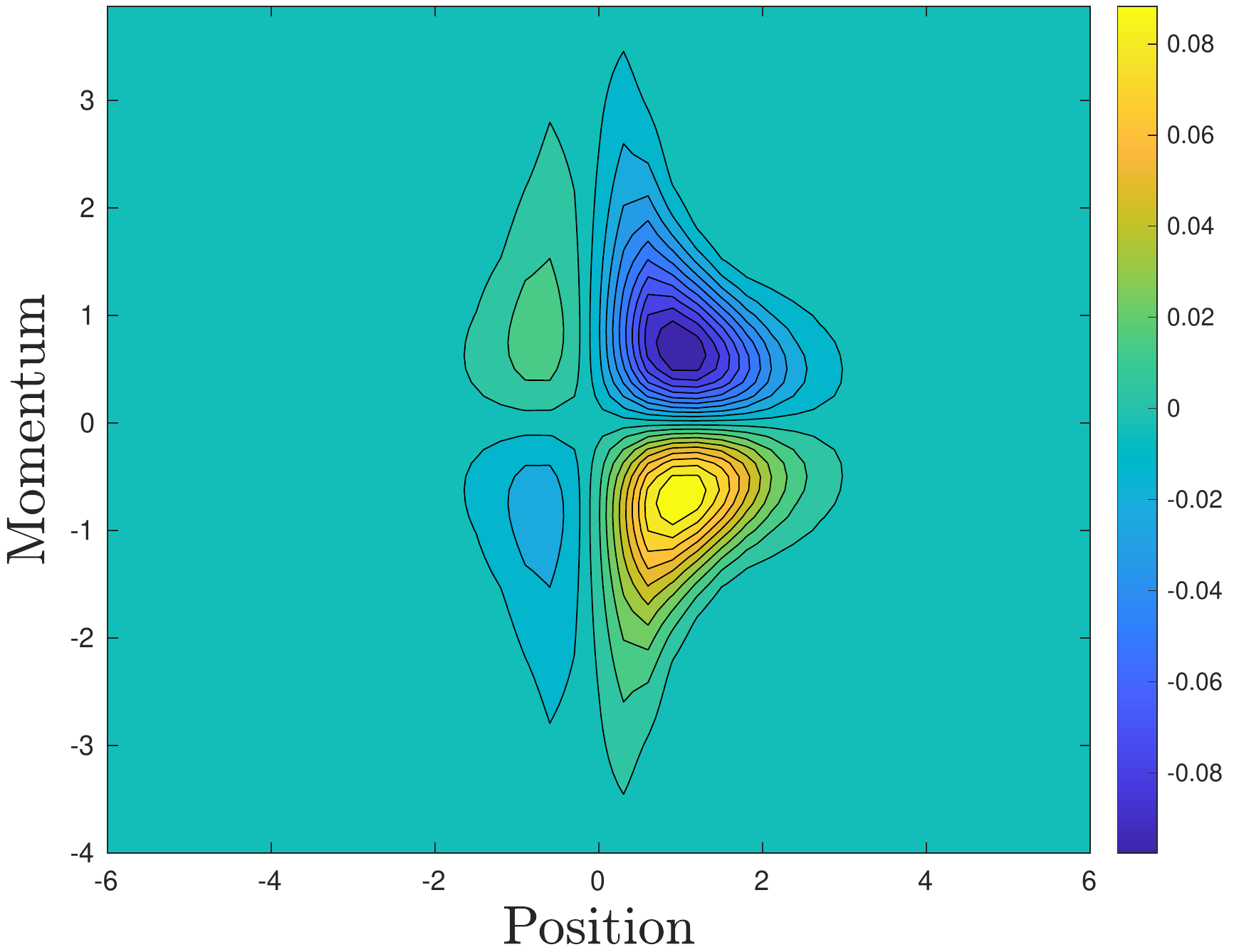}}}
\caption{\small  A comparison between TKM and PSM for a $\pdo$ with an initial Gaussian wavepacket.}
 \label{initial_PDO_comp}
\end{figure}

We first calculate $\pdo$ under TKM or PSM for a $\mathcal{X}$-grid mesh $[-6, 6]^3$ with $N_x = 41, \Delta x = 0.3$ and $\mathcal{K}$-grid mesh $[-4, 4]^3$ with $N_k = 64, \Delta k = 0.125$. In order to get rid of the blow-up in the formula \eqref{PSM_approximation}, we try to adopt two ways. The first is to shift $\mathcal{X}$-grid mesh to $[-6 + \delta x, 6 + \delta x]^3$ with a small spacing $\delta x$. The second is to set $\delta x = 0$ and let $V(\bx) = 0$ when $|\bx| = 0$. A comparison among the initial $\pdo$ under different strategies is given in Figure \ref{initial_PDO_comp}. At first glance, no evident differences are observed in the numerical results under TKM or PSM. 

\begin{figure}[!h]
\centering
\includegraphics[width=0.48\textwidth,height=0.27\textwidth]{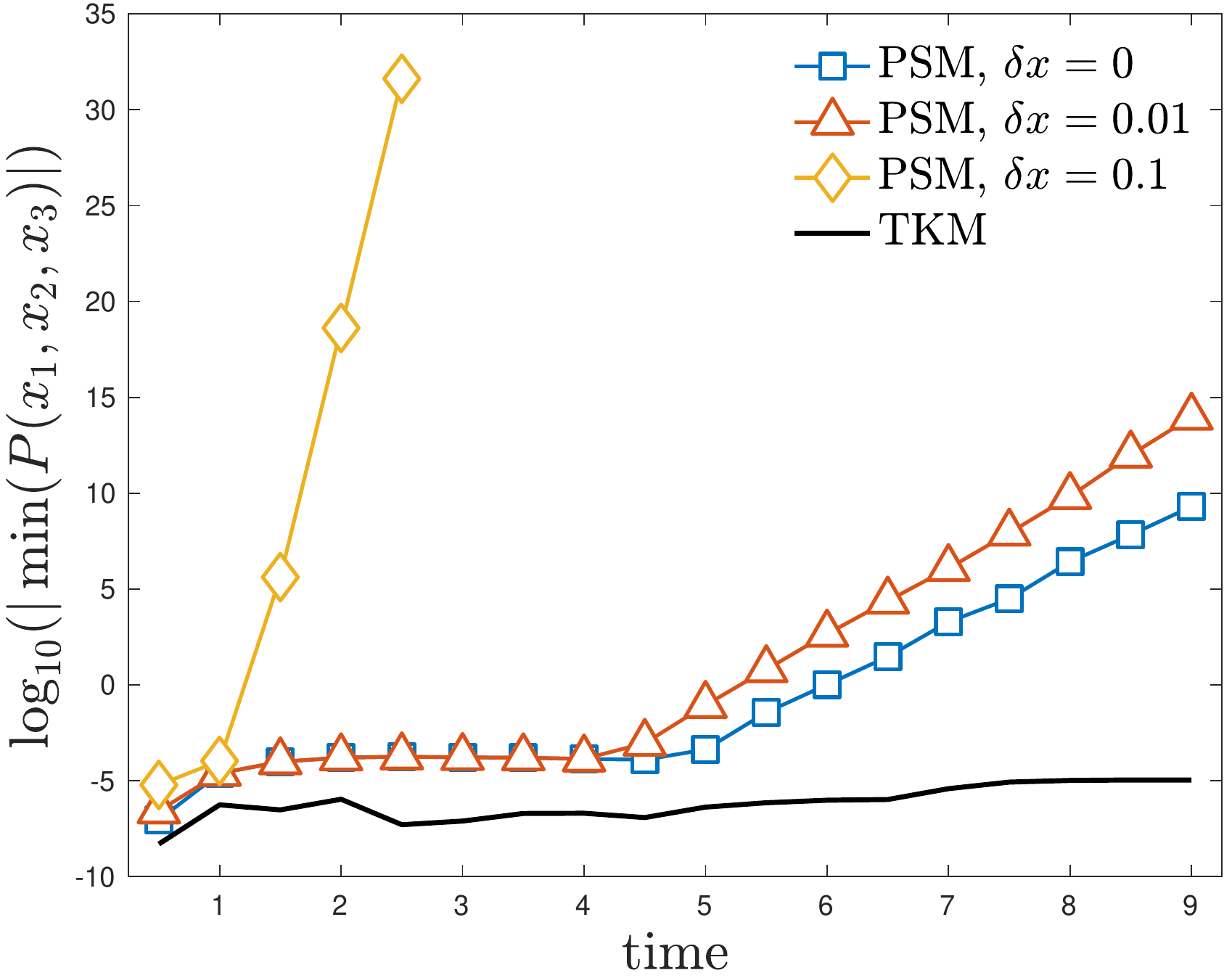}
\caption{\small The time evolution of numerical $\min_{\mathcal{X}} P(x_1, x_2, x_3)$ in log10 scale. PSM may introduce very large artificial negative parts and finally suffers from numerical instability. \label{PSM_instability}}
\end{figure}

However, when simulating the Wigner dynamics with $\mathcal{X}\times \mathcal{K} = [-6, 6]^3\times [-4, 4]^3$ and $N_x = 41, \Delta x = 0.3$, $N_k = 32, \Delta k = 0.25$, we have found that TKM and PSM exhibit distinct performances. Specifically, PSM may suffer from large errors near singularity and numerical instability as it treats the singularity near the origin incorrectly. For the sake of illustration, we consider the spatial marginal density
\begin{equation}
P(x_1, x_2, x_3) = \iiint_{\mathbb{R}^3} f(\bx, k_1, k_2, k_3) \D k_1 \D k_2 \D k_3.
\end{equation}
The spatial marginal density is proved to be positive semi-definite. Therefore, the negative value of numerical solution can be used an indicator for accuracy and stability and is visualized in Figure \ref{PSM_instability}. Although the spectral method might not preserve the positivity of the spatial marginal density, the errors remain at a stable level when TKM is adopted. By contrast, PSM may introduce very large artificial negative parts and finally results in numerical instability.

\begin{figure}[!h]
\centering
\subfigure[$t = 2$a.u. (left: TKM, middle: PSM with $\delta x = 0$, right: PSM with $\delta x = 0.01$).]{
{\includegraphics[width=0.32\textwidth,height=0.18\textwidth]{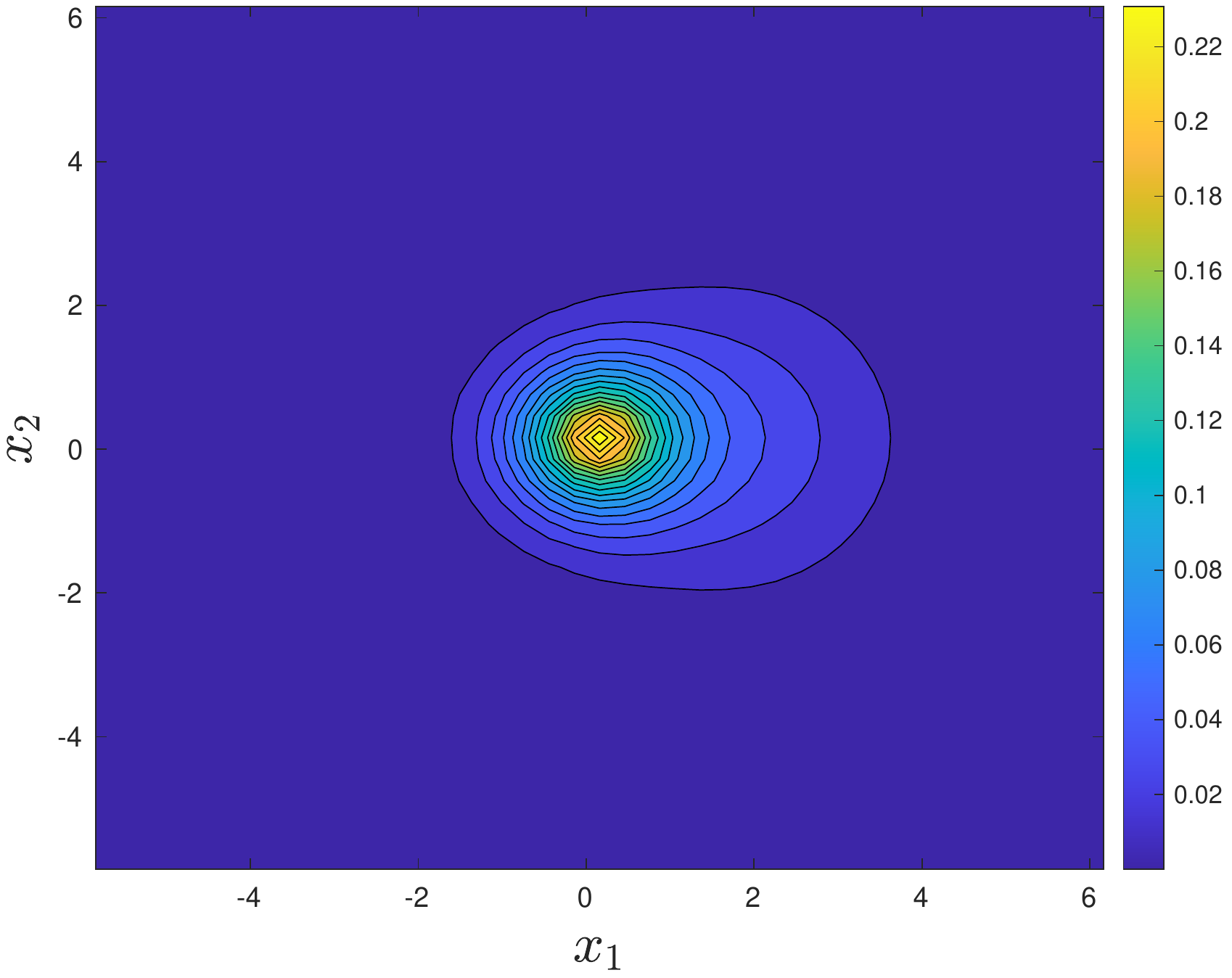}}
{\includegraphics[width=0.32\textwidth,height=0.18\textwidth]{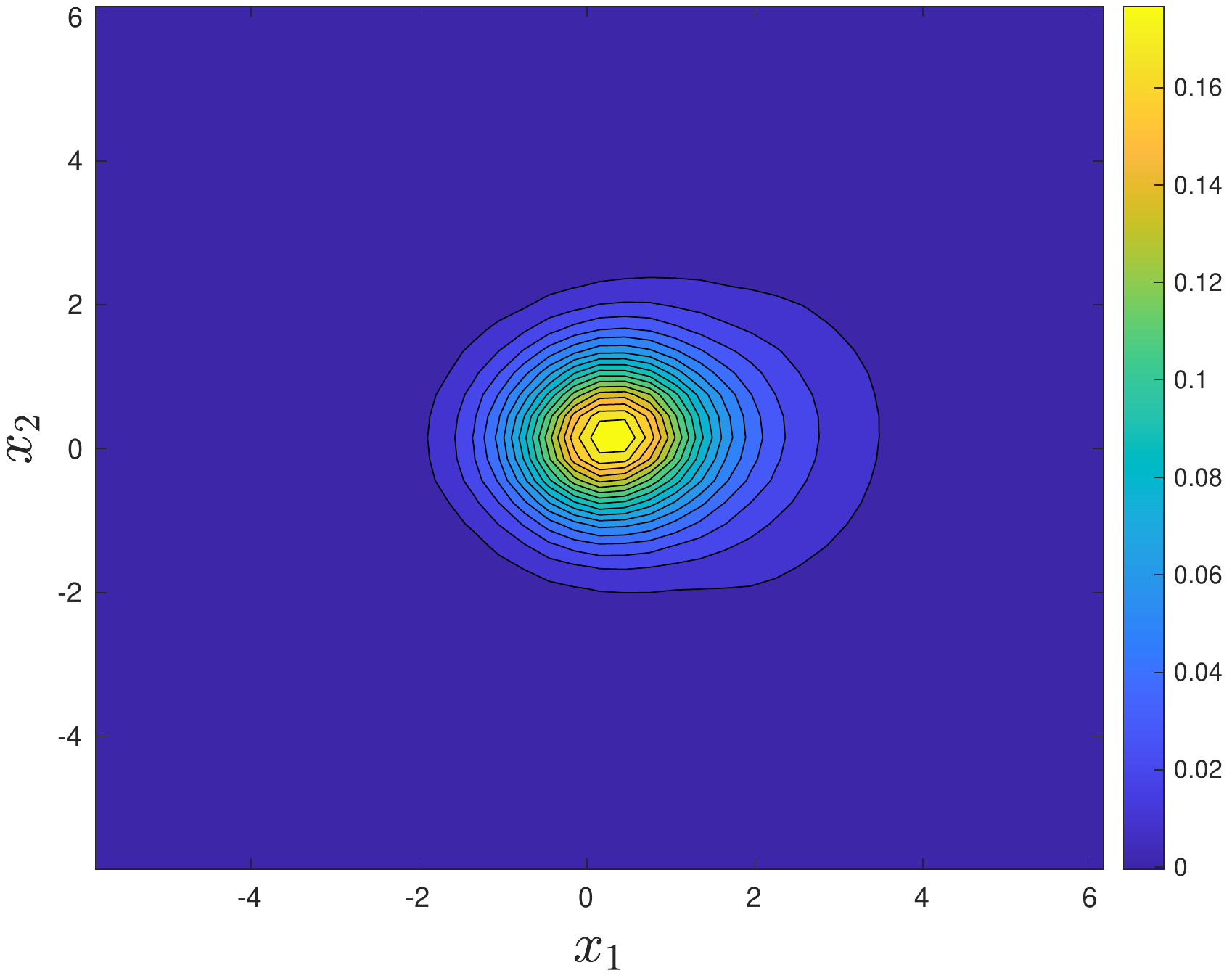}}
{\includegraphics[width=0.32\textwidth,height=0.18\textwidth]{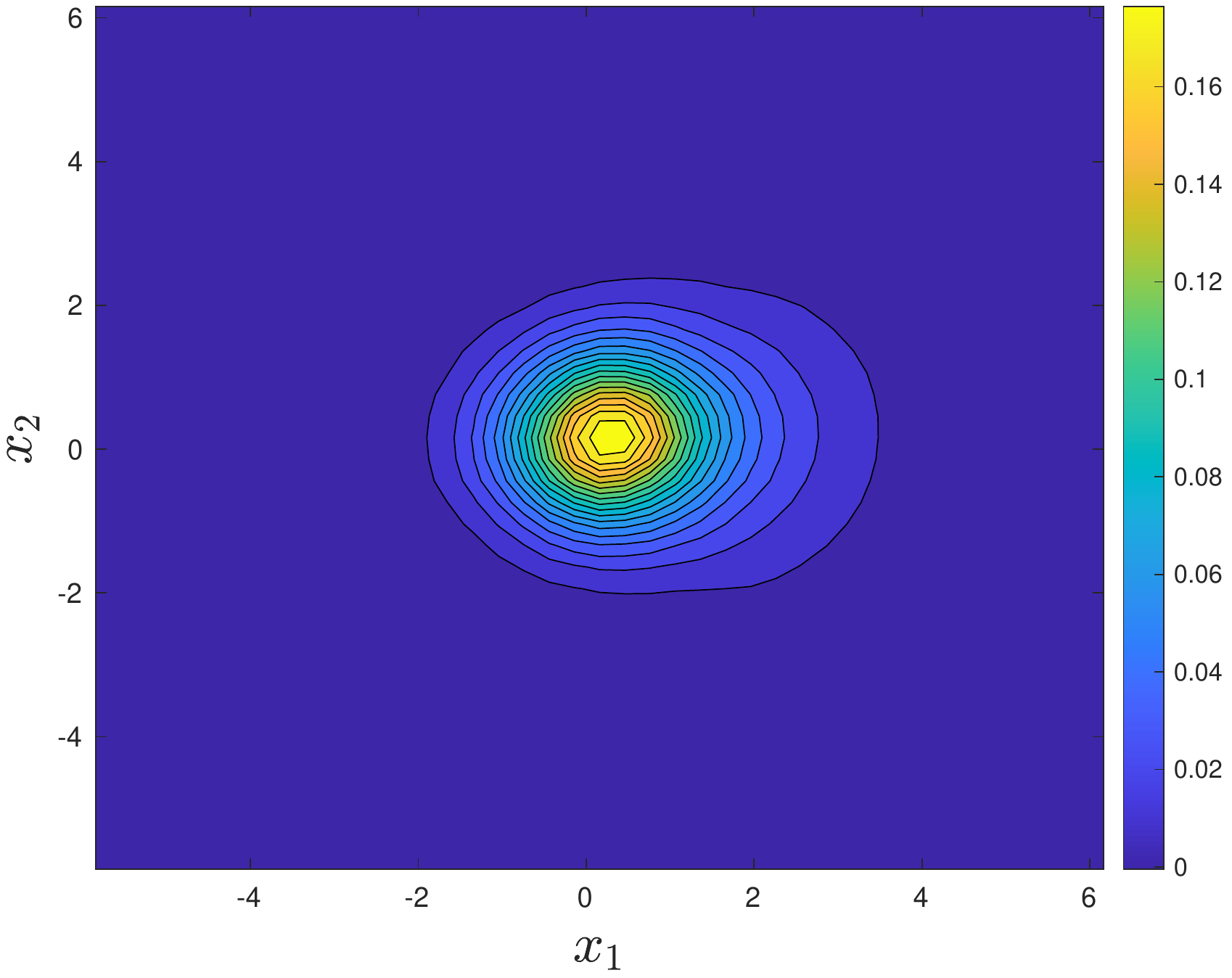}}}
\\
\centering
\subfigure[$t = 4$a.u. (left: TKM, middle: PSM with $\delta x = 0$, right: PSM with $\delta x = 0.01$).]{
{\includegraphics[width=0.32\textwidth,height=0.18\textwidth]{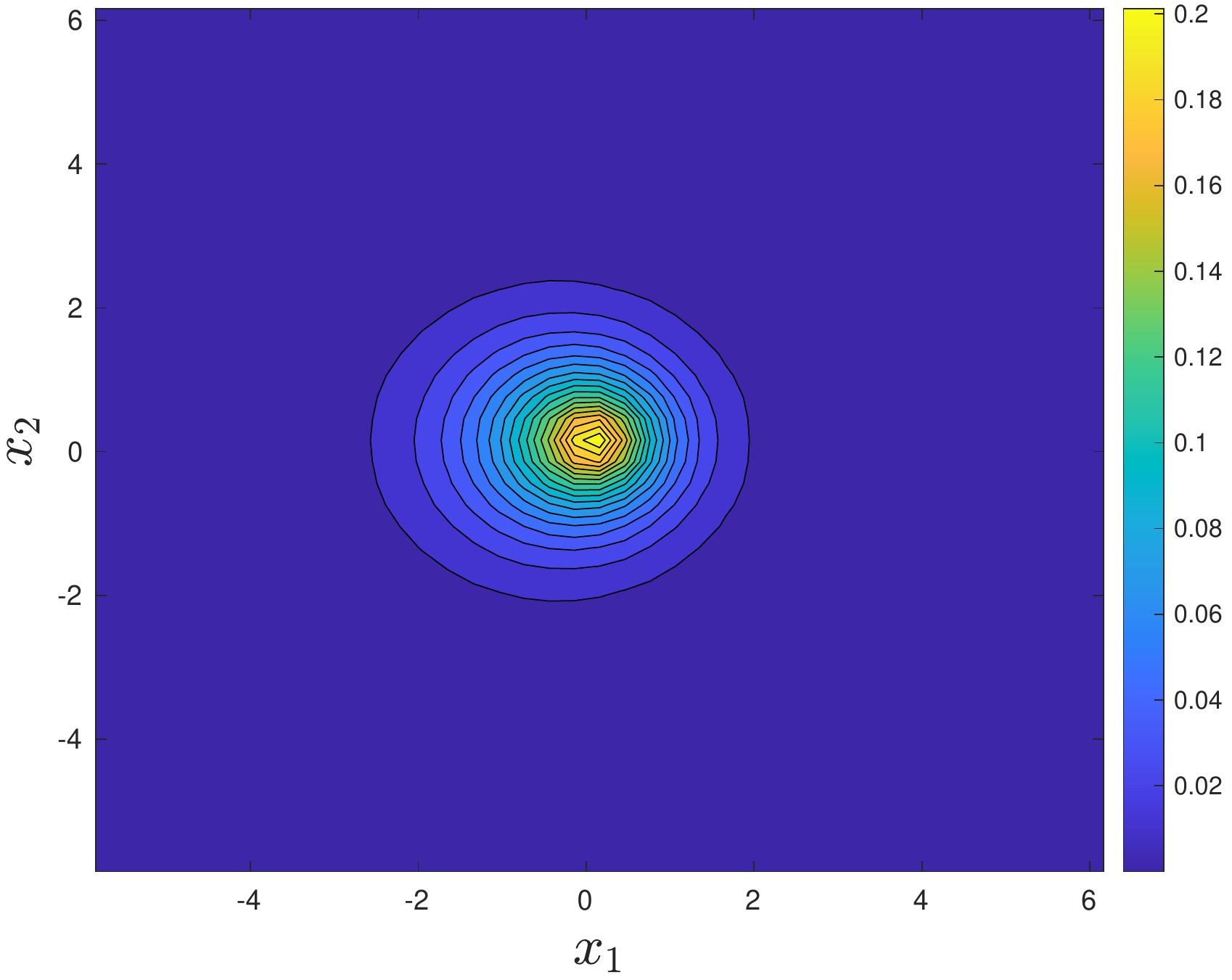}}
{\includegraphics[width=0.32\textwidth,height=0.18\textwidth]{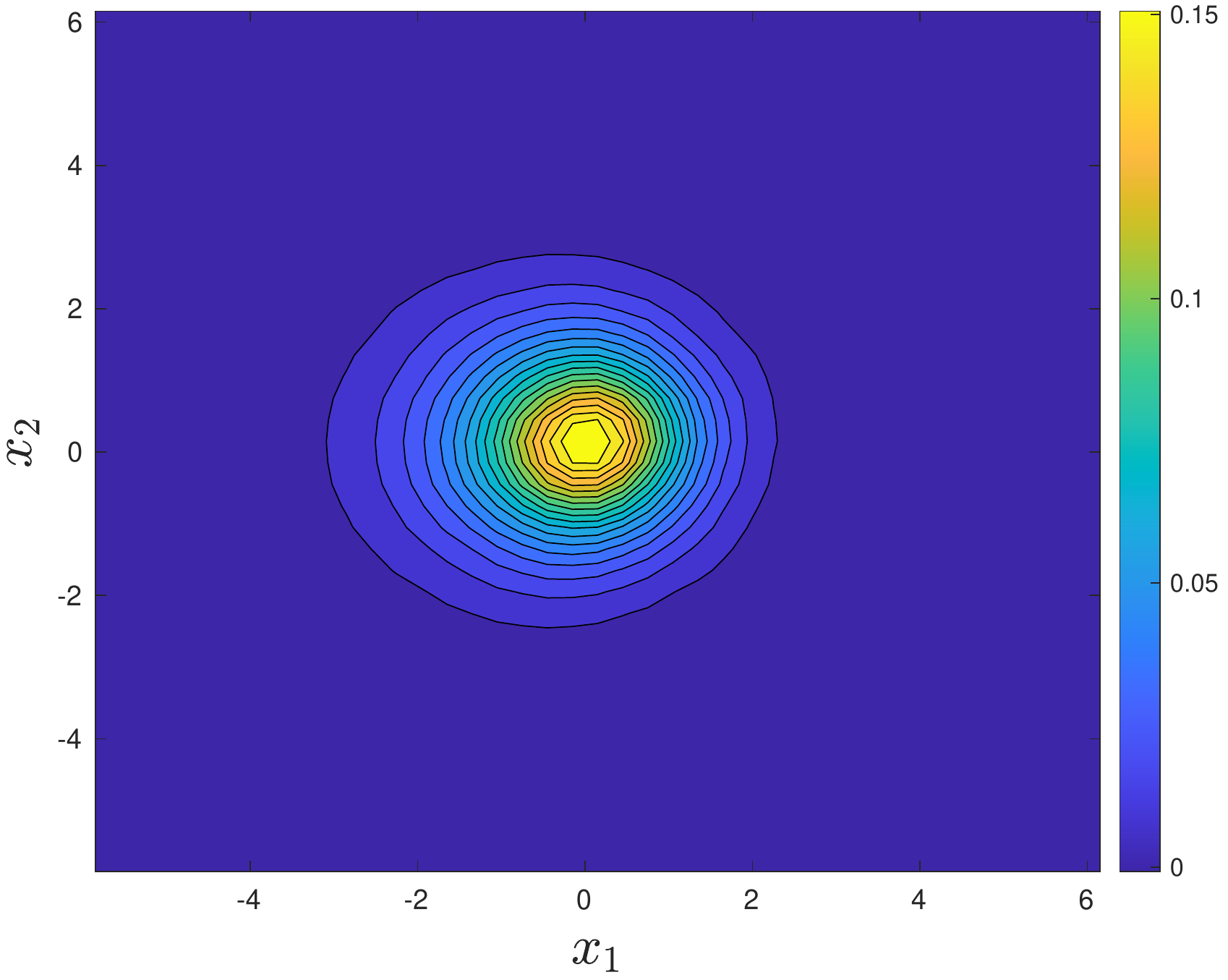}}
{\includegraphics[width=0.32\textwidth,height=0.18\textwidth]{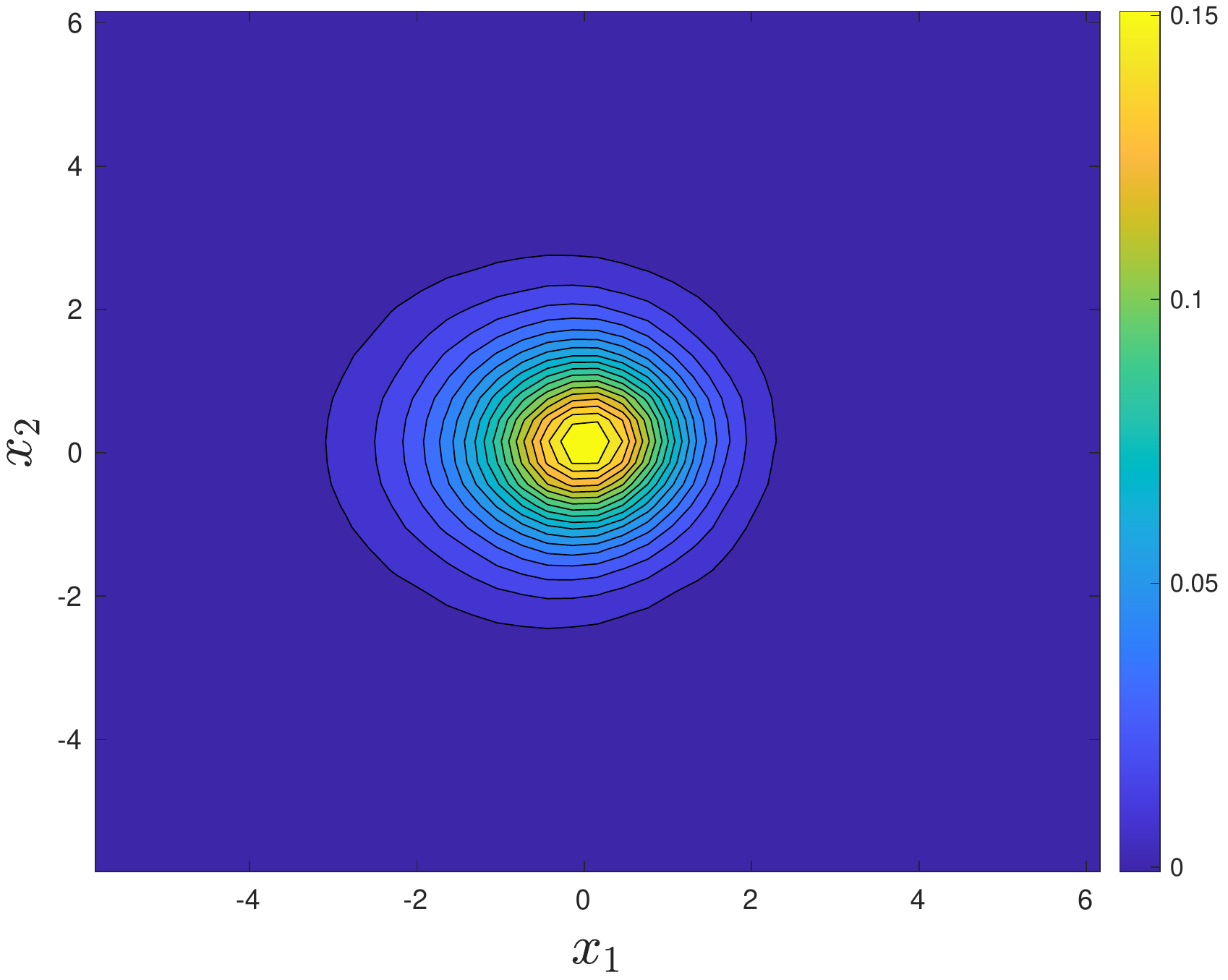}}}
\\
\centering
\subfigure[$t = 6$a.u. (left: TKM, middle: PSM with $\delta x = 0$, right: PSM with $\delta x = 0.01$).]{
{\includegraphics[width=0.32\textwidth,height=0.18\textwidth]{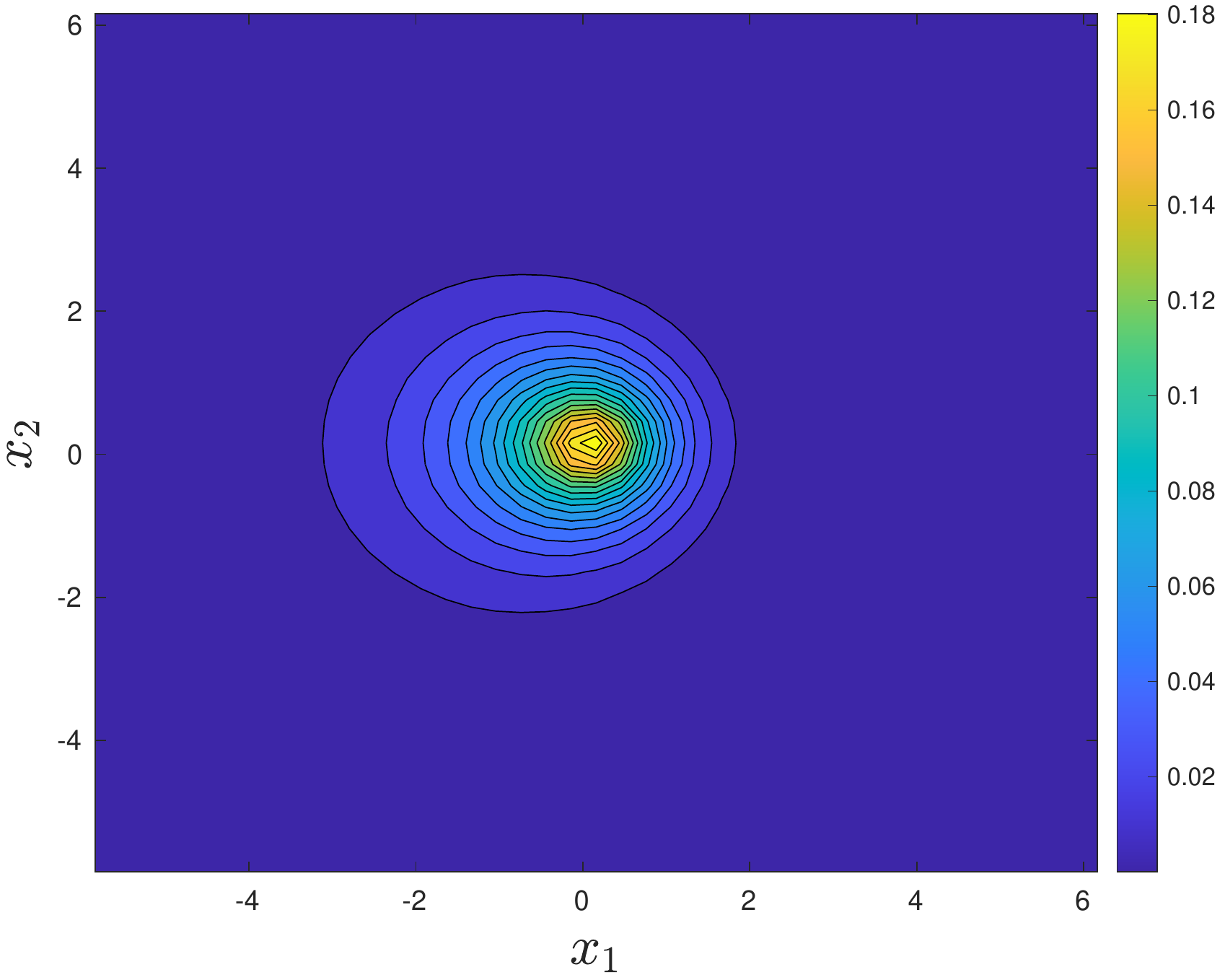}}
{\includegraphics[width=0.32\textwidth,height=0.18\textwidth]{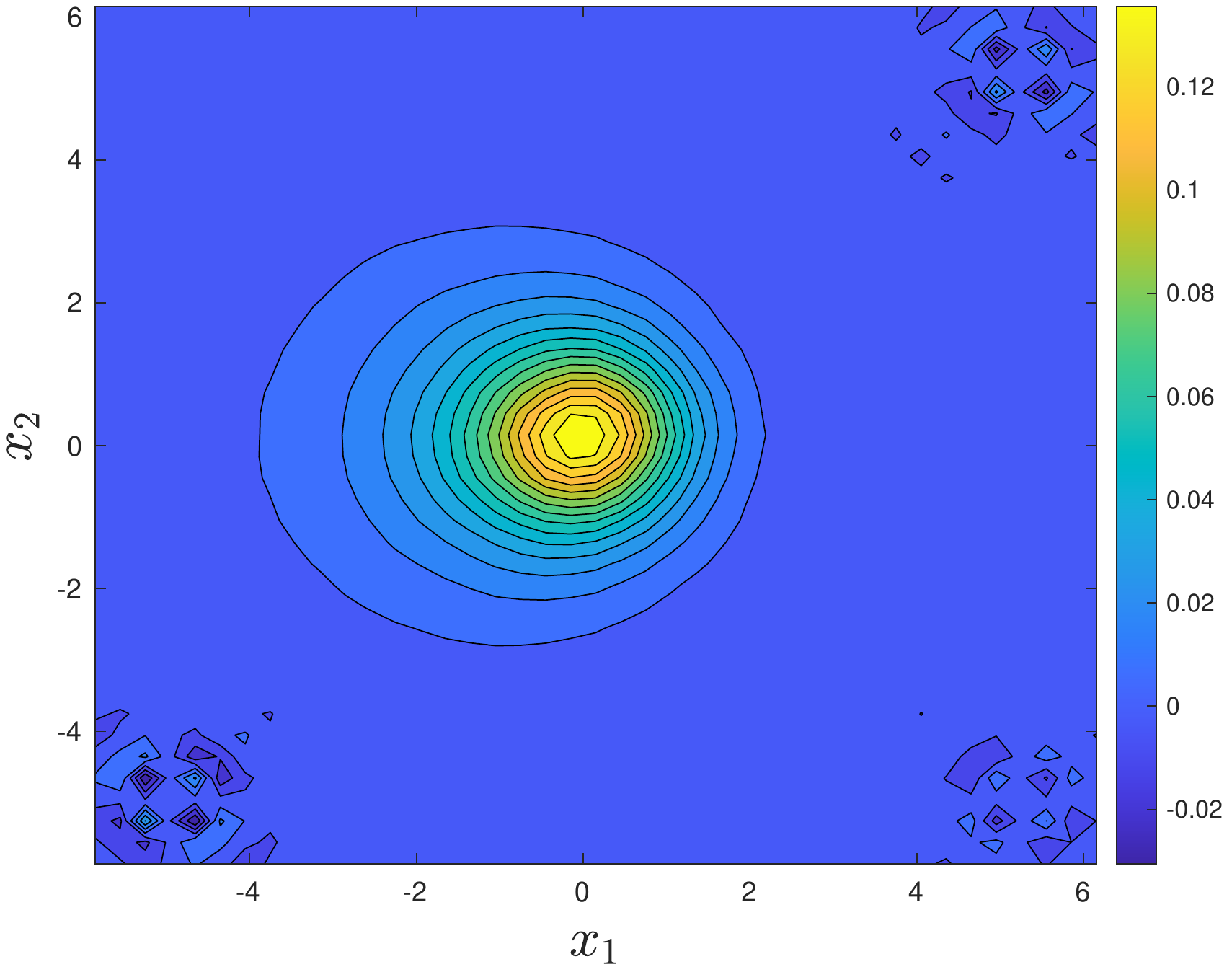}}
{\includegraphics[width=0.32\textwidth,height=0.18\textwidth]{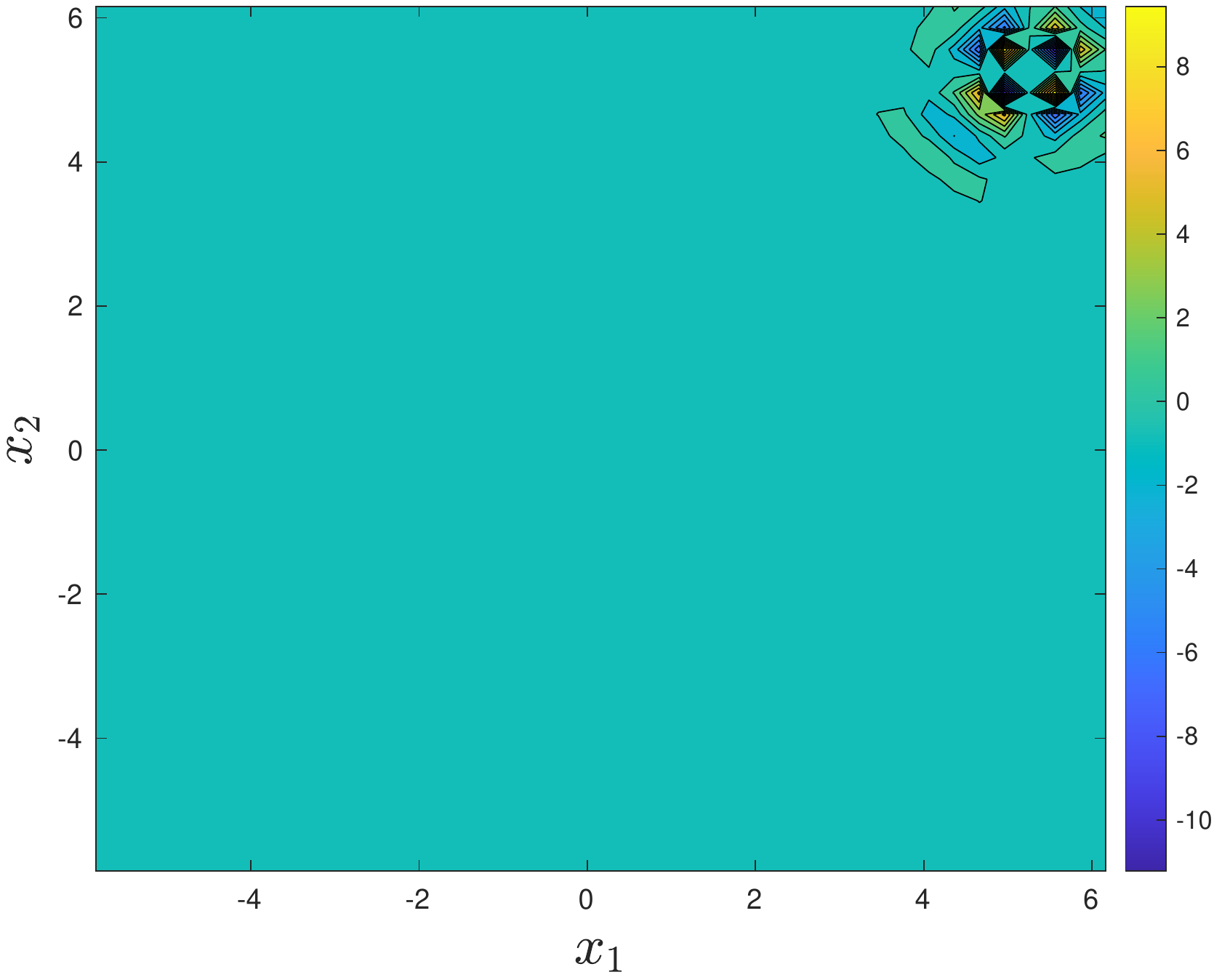}}}
\caption{\small  Visualization of the spatial marginal distribution projected onto $(x_1$-$x_2$) plane under TKM and PSM. PSM might not produce correct numerical results and suffers from  instability. \label{comp_TKM_PSM}}
\end{figure}

In Figure \ref{comp_TKM_PSM}, we visualize the spatial marginal distribution projected onto $(x_1$-$x_2$) plane. It is found that the peak of spatial marginal distribution has been evidently smoothed out by PSM at $2$ a.u. and artificial negative valleys are clearly seen at $6$a.u. This coincides with the observation in Figure \ref{PSM_instability} that PSM suffers from instability soon after $6$a.u.

\section{Comparison among exponential integrators and splitting method}

Finally, it needs to make a thorough comparison among various integrators, which in turn provides a guiding principle in choosing an appropriate integrator for our 6-D simulations. To this end, we provide two examples with exact solution. The first is the quantum harmonic oscillator in 2-D phase space and the second is the Hydrogen Wigner function of 1s state in 6-D phase space (see Example \ref{ex_qcs}).

The performance metrics include the $L^2$-error $\varepsilon_{2}(t)$:
\begin{align}
\varepsilon_{2}(t) &= \left[\iint_{\mathcal{X}\times \mathcal{K}} \left(f^{\textup{ref}}\left(\bm{x},\bm{k},t\right)-f^{\textup{num}}\left(\bm{x},\bm{k},t\right)\right)^{2}\textup{d}\bm{x}\textup{d} \bm{k}\right]^{\frac{1}{2}},\label{eq:e2}
\end{align}
the maximal error $\varepsilon_{\infty}(t)$:
\begin{align}
\varepsilon_{\infty}(t) &=\max_{(\bm{x},\bm{k})\in\mathcal{X}\times \mathcal{K}}\big |f^{\textup{ref}}\left(\bm{x},\bm{k},t\right)-f^{\textup{num}}\left(\bm{x},\bm{k},t\right) \big |, \label{eq:ef}
\end{align}
and the deviation of total mass $\varepsilon_{\textup{mass}}(t)$:
\begin{align}
\varepsilon_{\textup{mass}}(t) &= \Big |\iint_{\mathcal{X}\times \mathcal{K}} f^{\textup{num}}\left(\bm{x},\bm{k},t\right)\textup{d}\bm{x}\textup{d} \bm{k}-\iint_{\Omega} f^{\textup{ref}}\left(\bm{x},\bm{k},t=0\right)\textup{d}\bm{x}\textup{d} \bm{k} \Big |,\label{eq:emass}
\end{align}
where $f^{\textup{ref}}$ and $f^{\textup{num}}$ denote the reference and numerical solution, respectively, and $\mathcal{X}\times \mathcal{K}$ is the computational domain. In practice, the integral can be replaced by the average over all grid points. Besides, the relative maximal error and relative $L^2$-error are obtained by $\frac{\varepsilon_{\infty}(t)}{\max(|f(\bx, \bk, 0)|)}$ and $\varepsilon_{2}(t)/ \sqrt{\iint(|f(\bx, \bk, 0)|^2 \D \bx \D \bk)}$, respectively.

Our main observations are summarized as follows. 

\begin{itemize}

\item[1.] In order to ensure the accuracy of temporal integration, it is recommended to use LPC1, instead of splitting scheme or multi-stage schemes.

\item[2.] The operator splitting scheme is still useful in practice, as it saves half of the cost in calculation of nonlocal terms. 

\item[3.] It is suggested to choose the stencil length $n_{nb} = 15$ for PMBC to maintain the accuracy, while $n_{nb} < 10$ might lead to an evident loss of total mass.

\end{itemize}

The one-stage Lawson predictor-corrector scheme exhibit the best performance. Actually, the advantage of the Lawson scheme in both accuracy and stability has also been reported in the Boltzmann community \cite{CrouseillesEinkemmerMassot2020} recently.

\subsection{Quantum harmonic oscillator in 2-D phase space}

The third example is the quantum harmonic oscillator  $V(x) = \frac{1}{2} m \omega x^2$. In this situation, $\pdo$ reduces to the first-order derivative,
\begin{equation}\label{Wigner_harmonic}
\frac{\partial }{\partial t} f(x, k, t) + \frac{\hbar k}{m}  \nabla_{x} f(x, k, t) - \frac{1}{\hbar} \nabla_{x} V(x)  \nabla_{k} f(x, k, t) = 0.
\end{equation}
The exact solution can be solved by $f(x, k, t) = f(x(t), k(t), 0)$,
where $(x(t), k(t))$ obey a (reverse-time) Hamiltonian system ${\partial x}/{\partial t} =  -{\hbar k}/{m}, {\partial k}/{\partial t} = {m\omega x}/{\hbar}$,
and has the following form
\begin{equation}
\begin{split}
&x(t) = \cos \left(\sqrt{\omega} t\right) x(0) - \frac{\hbar}{m \sqrt{\omega}} \sin \left(\sqrt{ \omega}t\right) k(0), \\
&k(t) =\frac{m\sqrt{\omega}}{\hbar}\sin \left(\sqrt{\omega} t\right) x(0) +  \cos \left(\sqrt{ \omega}t\right) k(0).
\end{split}
\end{equation}

\begin{example} 
Consider a Quantum harmonic oscillator  $V(x) =  \frac{m \omega x^2}{2}$ and a Gaussian wavepacket $f_0(x, k) = \pi^{-1} \me^{-\frac{(x-1)^2}{2} - 2k^2}$ adopted as the initial condition. Here we choose and $\omega = (\pi/5)^2$ so that the wavepacket returns back to the initial state at the final time $T = 10$.
\end{example}

\begin{figure}[!h]
\centering
\subfigure[$\Delta x = 0.3$. \label{harmonic_comp_dx_03}]{
{\includegraphics[width=0.32\textwidth,height=0.22\textwidth]{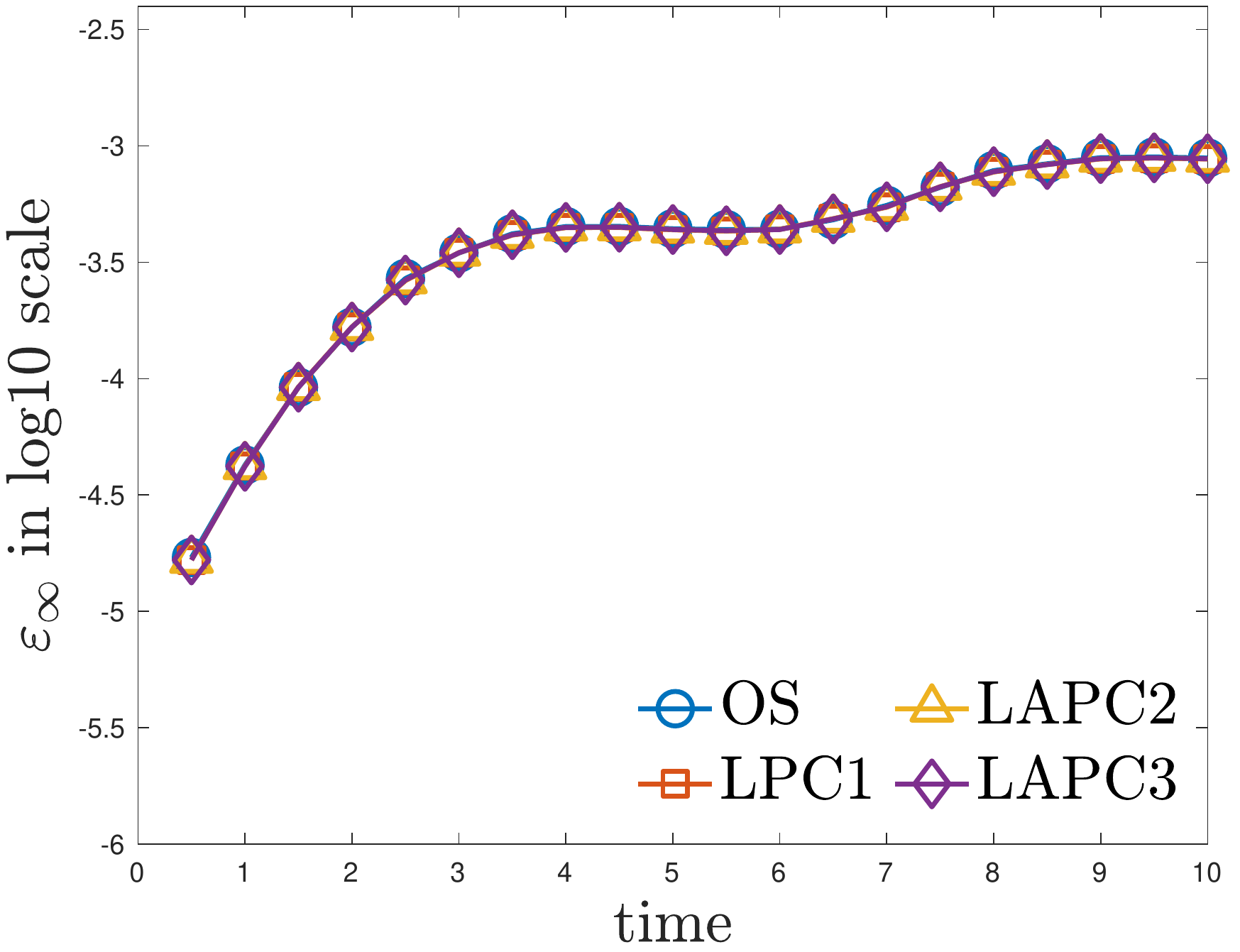}}
{\includegraphics[width=0.32\textwidth,height=0.22\textwidth]{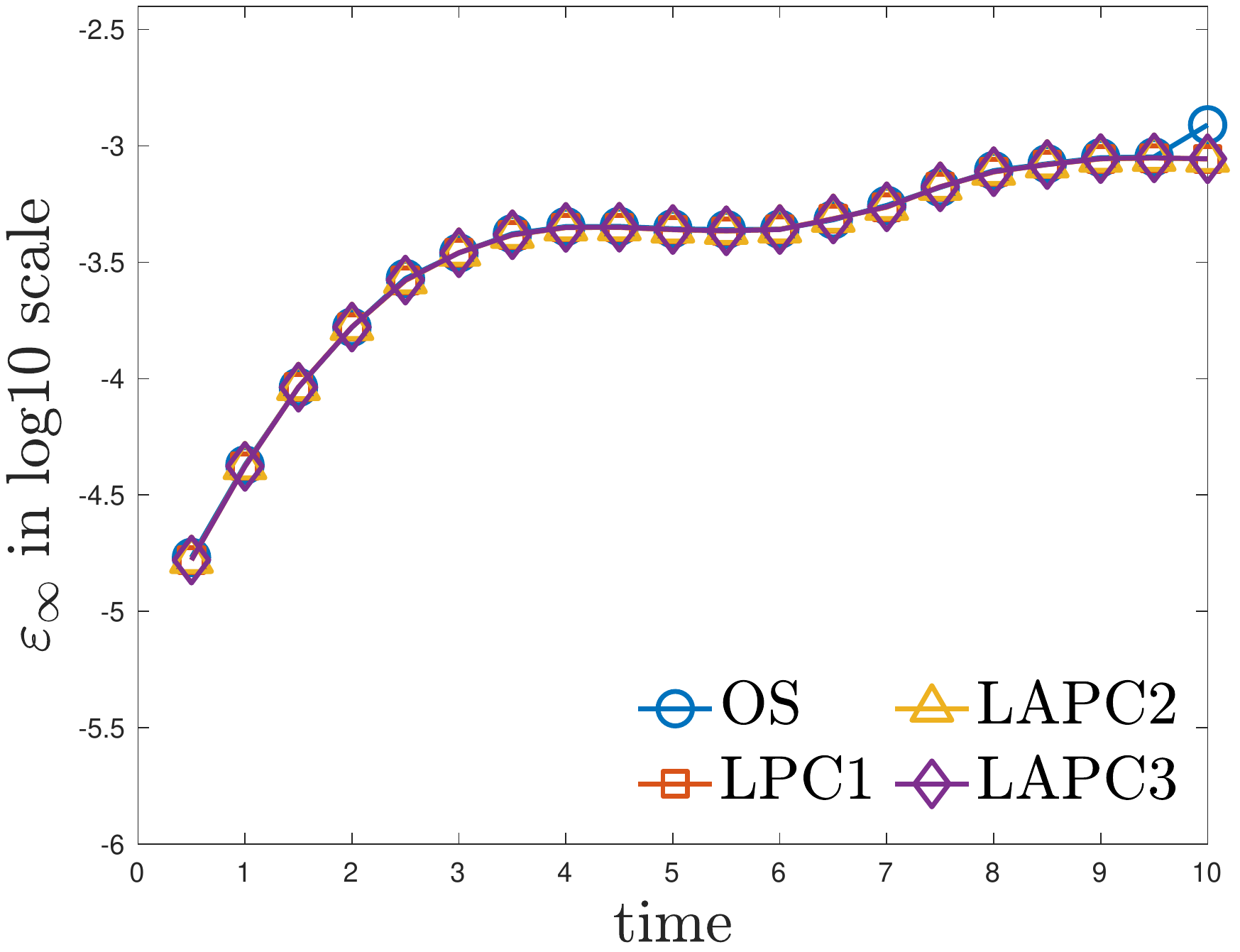}}
{\includegraphics[width=0.32\textwidth,height=0.22\textwidth]{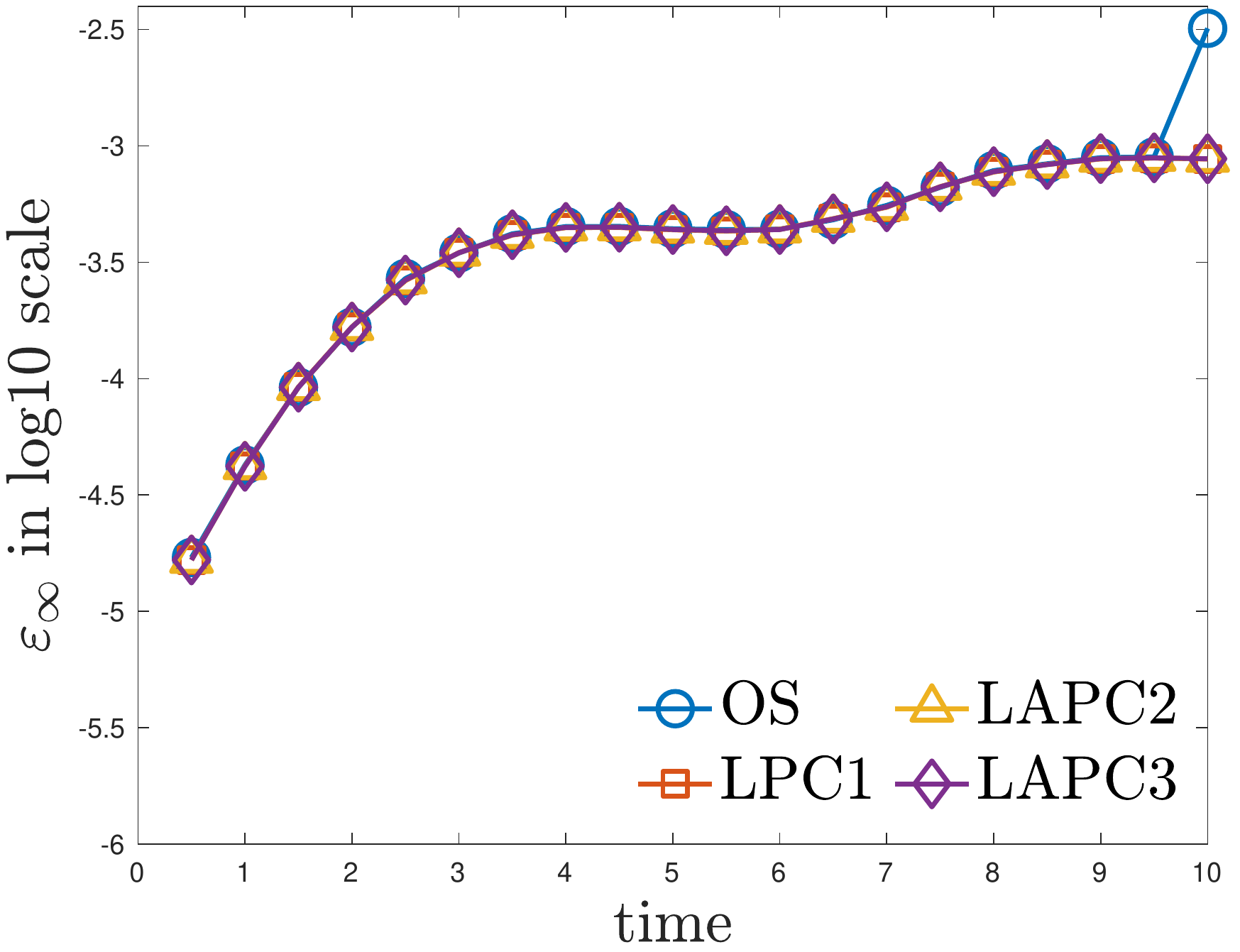}}}
\centering
\subfigure[$\Delta x = 0.2$. \label{harmonic_comp_dx_02}]{
{\includegraphics[width=0.32\textwidth,height=0.22\textwidth]{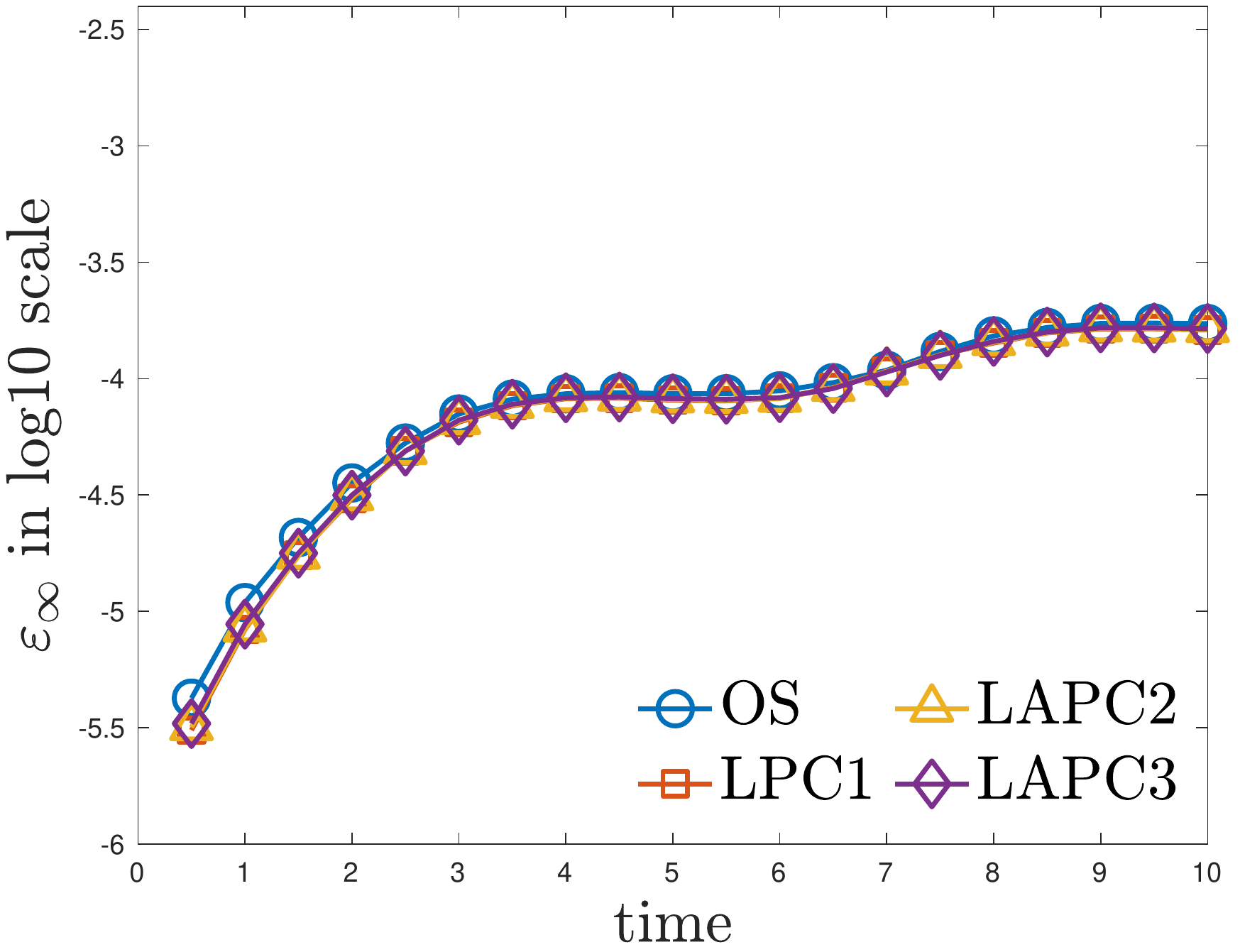}}
{\includegraphics[width=0.32\textwidth,height=0.22\textwidth]{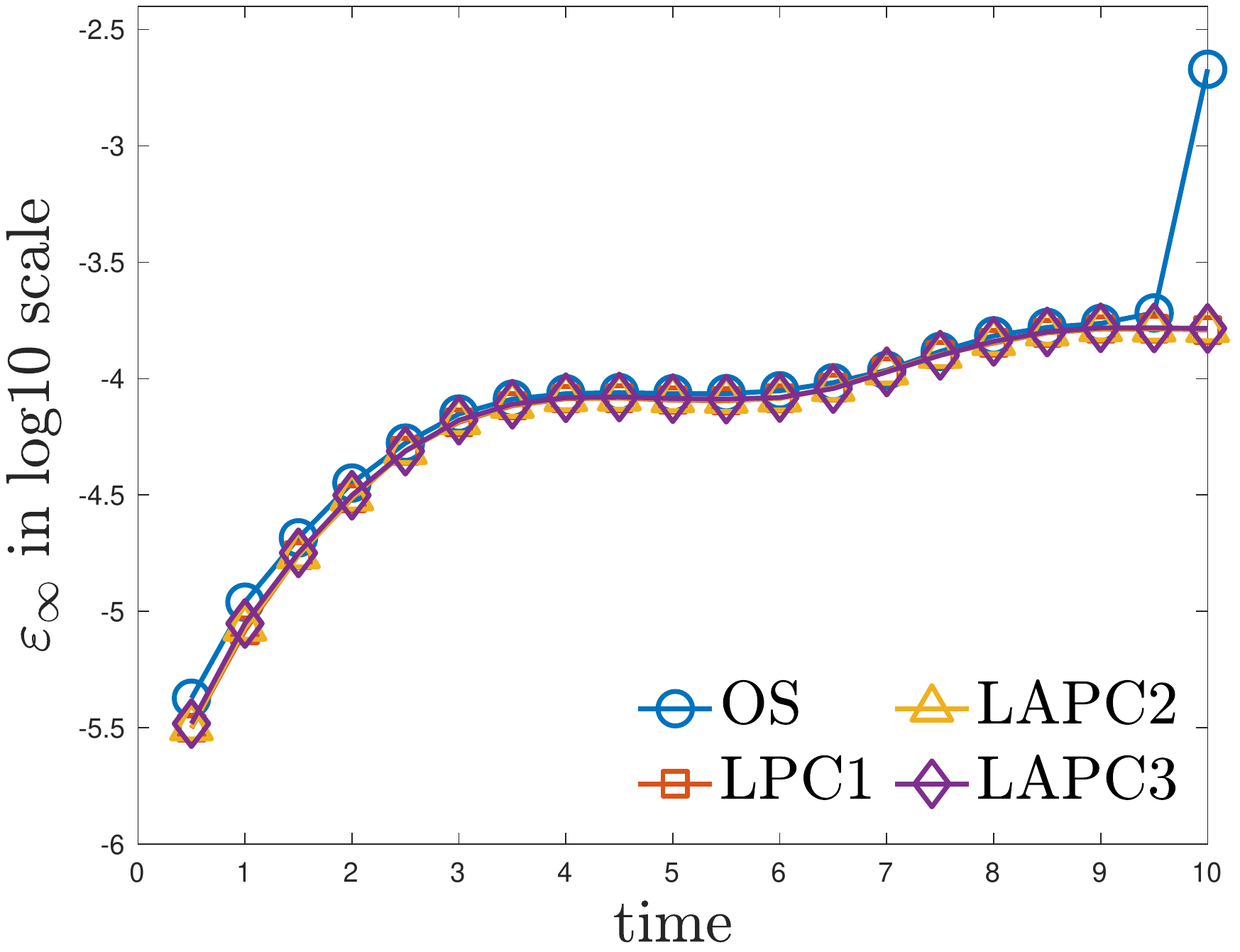}}
{\includegraphics[width=0.32\textwidth,height=0.22\textwidth]{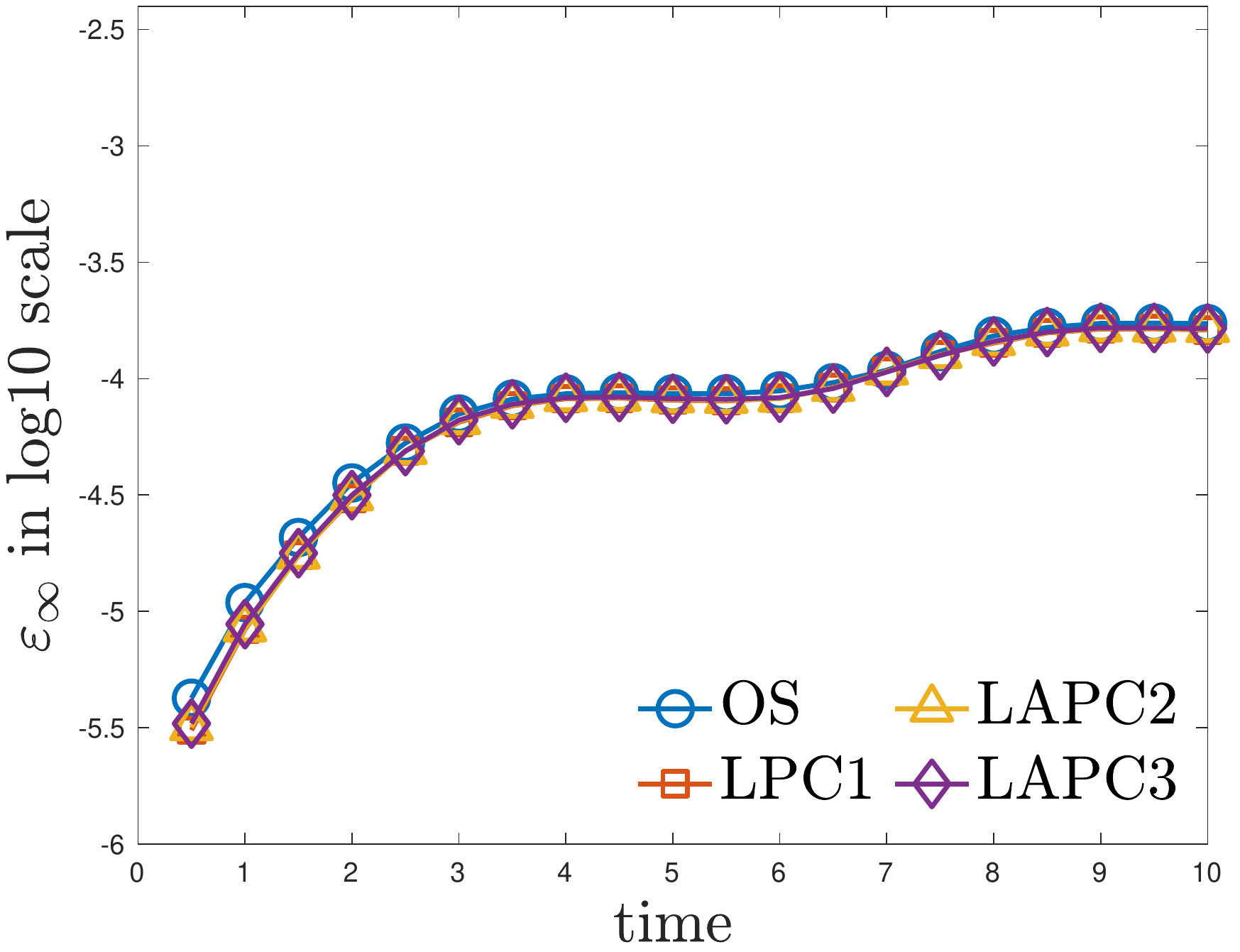}}}
\\
\centering
\subfigure[$\Delta x = 0.1$.]{
{\includegraphics[width=0.32\textwidth,height=0.22\textwidth]{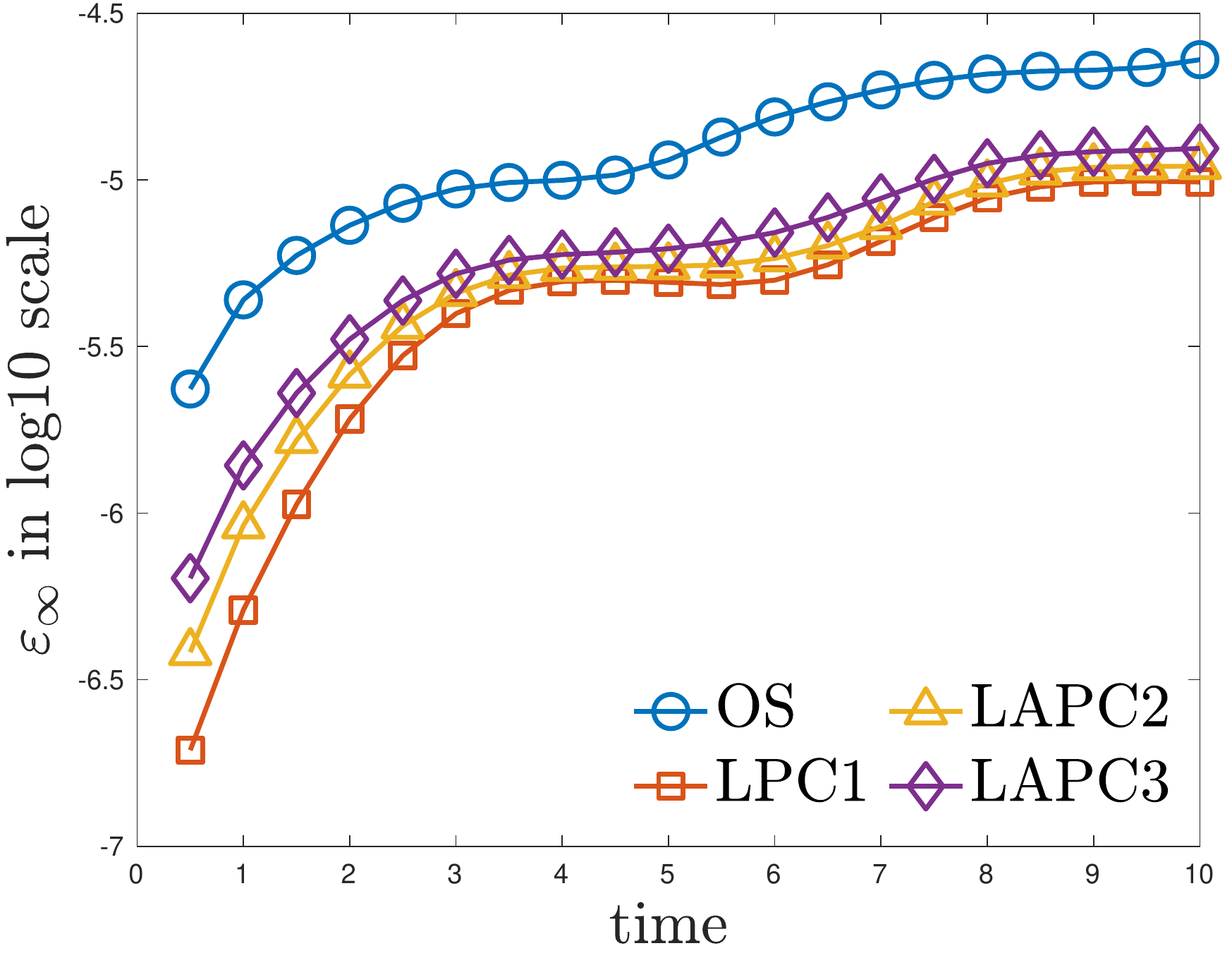}}
{\includegraphics[width=0.32\textwidth,height=0.22\textwidth]{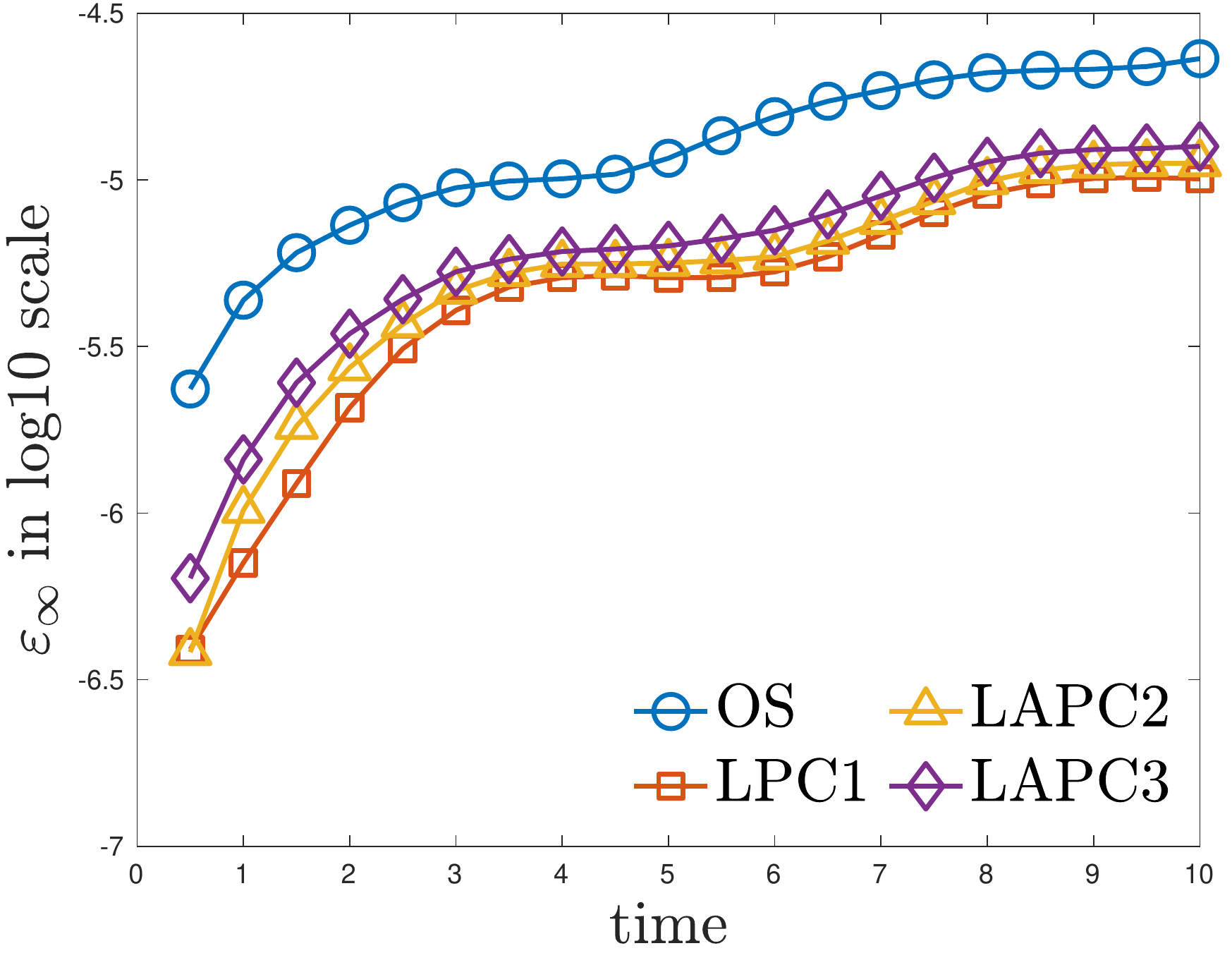}}
{\includegraphics[width=0.32\textwidth,height=0.22\textwidth]{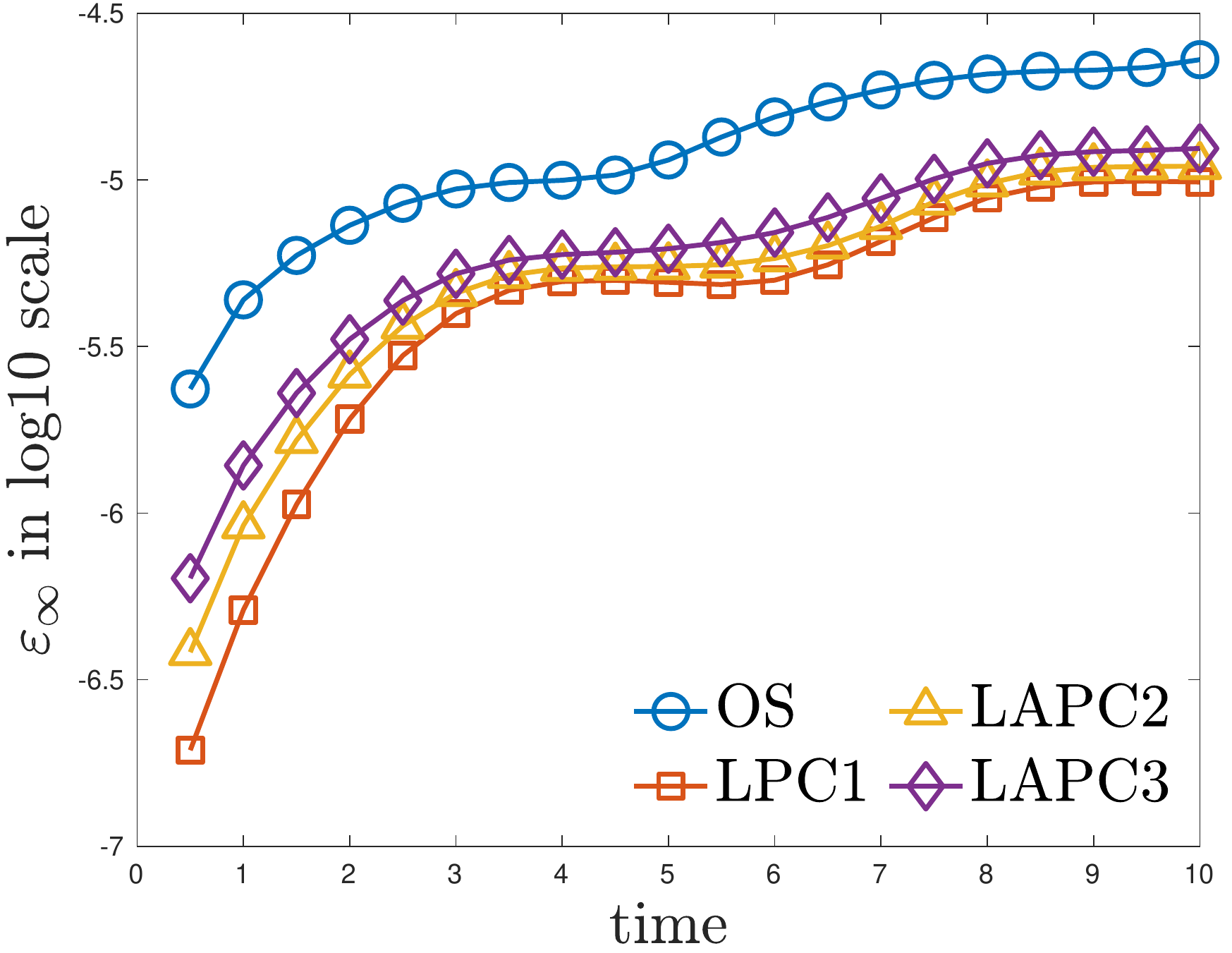}}}
\\
\centering
\subfigure[$\Delta x = 0.05$. \label{harmonic_comp_dx_005}]{
{\includegraphics[width=0.32\textwidth,height=0.22\textwidth]{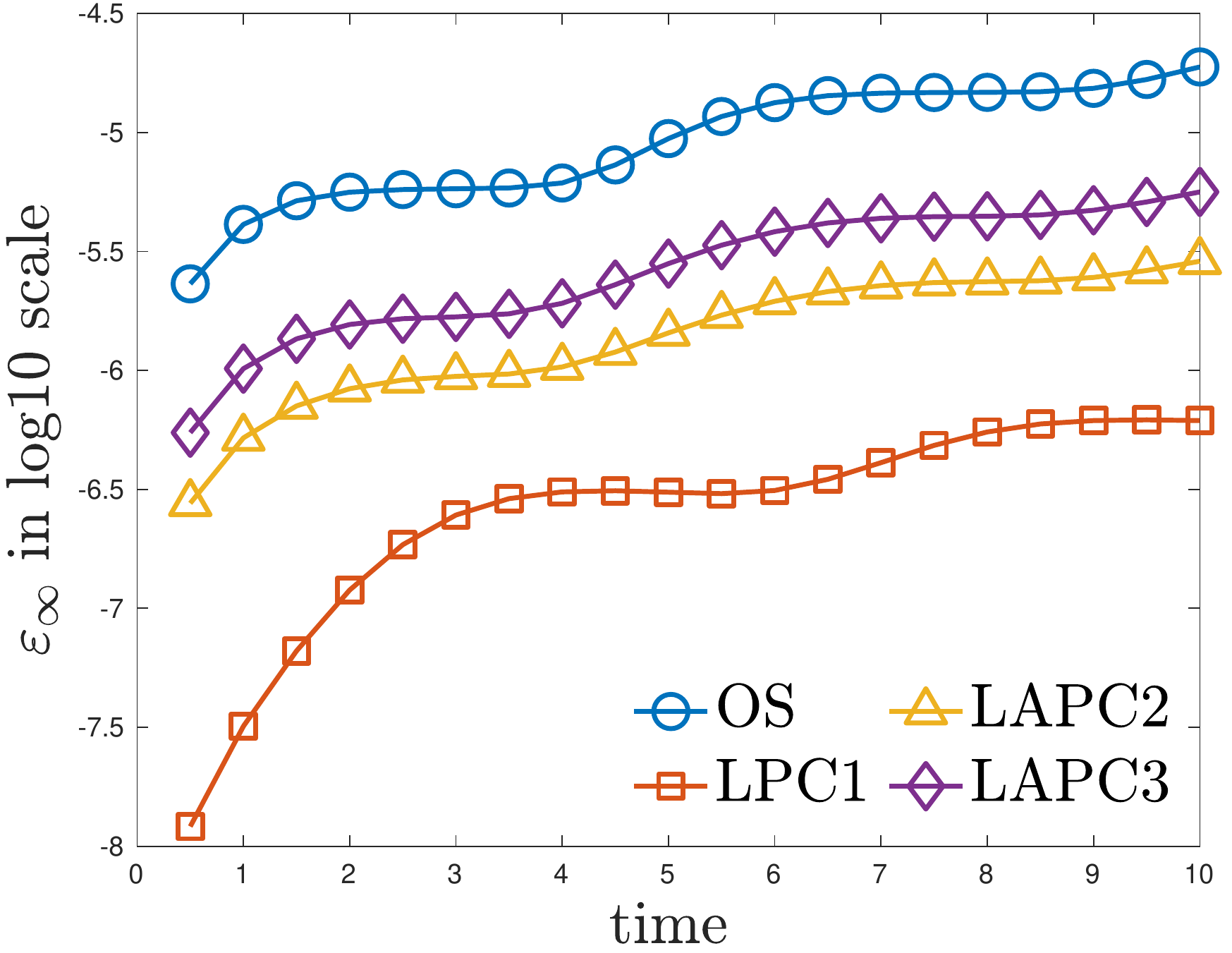}}
{\includegraphics[width=0.32\textwidth,height=0.22\textwidth]{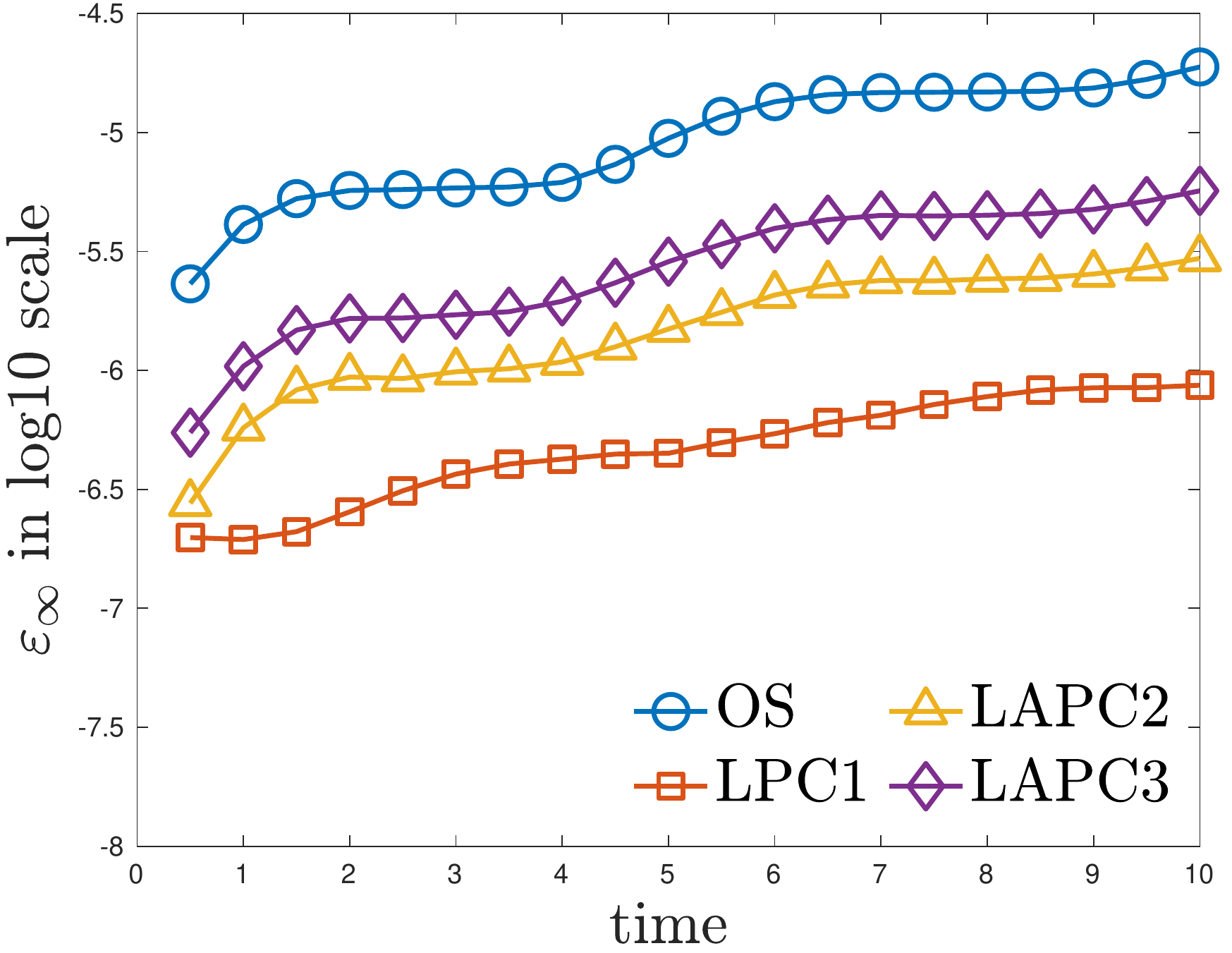}}
{\includegraphics[width=0.32\textwidth,height=0.22\textwidth]{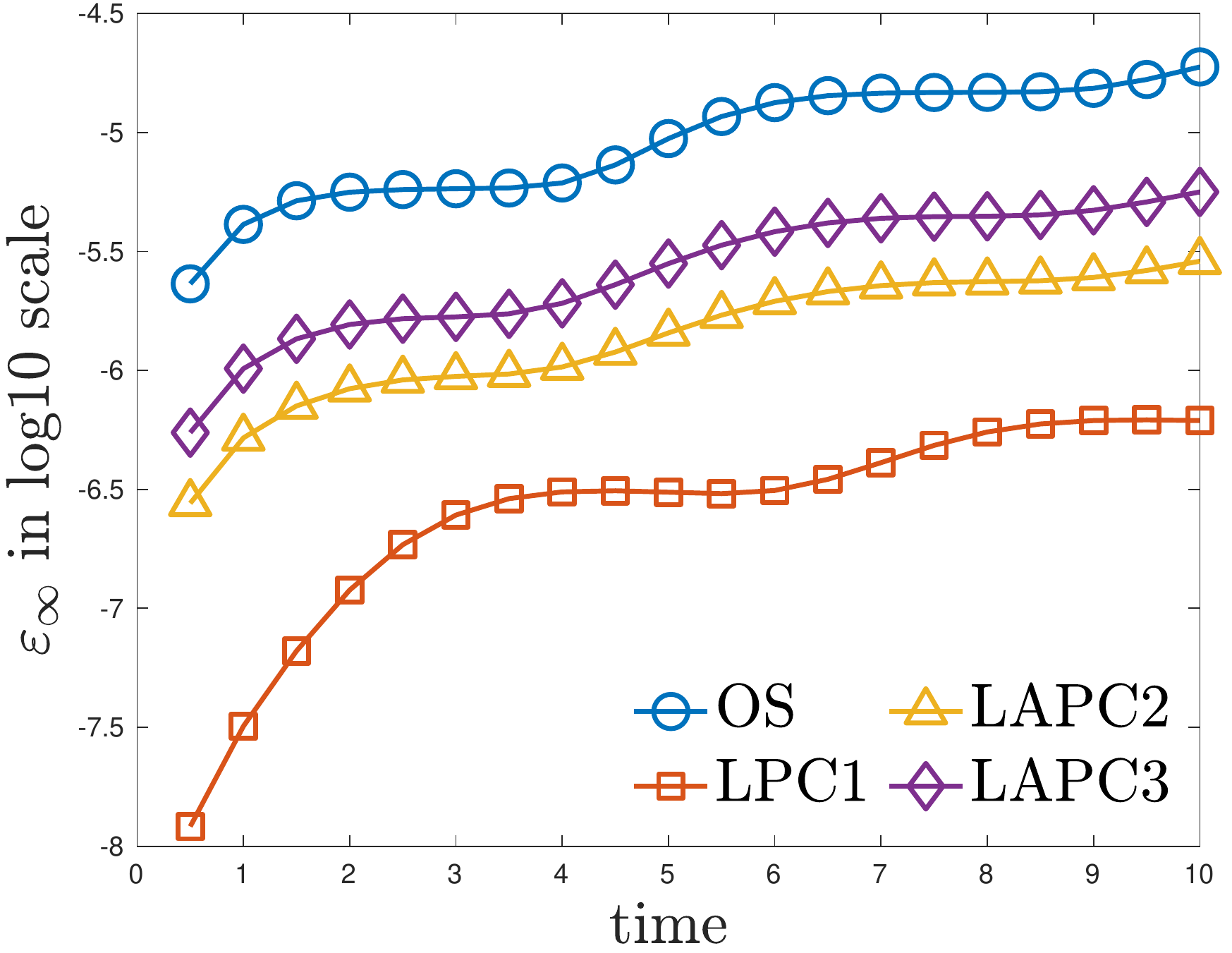}}}
\\
\centering
\subfigure[$\Delta x = 0.025$. \label{harmonic_comp_dx_0025}]{
{\includegraphics[width=0.32\textwidth,height=0.22\textwidth]{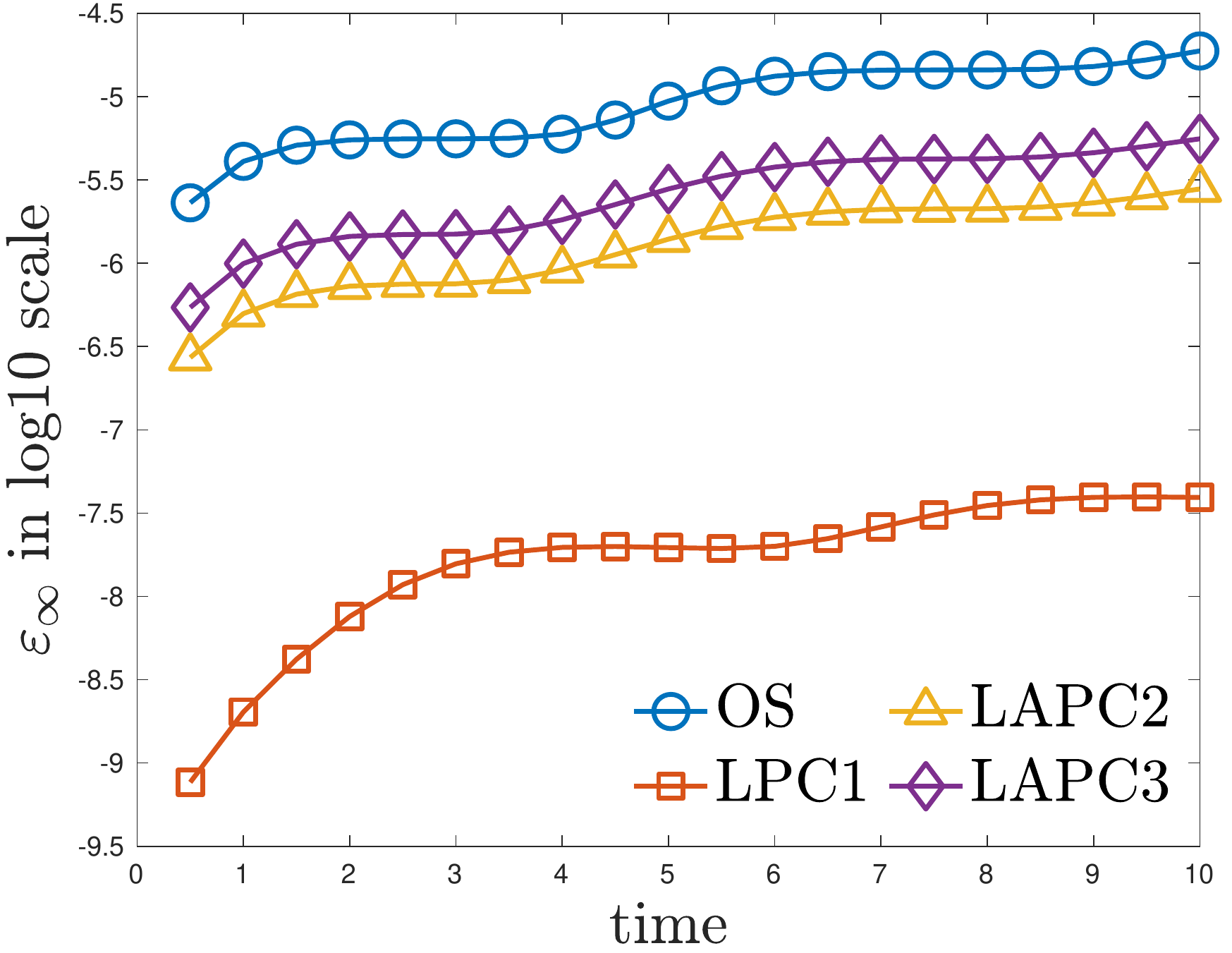}}
{\includegraphics[width=0.32\textwidth,height=0.22\textwidth]{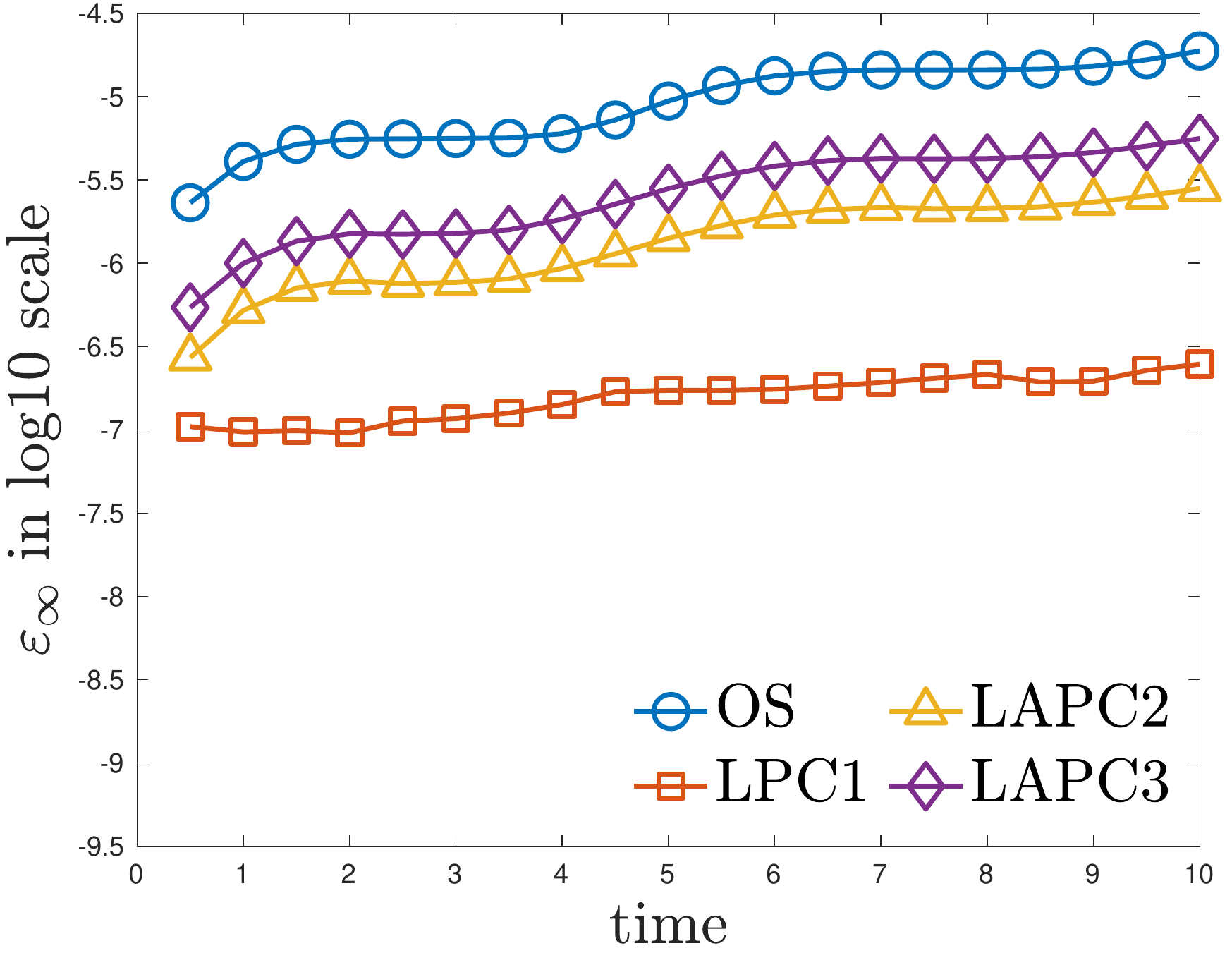}}
{\includegraphics[width=0.32\textwidth,height=0.22\textwidth]{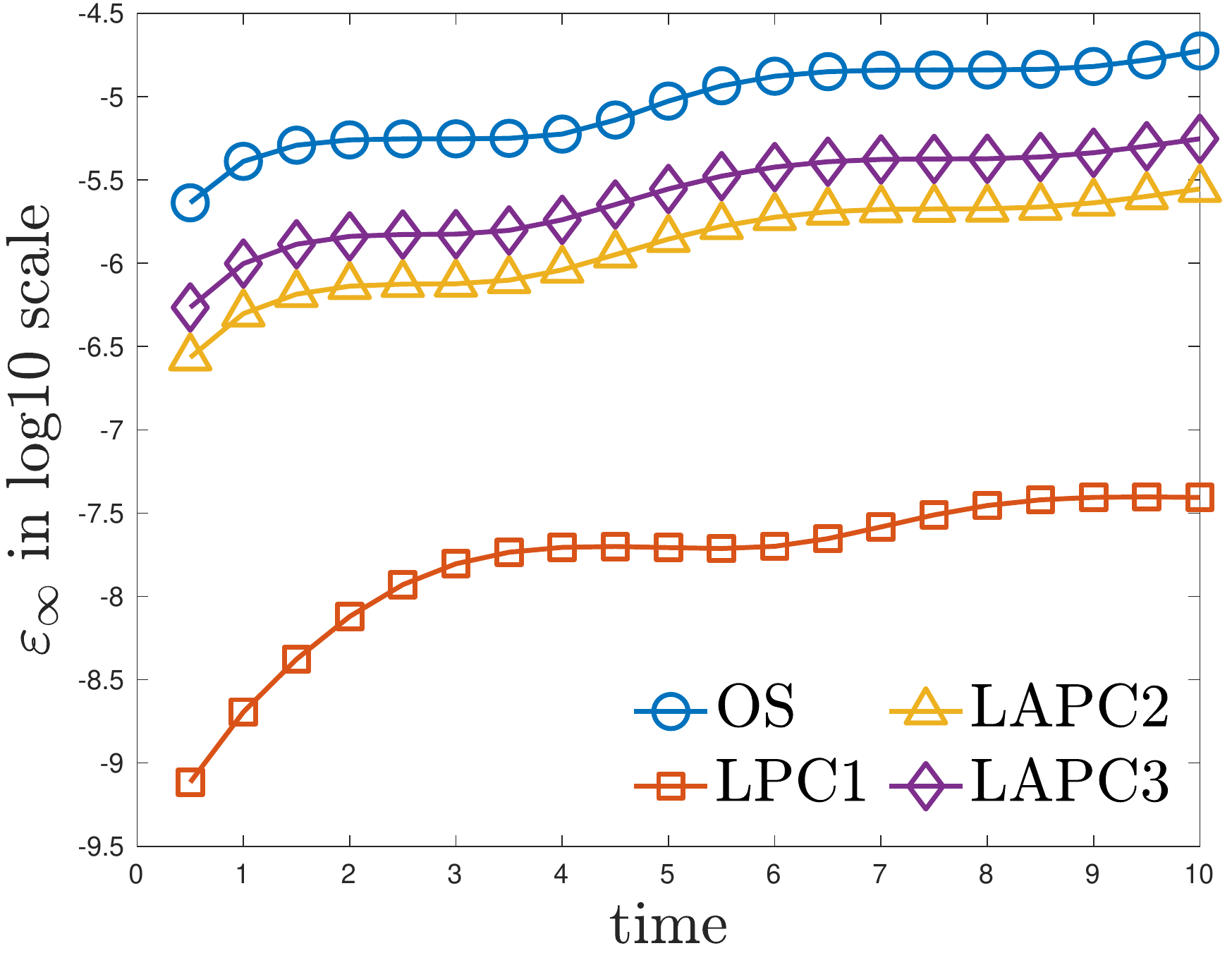}}}
\caption{\small  Quantum harmonic oscillator:  A comparison among different integrators under serial and parallel implementations. (left: serial, middle: $n_{nb=10}$, right: $n_{nb}=20$). LPC1 definitely outperforms other integrators, especially when $\Delta x$ is small.
\label{harmonic_comp_integrator}}
\end{figure}


\begin{figure}[!h]
\centering
\subfigure[Operator splitting scheme (OS).]{{\includegraphics[width=0.32\textwidth,height=0.22\textwidth]{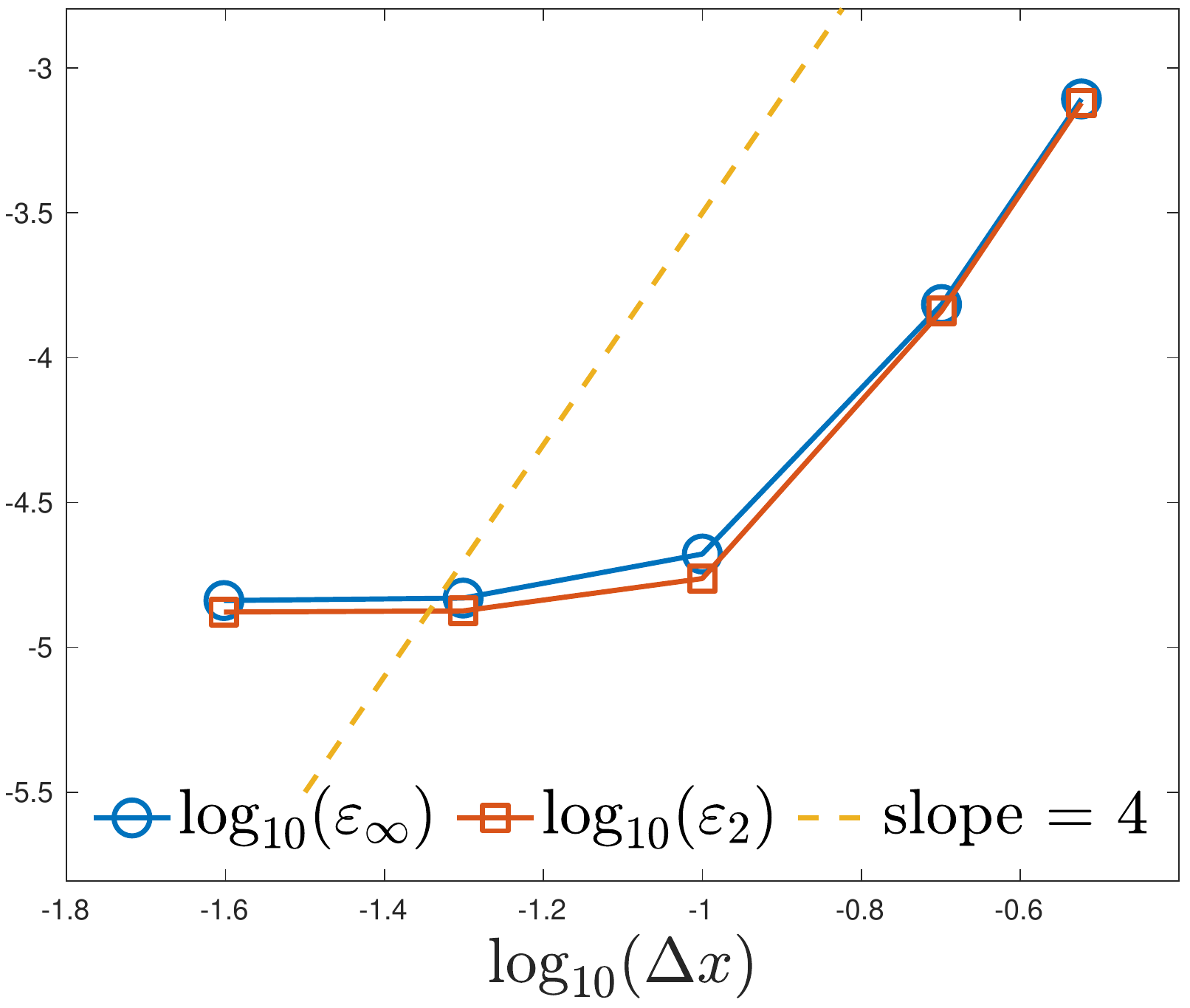}}
{\includegraphics[width=0.32\textwidth,height=0.22\textwidth]{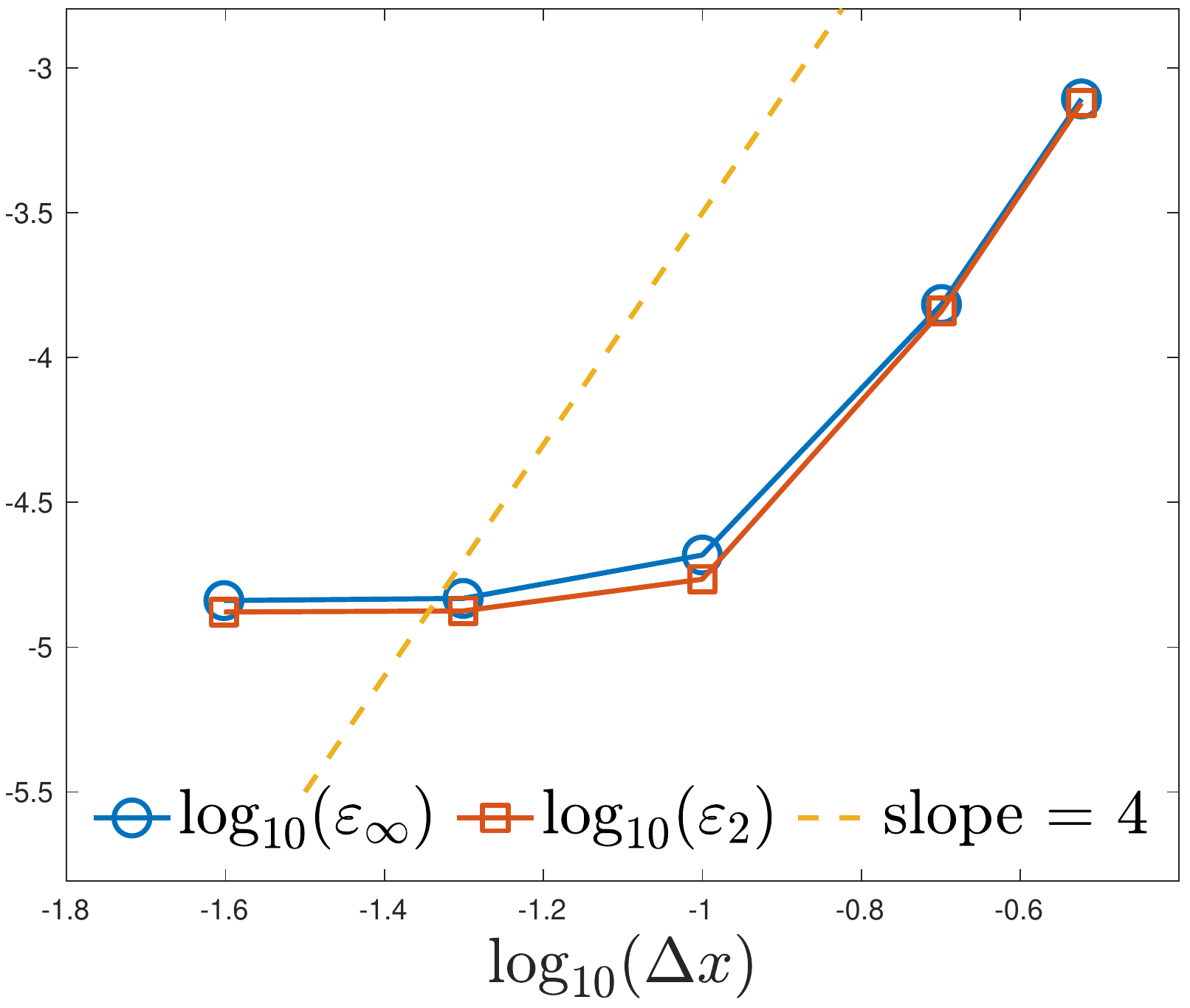}}
{\includegraphics[width=0.32\textwidth,height=0.22\textwidth]{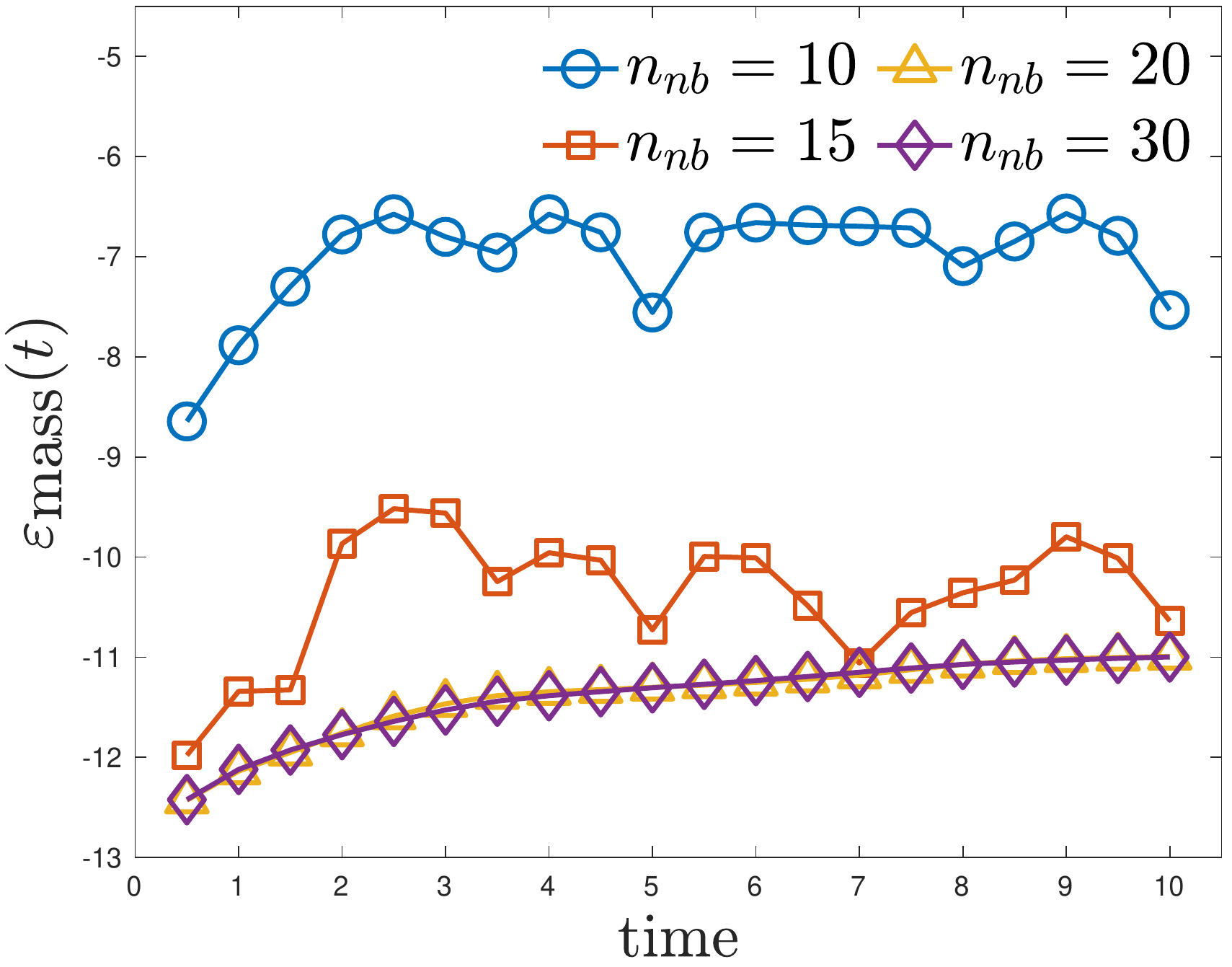}}}
\\
\centering
\subfigure[One-step Lawson predictor-corrector scheme (LPC1).]{{\includegraphics[width=0.32\textwidth,height=0.22\textwidth]{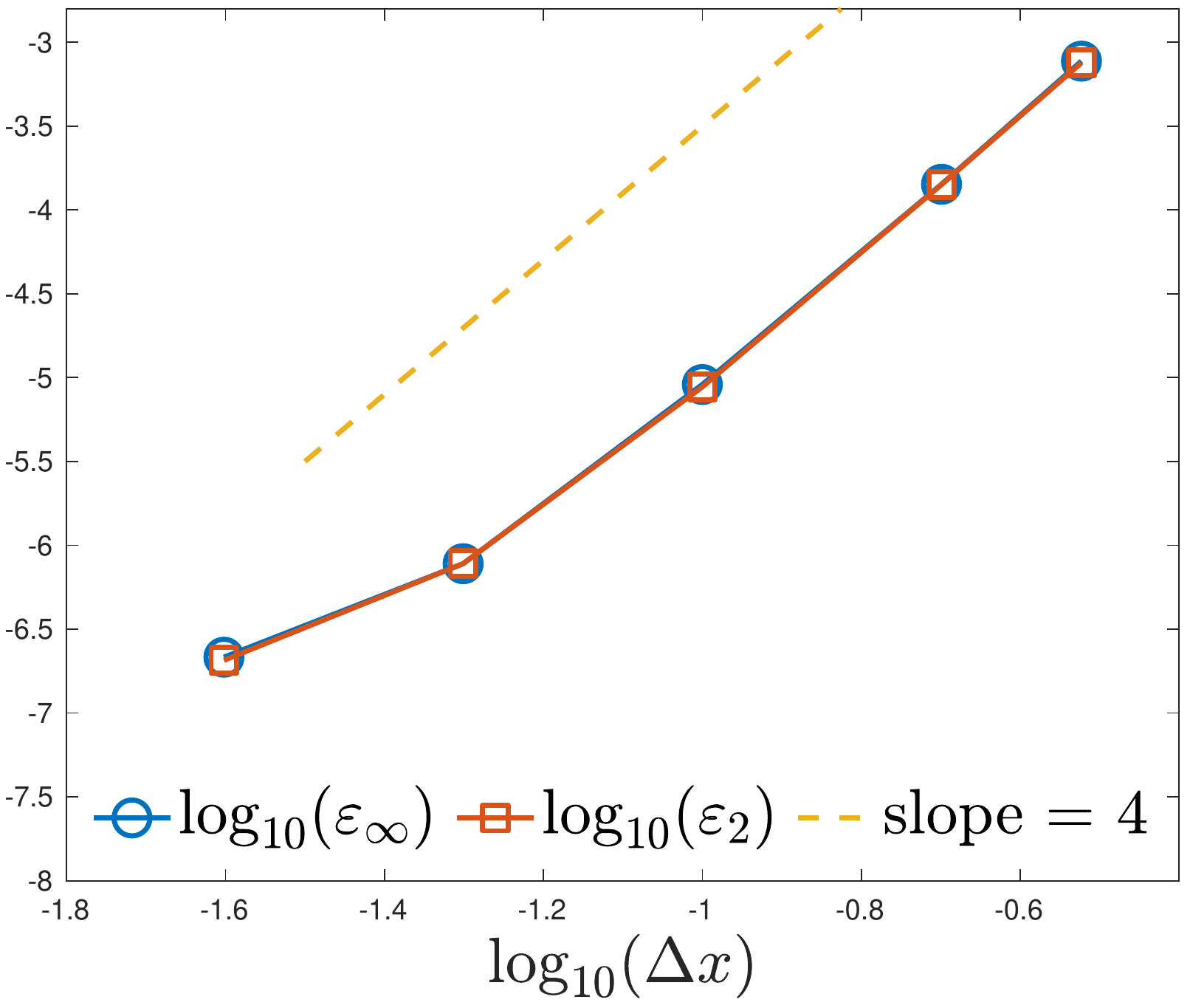}}
{\includegraphics[width=0.32\textwidth,height=0.22\textwidth]{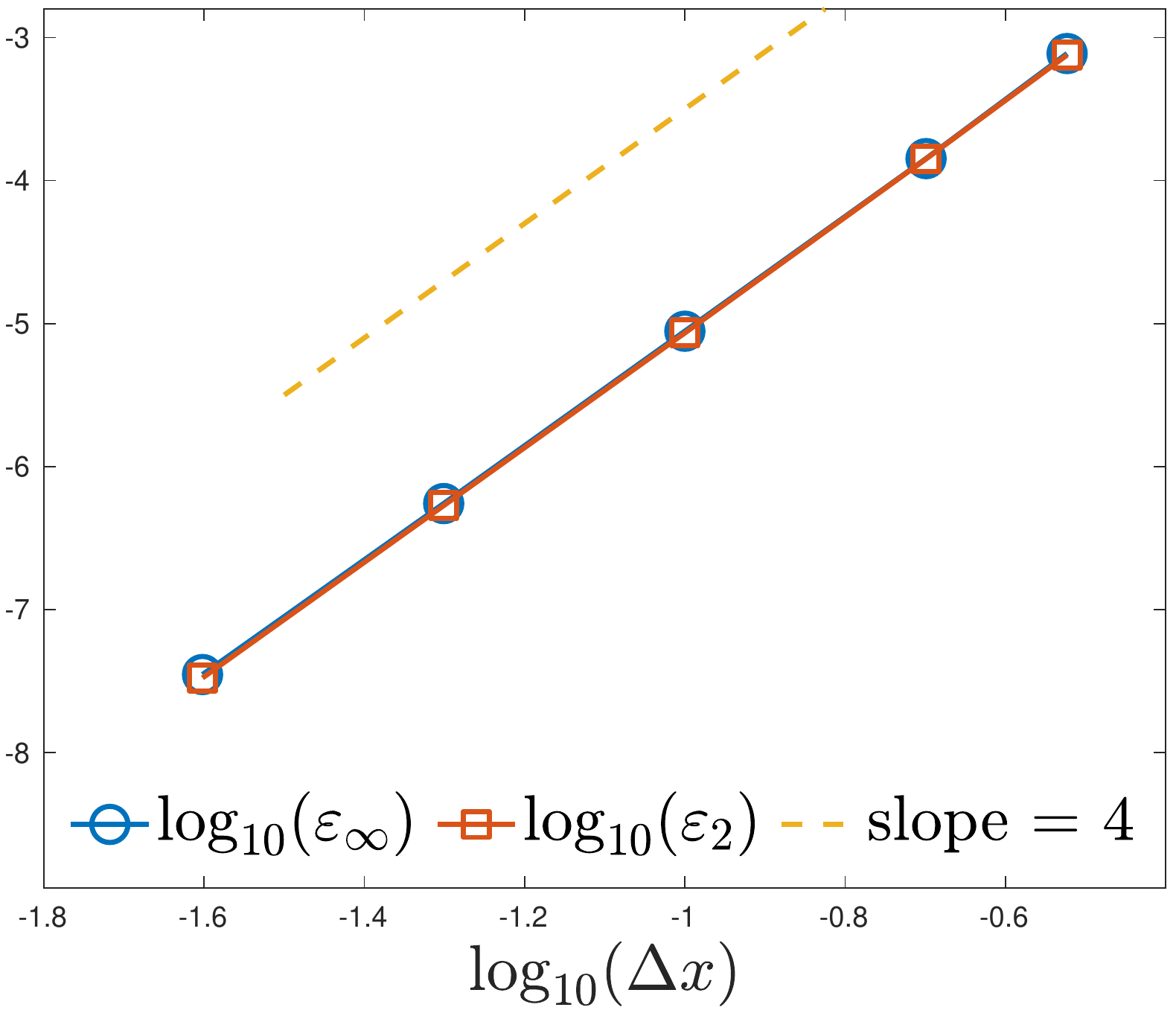}}
{\includegraphics[width=0.32\textwidth,height=0.22\textwidth]{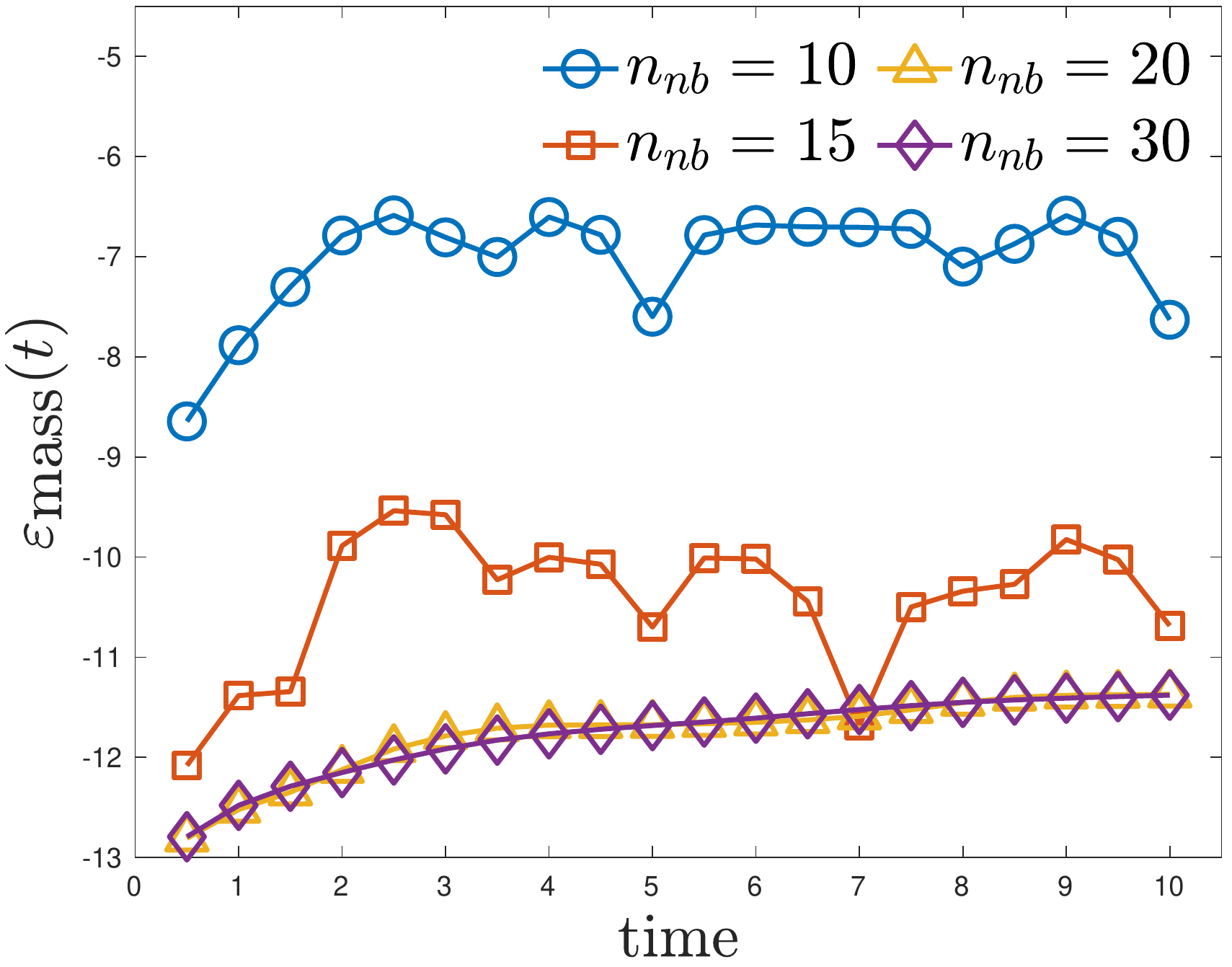}}}
\\
\centering
\subfigure[Two-step Lawson predictor-corrector scheme (LAPC2).]{{\includegraphics[width=0.32\textwidth,height=0.22\textwidth]{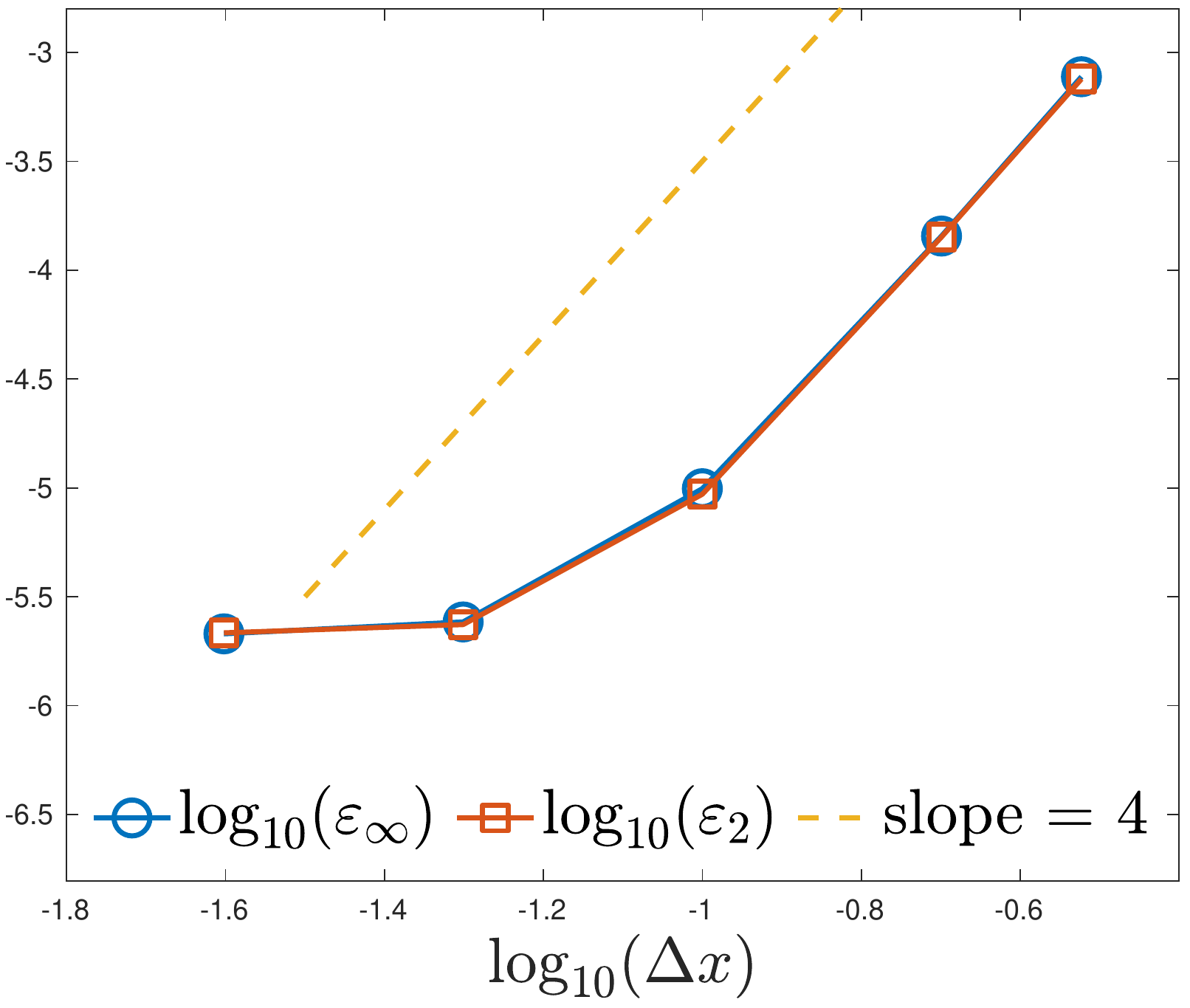}}
{\includegraphics[width=0.32\textwidth,height=0.22\textwidth]{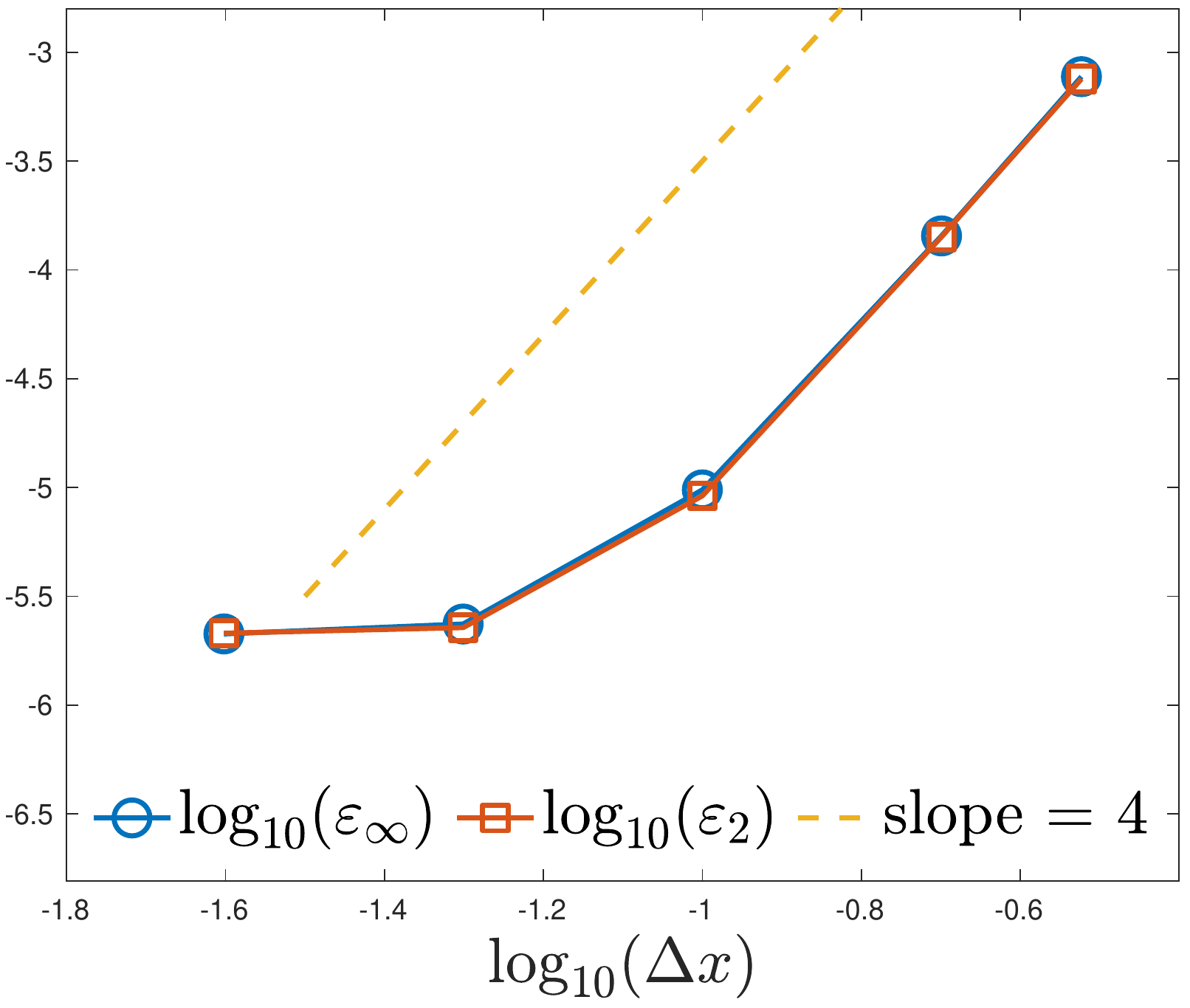}}
{\includegraphics[width=0.32\textwidth,height=0.22\textwidth]{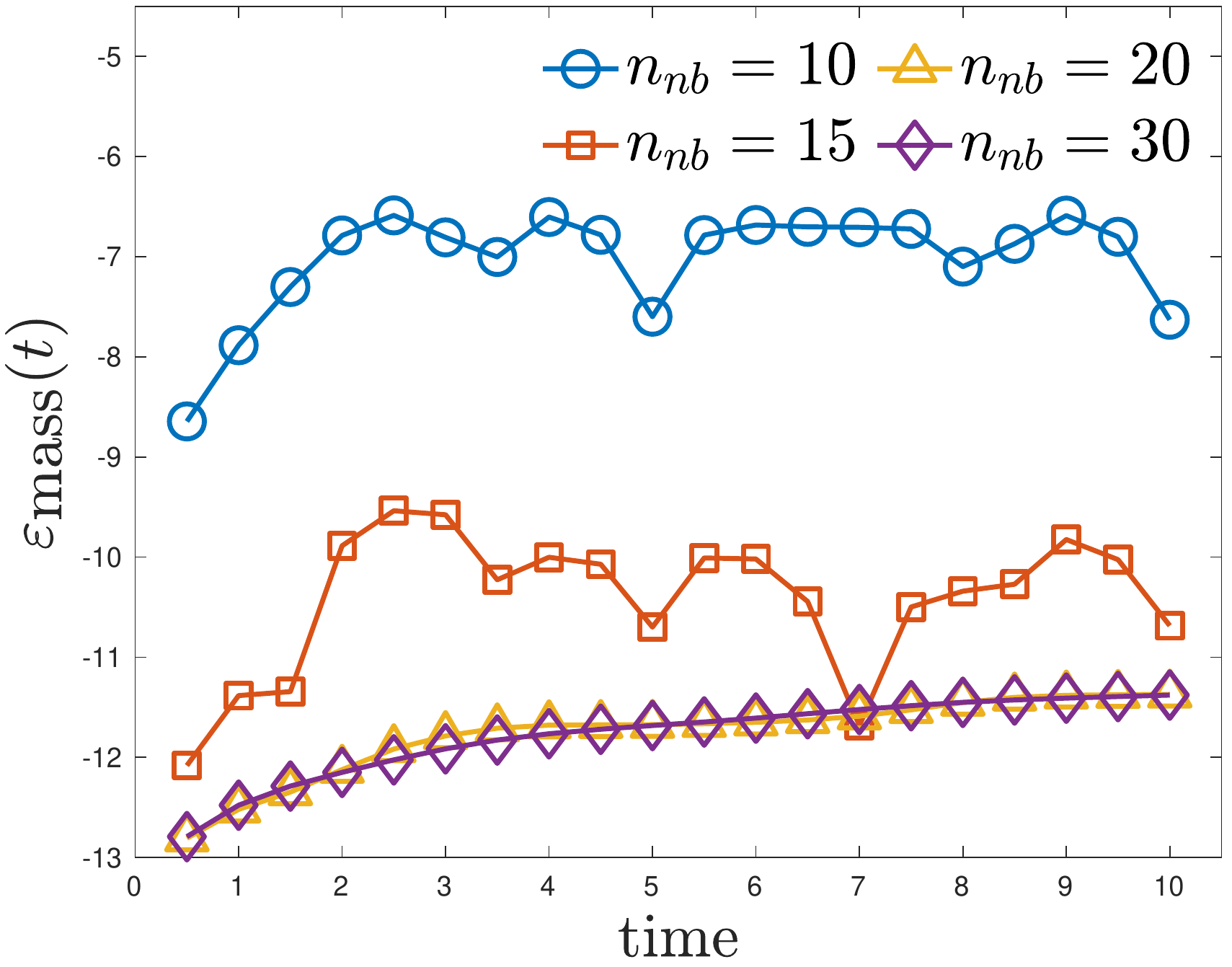}}}
\\
\centering
\subfigure[Three-step Lawson predictor-corrector scheme (LAPC3).]{{\includegraphics[width=0.32\textwidth,height=0.22\textwidth]{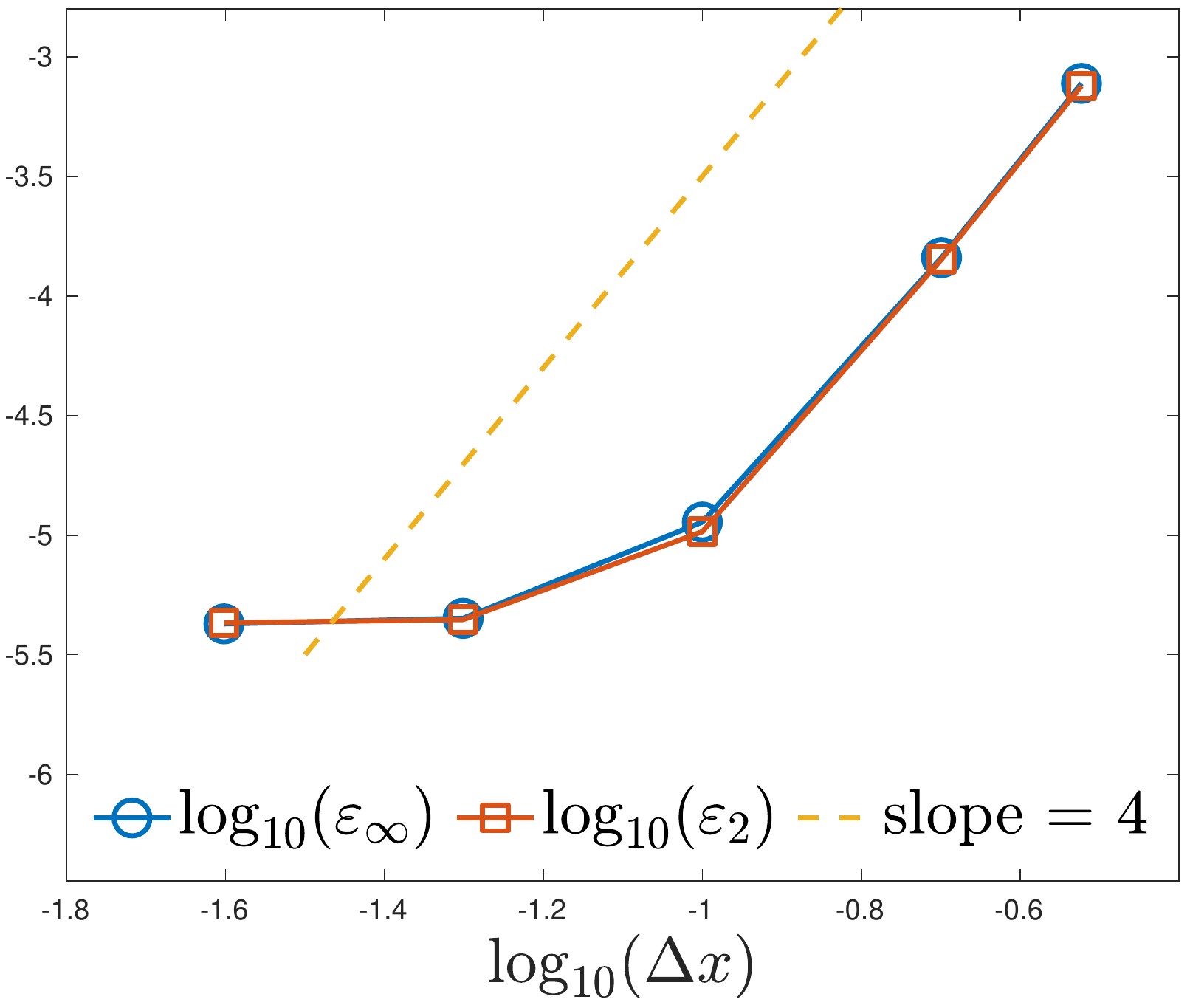}}
{\includegraphics[width=0.32\textwidth,height=0.22\textwidth]{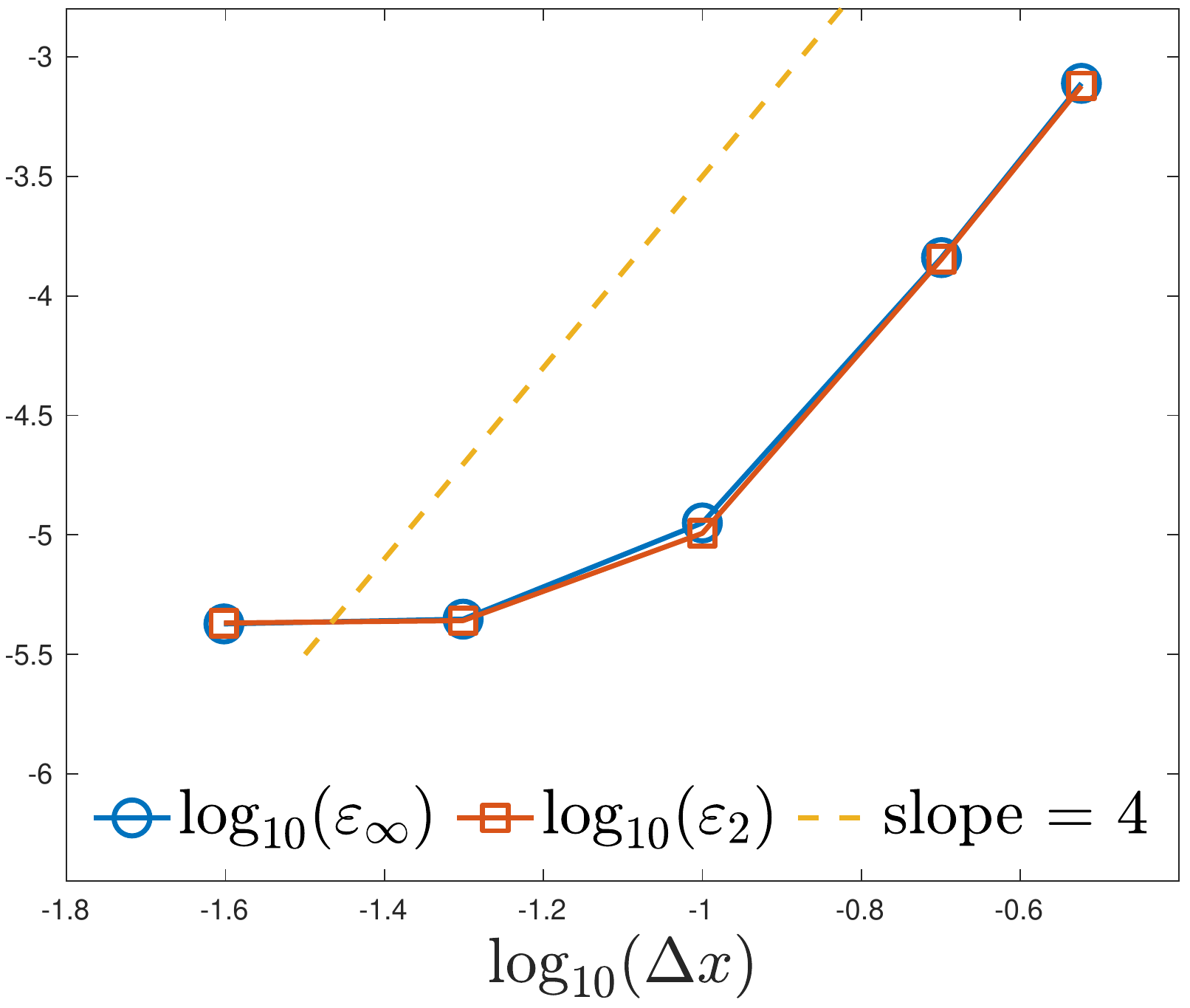}}
{\includegraphics[width=0.32\textwidth,height=0.22\textwidth]{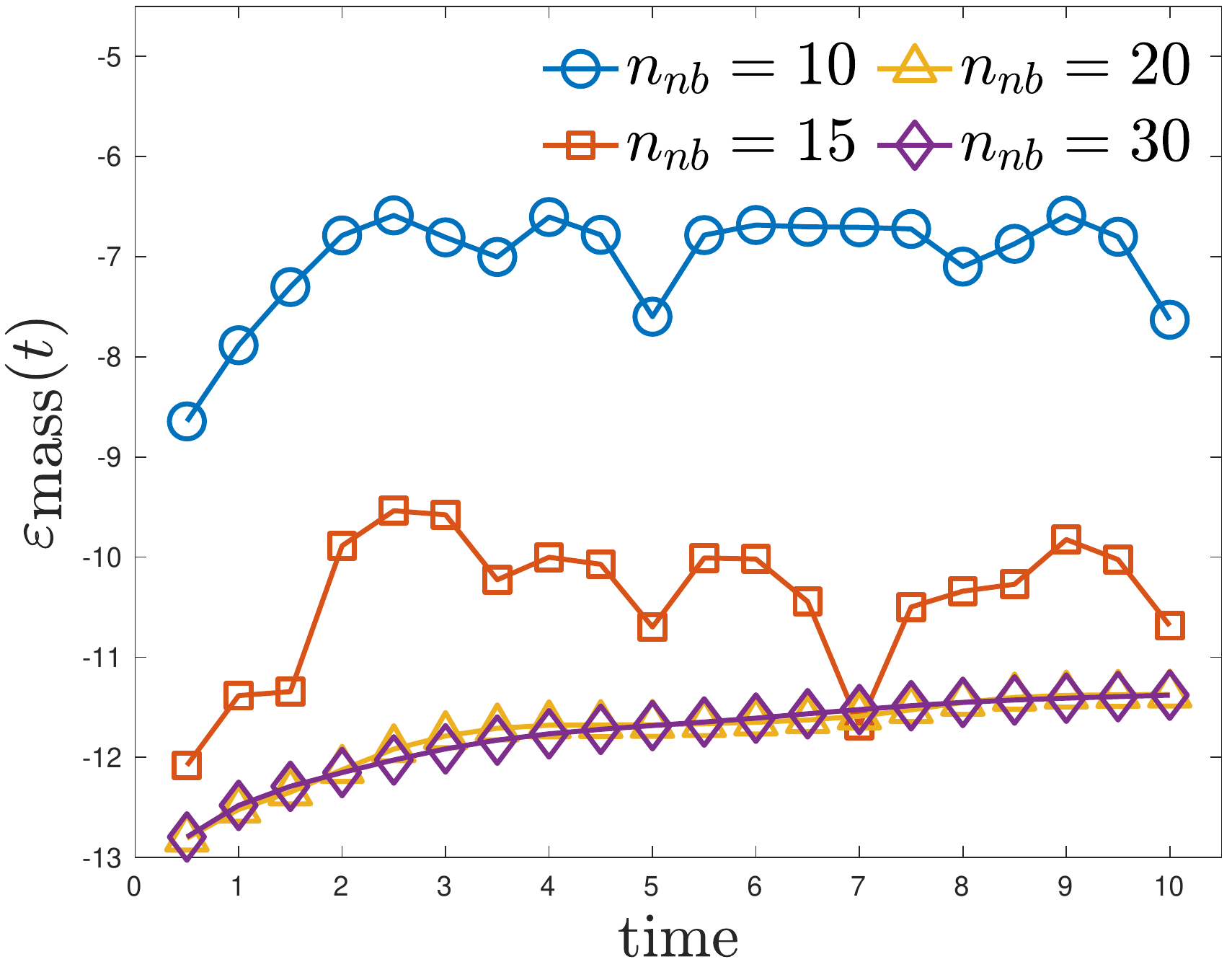}}}
\caption{\small Quantum harmonic oscillator: The convergence (left: $n_{nb=10}$, middle: $n_{nb}=20$) and $\varepsilon_{\textup{mass}}(t)$ (right) of different integrators. LPC1 can achieve fourth-order convergence in $\Delta x$, while other integrators may suffer from the reduction in convergence rate. PMBC indeed has some influences on both accuracy and mass conservation, but fortunately they can be eliminated when $n_{nb} \ge 20$.}
\label{harmonic_convergence}
\end{figure}

\begin{figure}[!h]
\centering
\includegraphics[width=0.48\textwidth,height=0.27\textwidth]{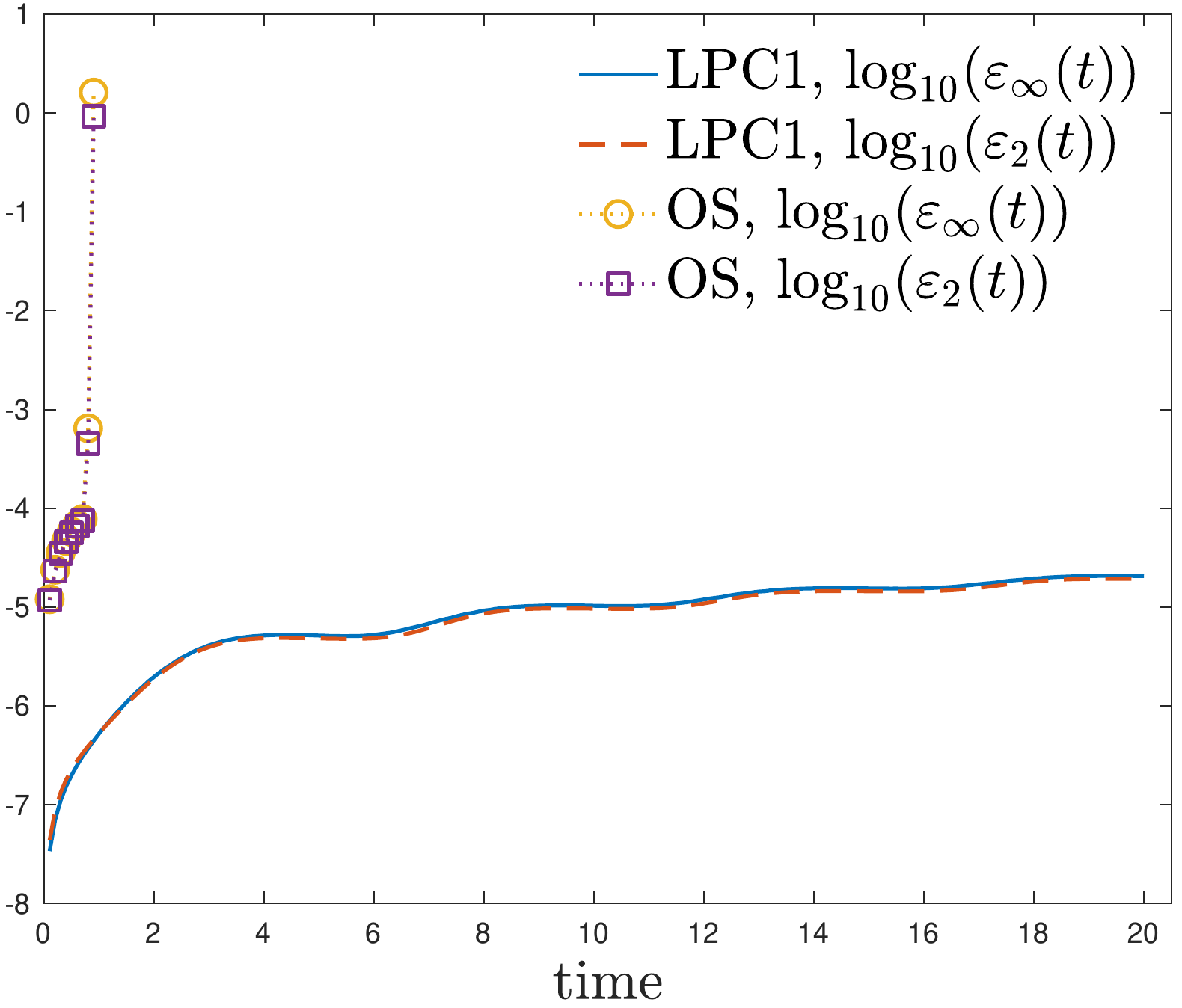}
\caption{\small Quantum harmonic oscillator: The Strang operator splitting suffers from numerical instability under time step $\tau = 0.0005$ and spatial spacing $\Delta x = 0.1$, while LPC1 is stable under such setting even up to $T = 20$. \label{OS_instability}}
\end{figure}

\begin{figure}[!h]
\centering
\subfigure[OS.\label{harmonic_visual_OS}]{
{\includegraphics[width=0.48\textwidth,height=0.27\textwidth]{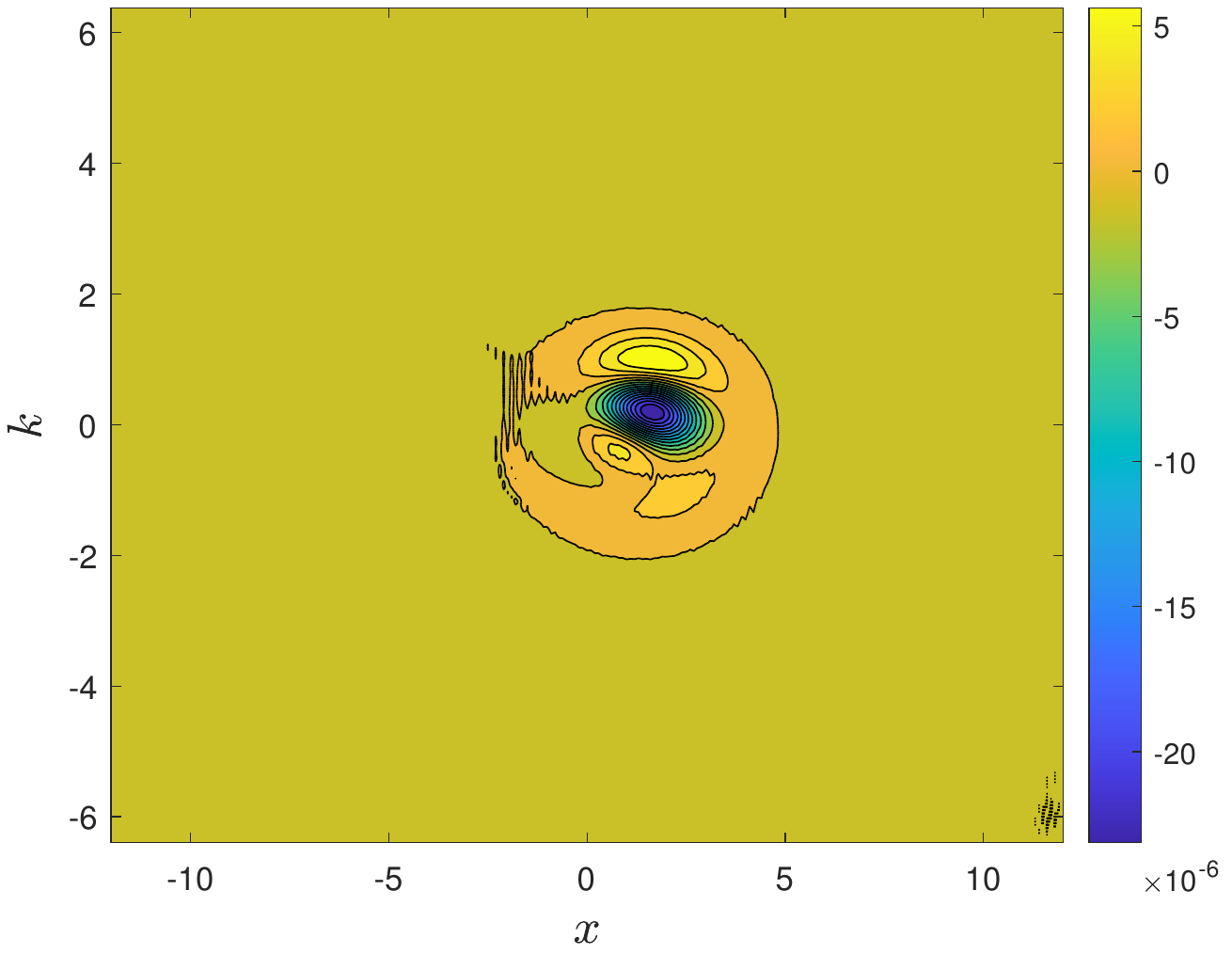}}}
\subfigure[LPC1.]{
{\includegraphics[width=0.48\textwidth,height=0.27\textwidth]{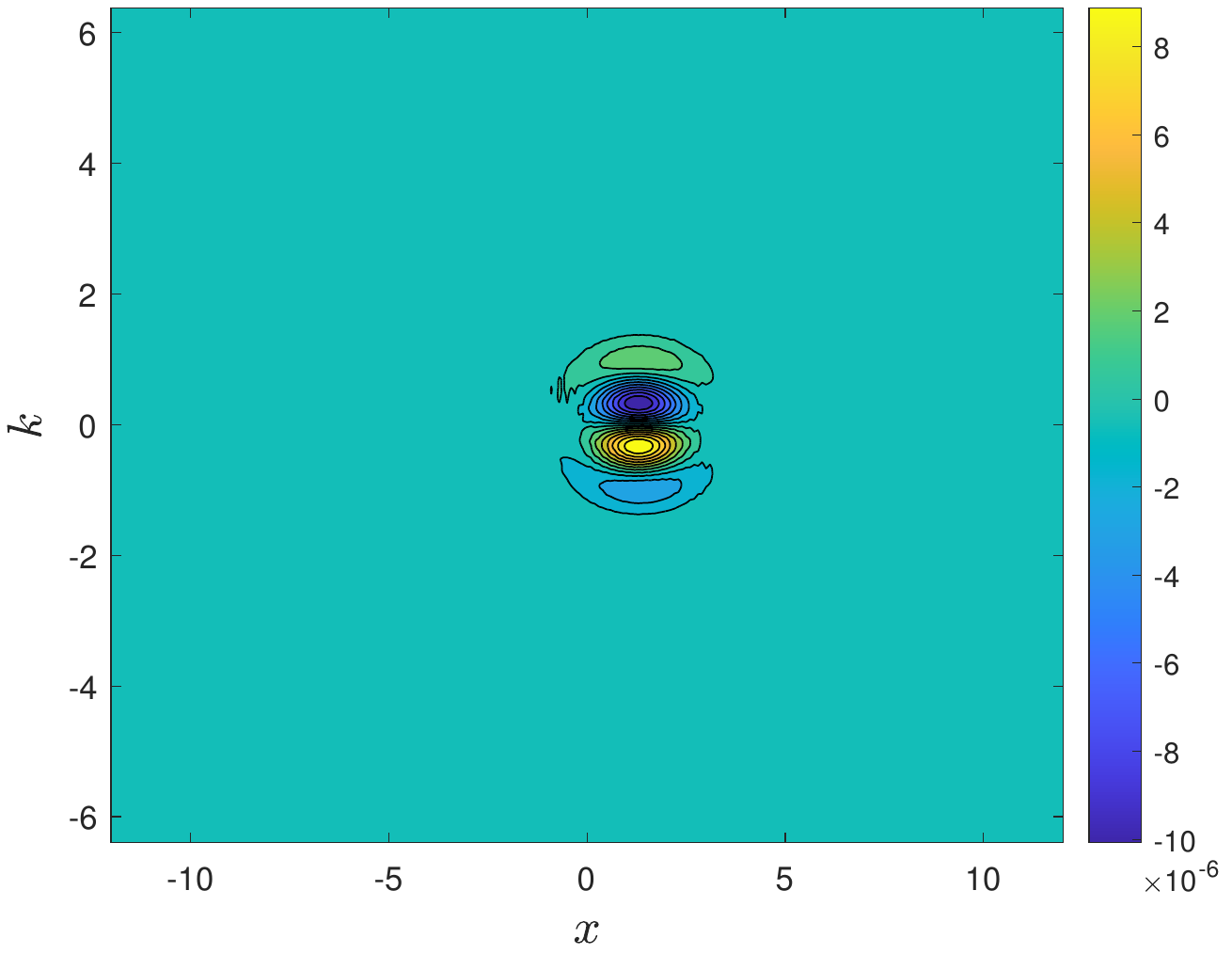}}}
\\
\centering
\subfigure[LAPC2.]{
{\includegraphics[width=0.48\textwidth,height=0.27\textwidth]{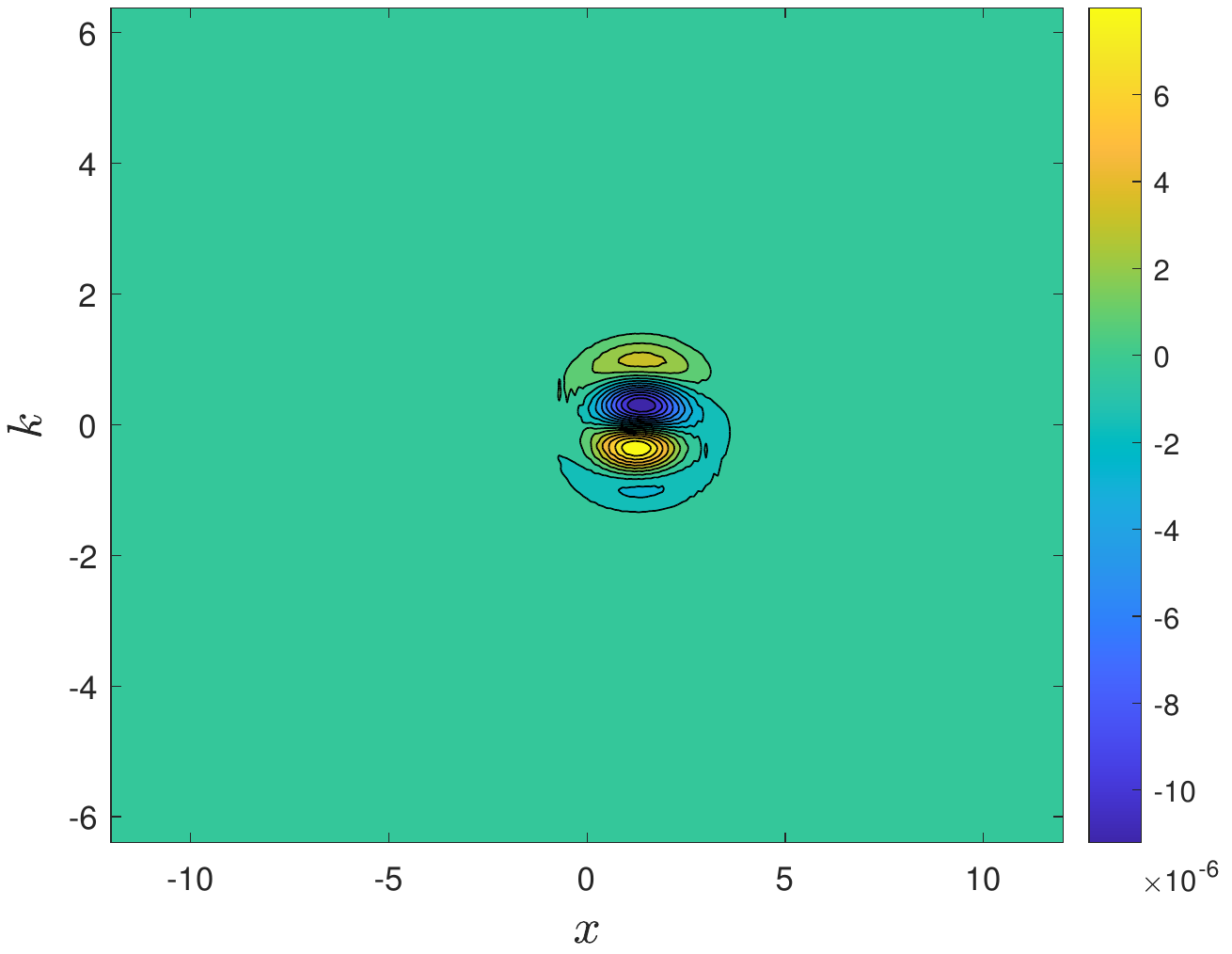}}}
\subfigure[LAPC3. \label{harmonic_visual_LAPC3}]{
{\includegraphics[width=0.48\textwidth,height=0.27\textwidth]{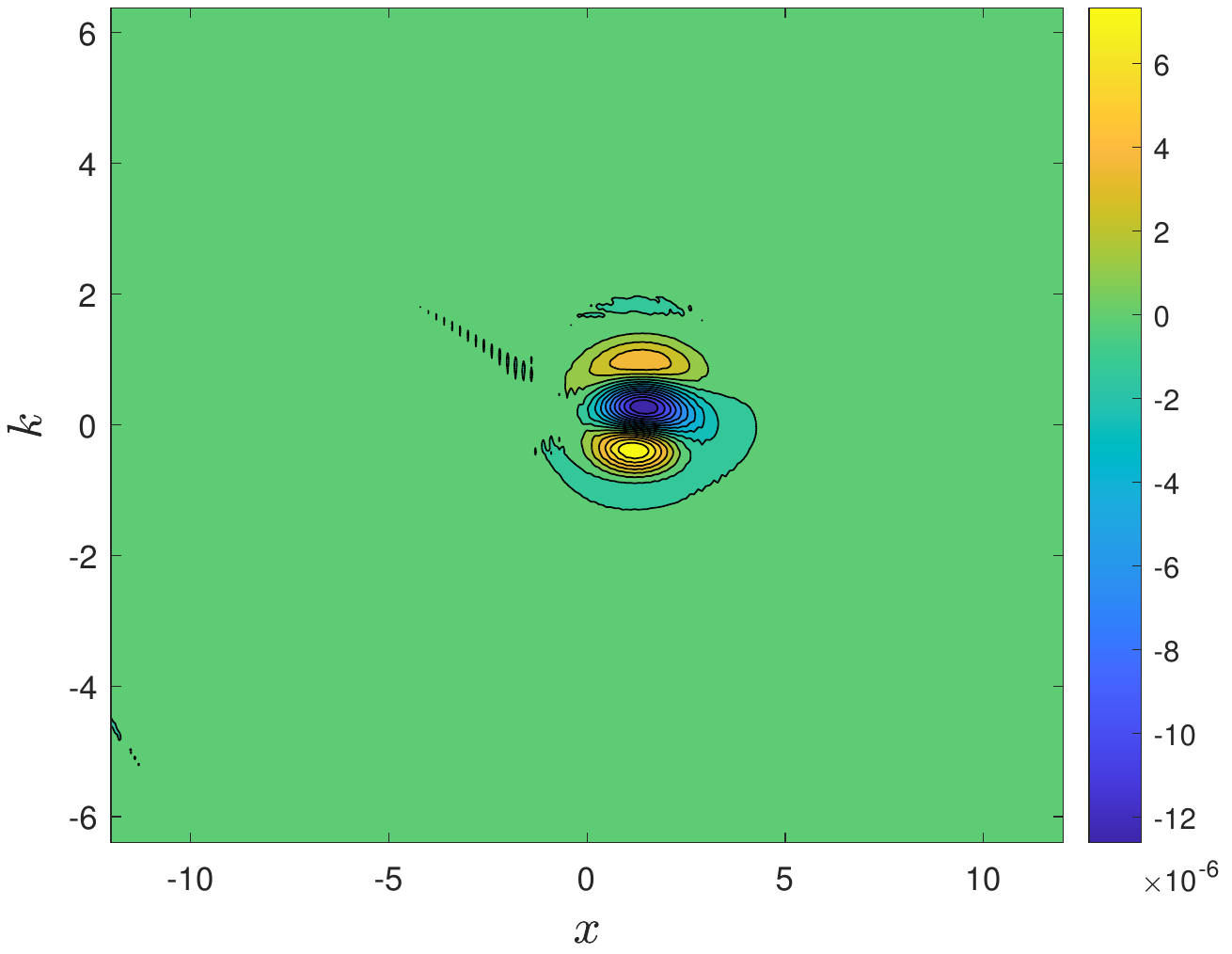}}}
\caption{\small  Quantum harmonic oscillator:  A visualization of numerical errors $f^{\textup{num}}(x, k, t) - f^{\textup{ref}}(x, k, t)$ at $t = 5$  induced by PMBC with $n_{nb} = 10$, $\Delta x =0.1$. It is seen that PMBC may bring in small oscillations at the junction of adjacent patches. The worse is the accumulation of errors near the boundary (see OS and LAPC3), which might lead to numerical instability for long-time evolution.
\label{harmonic_error_visualization}}
\end{figure}

The computational domain is $\mathcal{X} \times \mathcal{K} = [-12, 12] \times [-6.4, 6.4]$, which is evenly decomposed into 4 patches for MPI implementation. The natural boundary condition is adopted at two ends so that there is a slight loss of mass (about $10^{-13}$) up to $T=10$. Since we mainly focus on the convergence with respect to $\Delta x$ and $n_{nb}$,  simulations under $\Delta x = 0.025, 0.05, 0.1, 0.2,0.3$ and $n_{nb}=10, 15,20,30$ are performed, where other parameters are set as: the time step $\tau  = 0.00002$ to avoid numerical stiffness and $\Delta k = 0.025$ to achieve very accurate approximation to $\pdo$. A comparison of all integrators under different $\Delta x$ and $n_{nb}$ is presented in Figure \ref{harmonic_comp_integrator}, and numerical errors $f^{\textup{num}} - f^{\textup{ref}}$ are visualized in Figure \ref{harmonic_error_visualization}. The convergence with respect to $\Delta x$ and the mass conservation under different $n_{nb}$ are given in Figure \ref{harmonic_convergence}. From the results, we can make the following observations.

{\bf Comparison of non-splitting and splitting scheme:} It is clearly seen that LPC1 outperforms the splitting scheme and multi-stage non-splitting schemes in accuracy, especially when $\Delta x$ is small, because it avoids both the accumulation of the splitting errors and additional spline interpolation errors in multi-stage Lawson scheme. While for sufficiently large $\Delta x$, e.g., $\Delta x = 0.2$ or $0.3$, the performances of all integrators are comparable as the interpolation error turns out to be dominated.

{\bf Numerical stability:} The first order derivative in Eq.~\eqref{Wigner_harmonic} brings in strong numerical stiffness and puts a severe restriction on the time step $\tau$ in the parallel CHAracteristic-Spectral-Mixed (CHASM) scheme. Nevertheless, the non-splitting scheme seems to be more stable than the splitting scheme, and one-step scheme is more stable than multi-stage ones. In Figures~\ref{harmonic_comp_dx_03} and \ref{harmonic_comp_dx_02}, we can observe an abrupt reduction in accuracy for OS.  This is induced by the accumulation of errors near the boundary (see the small oscillations in Figures \ref{harmonic_visual_OS} and \ref{harmonic_visual_LAPC3}). In fact, LPC1 turns out to be stable up to $T=20$ even under a larger time step $\tau =0.0005$ and $\Delta x = 0.1$, while OS suffers from numerical instability under such setting (see Figure \ref{OS_instability}).

{\bf Convergence with respect to $\Delta x$:} The convergence rate is plotted in Figure \ref{harmonic_convergence}. Only LPC1 can achieve fourth order convergence in $\Delta x$, according with the theoretical value of the cubic spline interpolation. By contrast, for other schemes, the accumulation of errors induced by temporal integration and mixed interpolations contaminate the numerical accuracy, leading to a reduction in convergence order for small  $\Delta x$.

{\bf Influence of PMBCs:} From Figures \ref{harmonic_comp_dx_005} and \ref{harmonic_comp_dx_0025}, one can see that $n_{nb} = 10$ only bring in additional errors about $10^{-5}$, e.g., the small oscillations are found near the junctions of patches in Figure \ref{harmonic_error_visualization}. But such errors seem to be negligible when $n_{nb} \ge 15$, which also coincides with the observations made in \cite{MalevskyThomas1997}. However, the truncation of stencil indeed has a great influence on the mass conservation as seen in Figure  \ref{harmonic_convergence}, where $\varepsilon_{\textup{mass}}$ is about $10^{-6}$ when $n_{nb}=10$ or $10^{-9}$ when $n_{nb}=15$. Fortunately, its influence on total mass can be completely eliminated when $n_{nb} \ge 20$.

{\bf Efficiency:} For one-step evolution, OS requires spatial interpolations twice and calculation of $\pdo$ once, while LPC1 requires   spatial interpolations once and  calculation of $\pdo$ twice. Thus computational complexity of multi-stage schemes is definitely higher than that of OS and LPC1. 

\subsection{The Wigner function for the Hydrogen 1s state}

The Hydrogen Wigner function is  the stationary solution of the Wigner equation \eqref{eq.Wigner} with the pseudo-differential operator under the attractive Coulomb interaction $V(\bx) = -1/|\bx - \bx_A|$, 
\begin{equation}\label{def.pdo}
\Theta_{V}[f](\bx, \bk, t) = \frac{2}{c_{3, 1}\mi}  \int_{\mathbb{R}^3} \me^{2 \mi (\bx - \bx_A) \cdot \bk^{\prime}} \frac{1}{|\bk^{\prime}|^2}  ( f(\bx, \bk - \bk^{\prime}, t) - f(\bx, \bk+\bk^{\prime}, t) )\D \bk^{\prime}.
\end{equation}
The twisted convolution of the form \eqref{def.pdo} can be approximated by the truncated kernel method \cite{VicoGreengardFerrando2016,GreengardJiangZhang2018}.

For the 1s orbital, $\phi_{\textup{1s}}(\bx) = \frac{1}{2\sqrt{2} \pi^2} \exp( - |\bx|)$,  and the corresponding Wigner function reads 
\begin{equation}\label{1s_Wigner_function}
f_{\textup{1s}}(\bx, \bk) = \frac{1}{(2\pi)^3}\int_{\mathbb{R}^3} \phi_{\textup{1s}}(\bx - \frac{\by}{2}) \phi_{\textup{1s}}^\ast(\bx + \frac{\by}{2}) \me^{-\mi \bk \cdot \by} \D \by.
\end{equation}
 Although it is too complicated to obtain an explicit formula \cite{PraxmeyerMostowskiWodkiewicz2005}, the Hydrogen Wigner function of 1s state can be highly accurately approximated by the discrete Fourier transform of Eq.~\eqref{1s_Wigner_function}: For $\bk_{\bzeta} = \bzeta \Delta k$,
\begin{equation*}
f_{\textup{1s}}(\bx, \bk_{\bzeta}) \approx \sum_{\eta_1 = -\frac{N_{y}}{2}}^{\frac{N_{y}}{2}-1} \sum_{\eta_2 = -\frac{N_{y}}{2}}^{\frac{N_{y}}{2}-1} \sum_{\eta_3 = -\frac{N_{y}}{2}}^{\frac{N_{y}}{2}-1} \phi_{\textup{1s}}(\bx - \frac{\bm{\eta} \Delta y}{2}) \phi_{\textup{1s}}^\ast(\bx + \frac{\bm{\eta} \Delta y}{2}) \me^{-  \mi (\bzeta \cdot \bm{\eta} ) \Delta k \Delta y }  (\Delta y)^3.
\end{equation*}
By taking $\Delta y = \frac{2\pi}{N_k \Delta k}$, it can be realized by FFT with $N_y =128$. 

\begin{figure}[!h]
\centering
\subfigure[ $W_1(x, k)$ for 1s orbital.\label{1s_wigner}]{
\includegraphics[width=0.48\textwidth,height=0.27\textwidth]{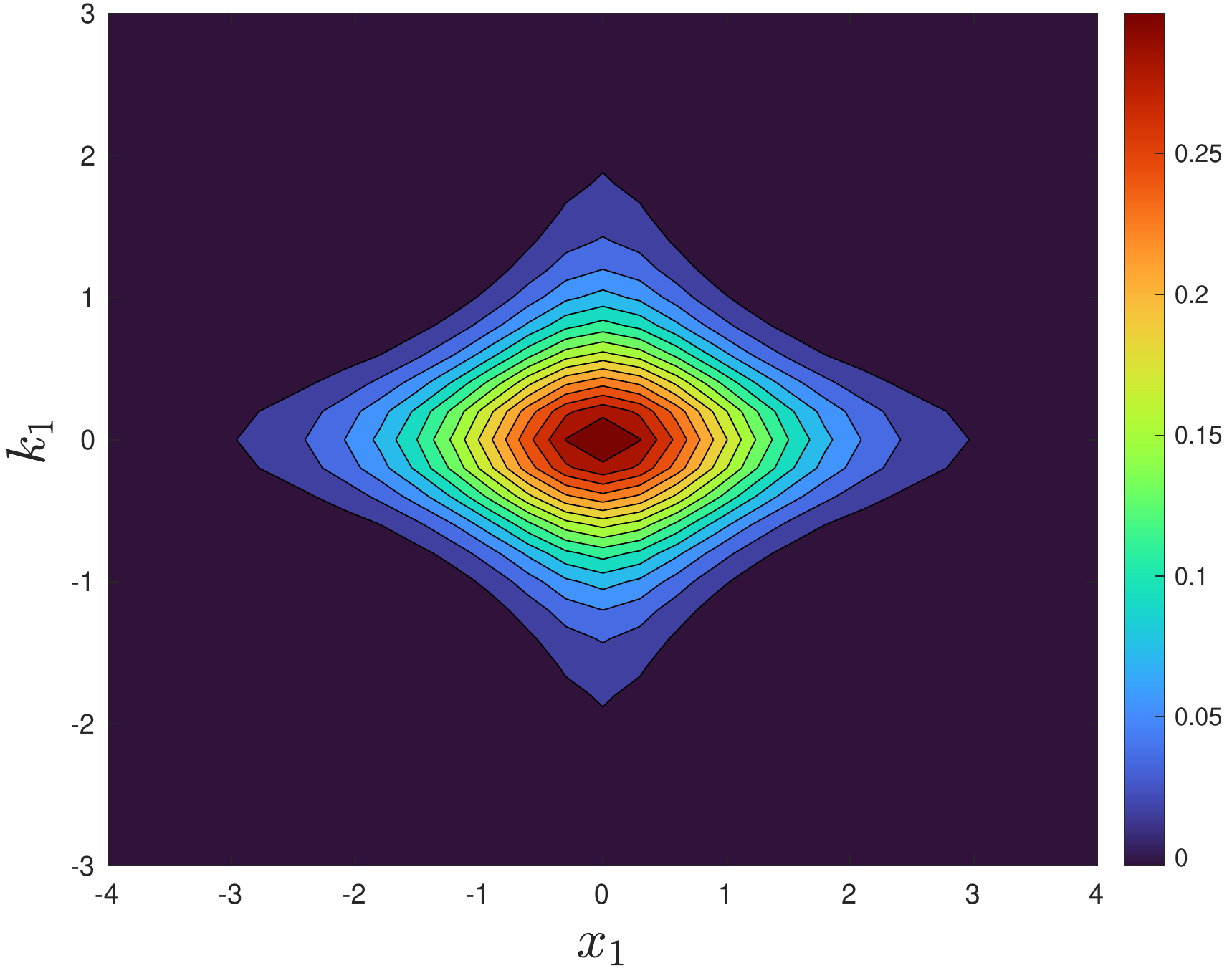}}
\subfigure[The heavy tail in momental space.\label{1s_tail}]{
\includegraphics[width=0.48\textwidth,height=0.27\textwidth]{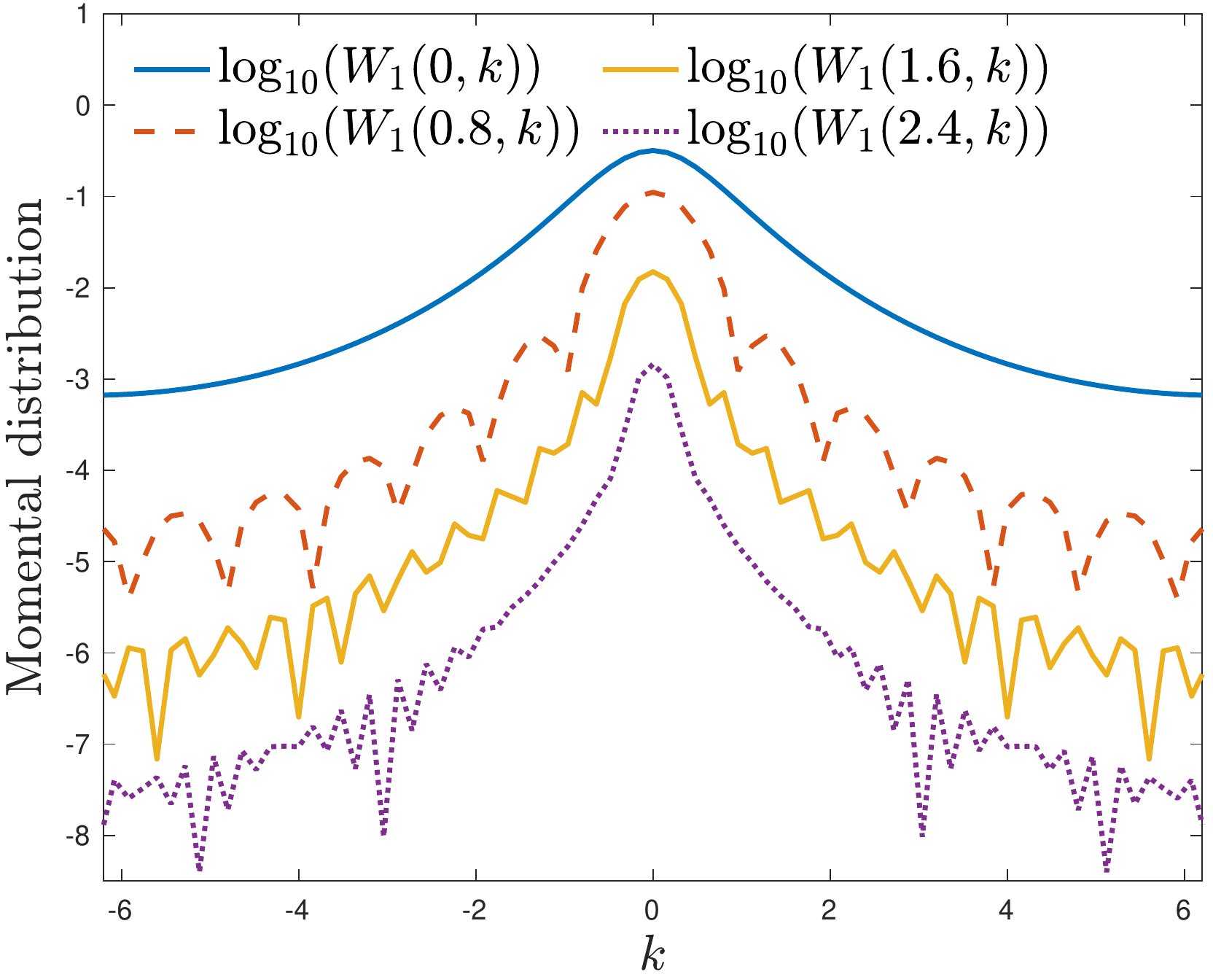}}
\caption{\small The Hydrogen 1s Wigner function: Plot of the reduced Wigner function $W_1(x, k)$. \label{W1_init}}
\end{figure}

The Hydrogen 1s Wigner function can be adopted as the initial and reference solutions for dynamical testing. Besides, for multidimensional case, the reduced Wigner function $W_1(x, k, t) $, defined by the projection of  $f$ onto $(x_1$-$k_1$) plane, is used for visualization. 
\begin{equation}\label{def.reduced_Wigner_function}
W_1(x, k, t) = \iint_{\mathbb{R}^2 \times \mathbb{R}^2} f(\bx, \bk, t) \D x_2 \D x_3 \D k_2 \D k_3.
\end{equation}
For 1s state, the reduced Wigner function is plotted in Figure \ref{1s_wigner}, which exhibits a heavy tail in $\bk$-space as shown in Figure \ref{1s_tail}.

The computational domain is $\mathcal{X} \times \mathcal{K} = [-9, 9]^3 \times [-6.4, 6.4]^3$ with a fixed spatial step size $\Delta x = 0.3$ ($N_{x_1} = N_{x_2} = N_{x_3} = 61$), which is evenly divided into $4\times 4 \times 4$ patches and distributed into $64$ processors, and each processor provides 4 threads for shared-memory parallelization using the OpenMP library. The natural boundary conditions are adopted at two ends. As the accuracy of spline interpolation has been already tested in the above 2-D example, we will investigate the convergence of nonlocal approximation under five groups: $N_k =8, 16, 32, 64, 80$ ($\Delta k  = 1.6, 0.8, 0.4, 0.2, 0.16$). Other parameters are set as: the stencil length in PMBC is $n_{nb} = 15$ and the time stepsize is $\tau = 0.025$.

\begin{figure}[!h]
\centering
\subfigure[LPC, $\varepsilon_\infty(t)$.]{
{\includegraphics[width=0.48\textwidth,height=0.27\textwidth]{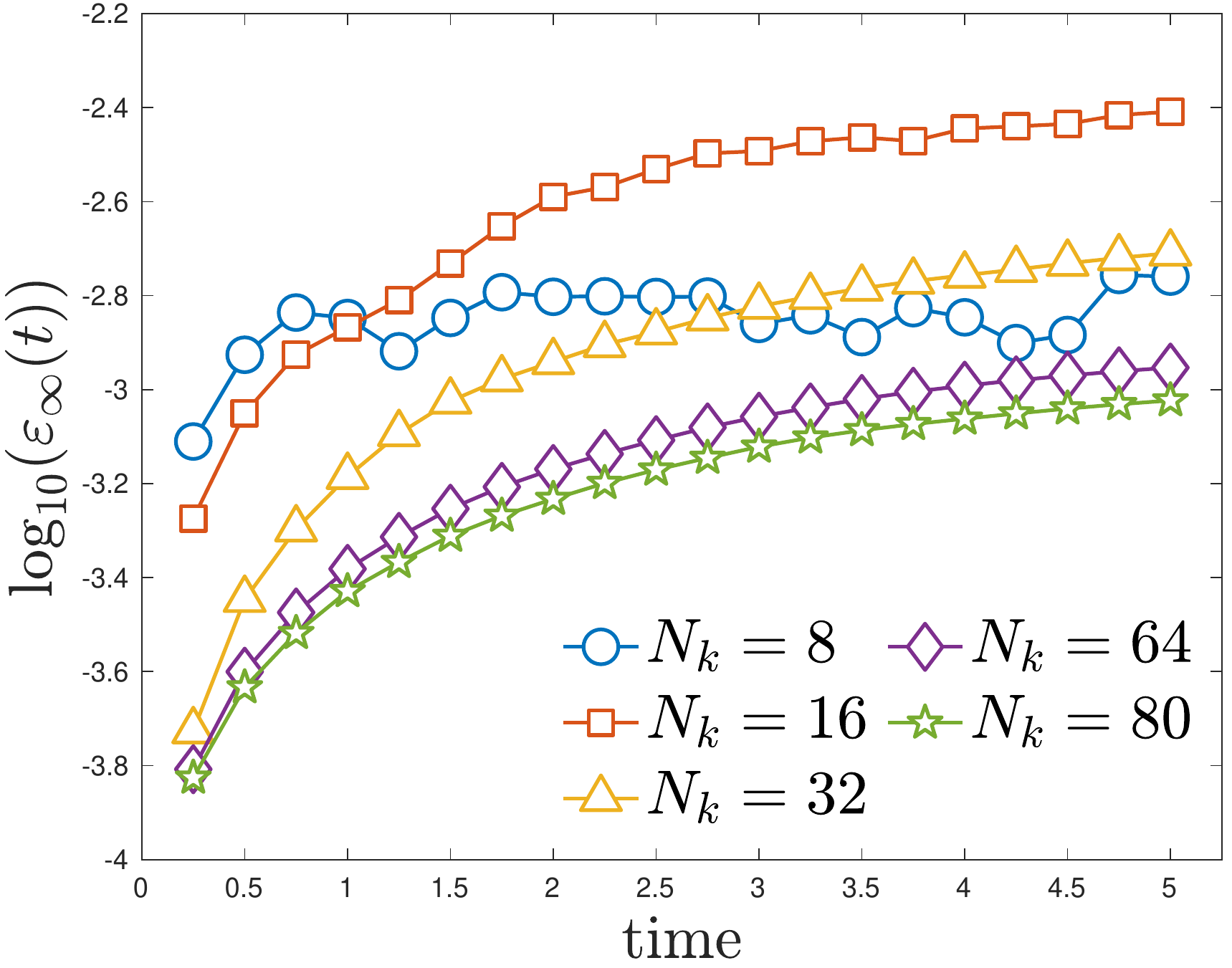}}}
\subfigure[OS, $\varepsilon_\infty(t)$. ]{
{\includegraphics[width=0.48\textwidth,height=0.27\textwidth]{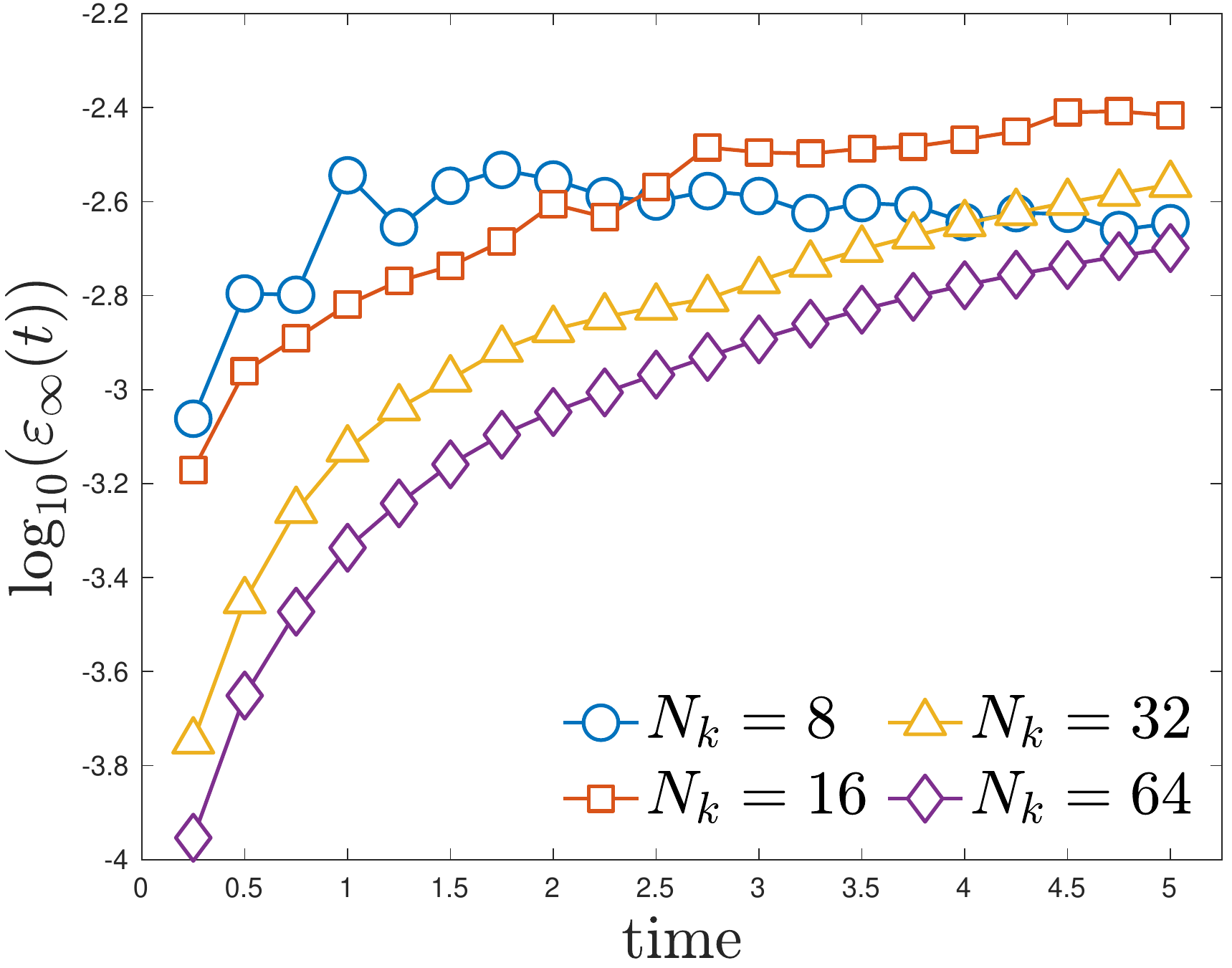}}}
\\
\centering
\subfigure[LPC, $\varepsilon_2(t)$.]{
{\includegraphics[width=0.48\textwidth,height=0.27\textwidth]{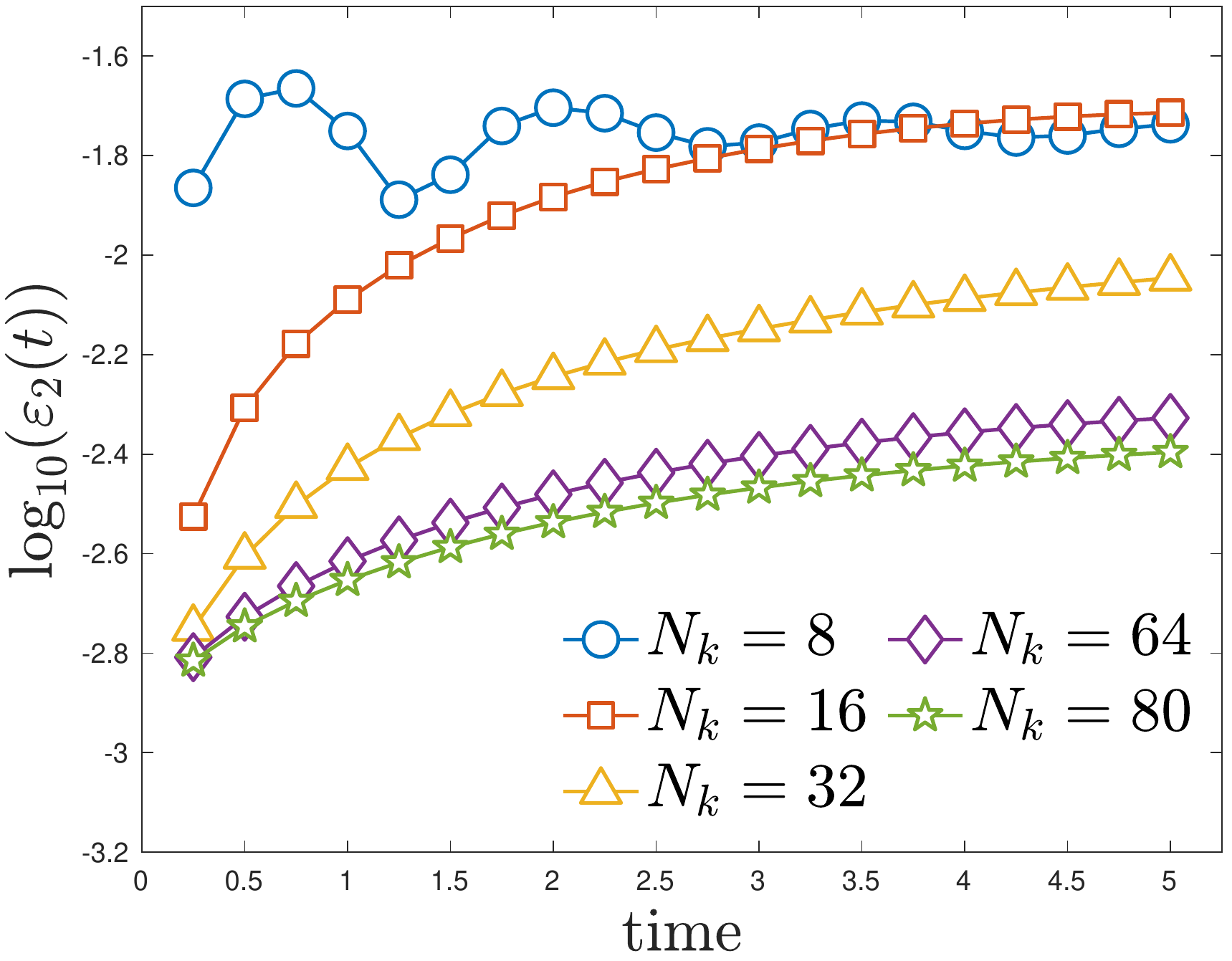}}}
\subfigure[OS, $\varepsilon_2(t)$ ]{
{\includegraphics[width=0.48\textwidth,height=0.27\textwidth]{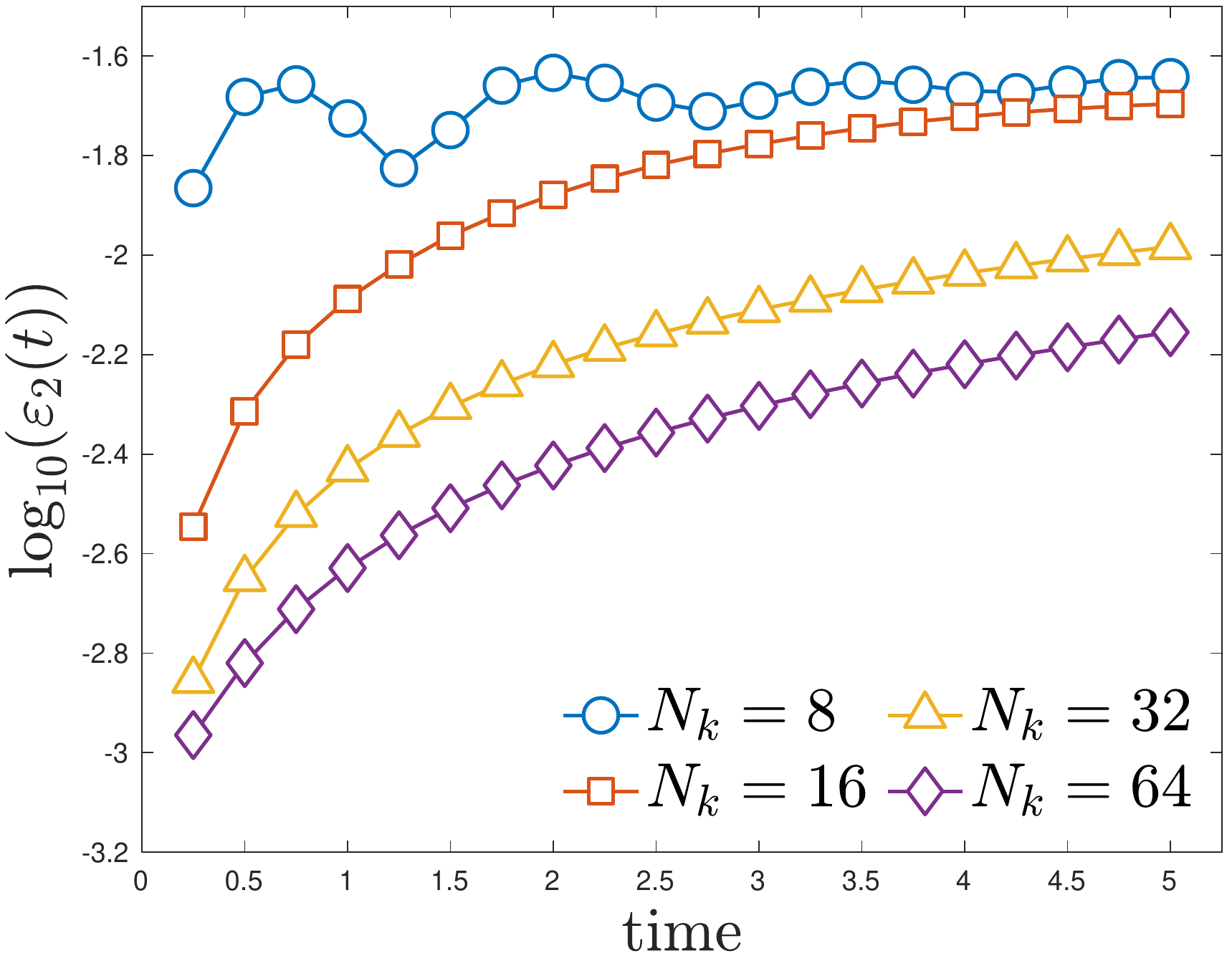}}}
\\
\centering
\subfigure[LPC1, deviation in total mass.\label{LPC_mass}]{
{\includegraphics[width=0.48\textwidth,height=0.27\textwidth]{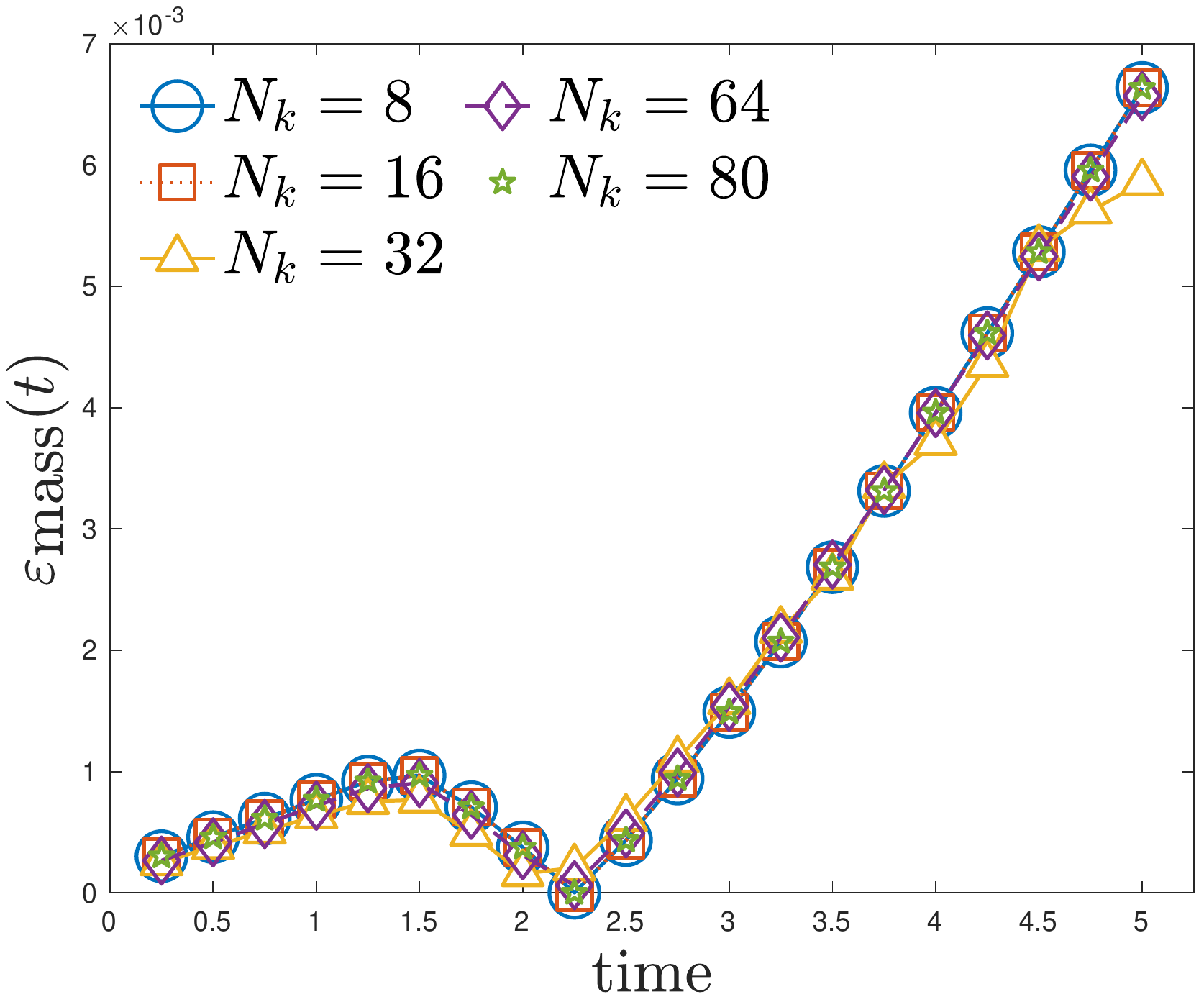}}}
\subfigure[OS, deviation in total mass. \label{OS_mass}]{
{\includegraphics[width=0.48\textwidth,height=0.27\textwidth]{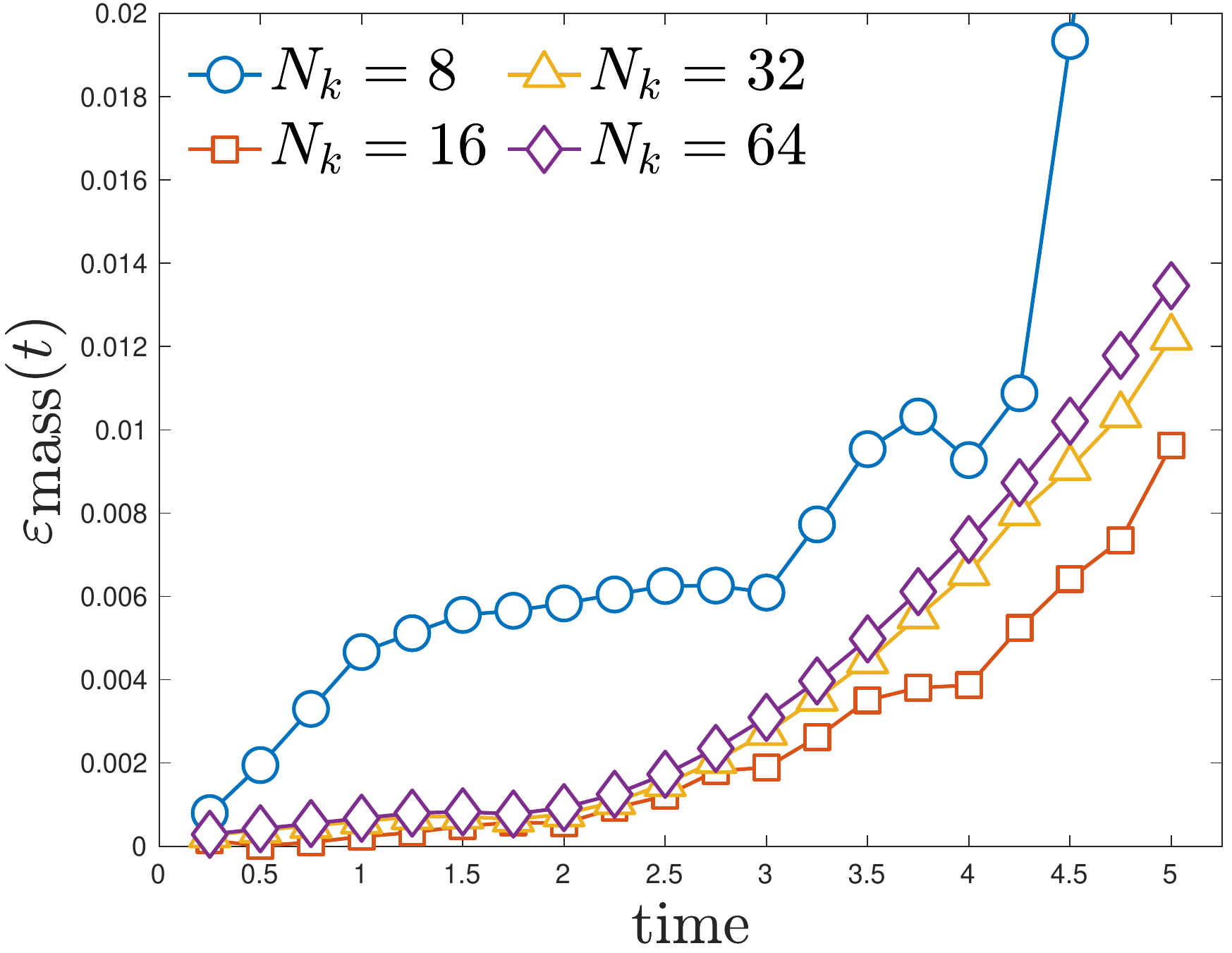}}}
\\
\centering
\subfigure[Convergence of $\varepsilon_{\infty}$ at $t=5$. \label{1s_convergnce_maxerr}]{
{\includegraphics[width=0.48\textwidth,height=0.27\textwidth]{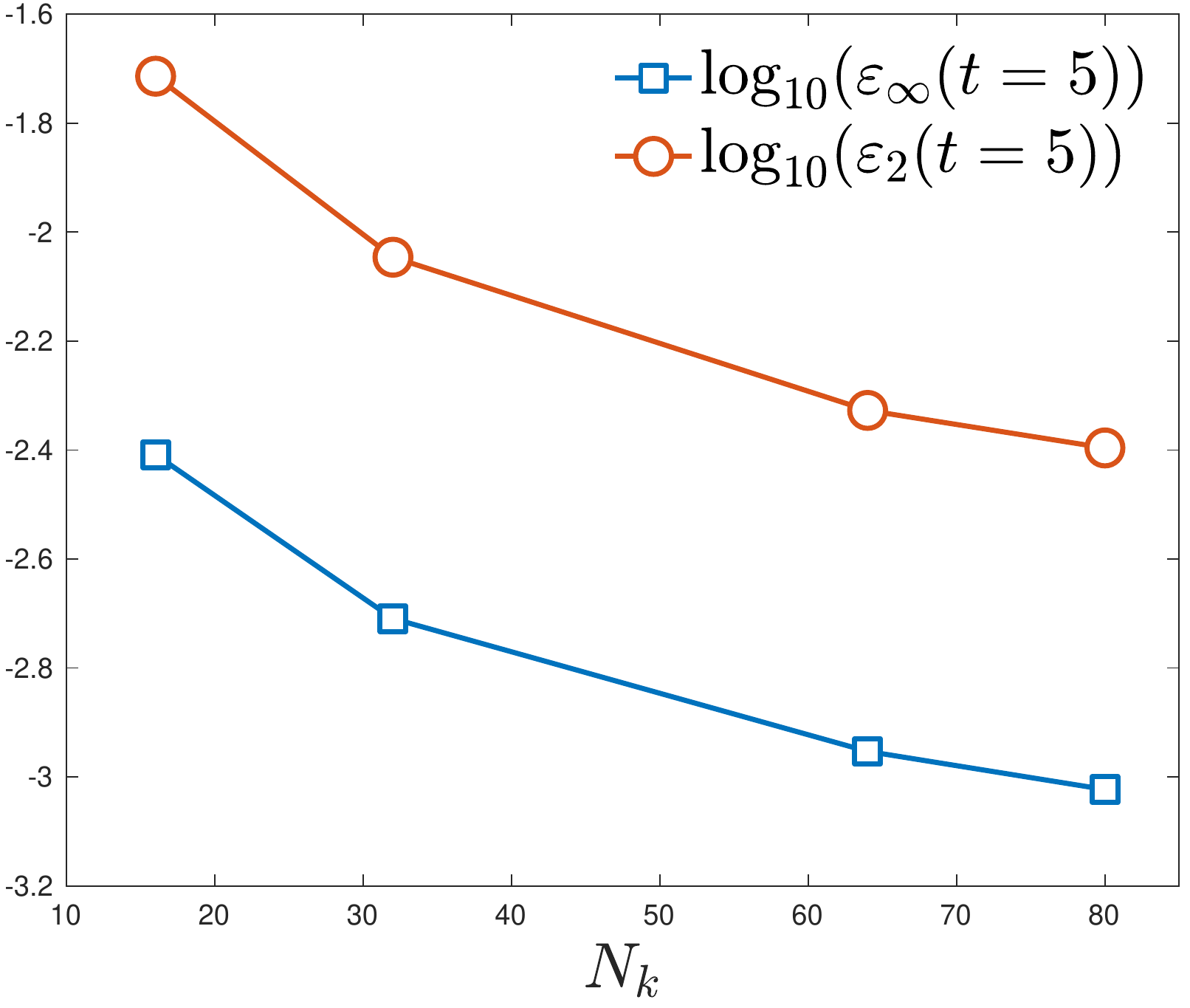}}}
\subfigure[Convergence of $\varepsilon_{2}$ at $t=5$. \label{1s_convergnce_L2err}]{
{\includegraphics[width=0.48\textwidth,height=0.27\textwidth]{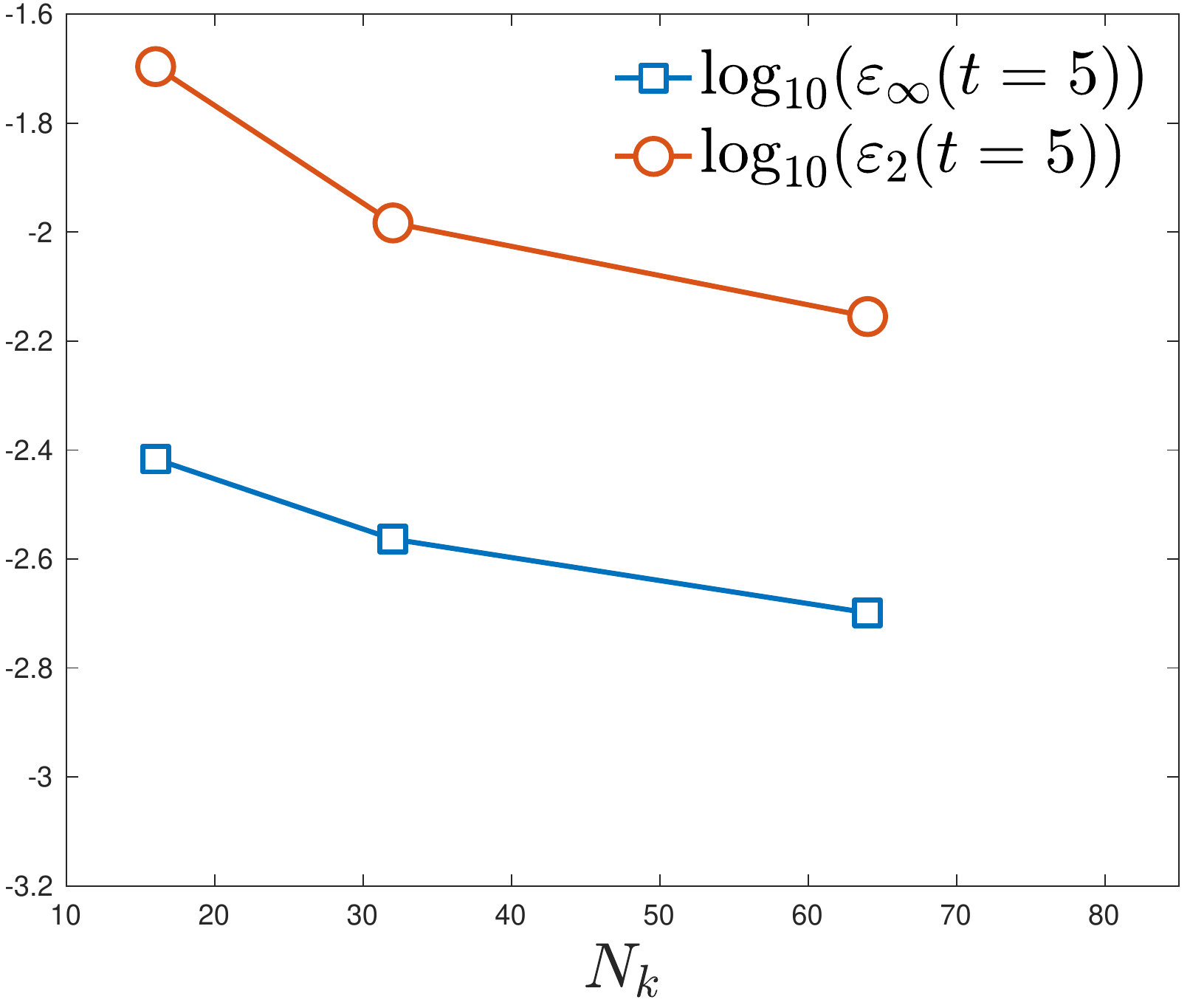}}}
\caption{\small The Hydrogen 1s Wigner function:  The performance of TKM under different $\Delta k$, with $\Delta x = 0.3$. The convergence of TKM is verified, albeit with lower convergence rate due to errors caused by the spline interpolation and truncation of $\bk$-space. In addition, LPC1 still outperforms OS in both accuracy and mass conservation.
 \label{1s_comp_integrator}}
\end{figure}

\begin{figure}[!h]
\centering
\subfigure[Time evolution of $\varepsilon_{\infty}$.]{
{\includegraphics[width=0.48\textwidth,height=0.27\textwidth]{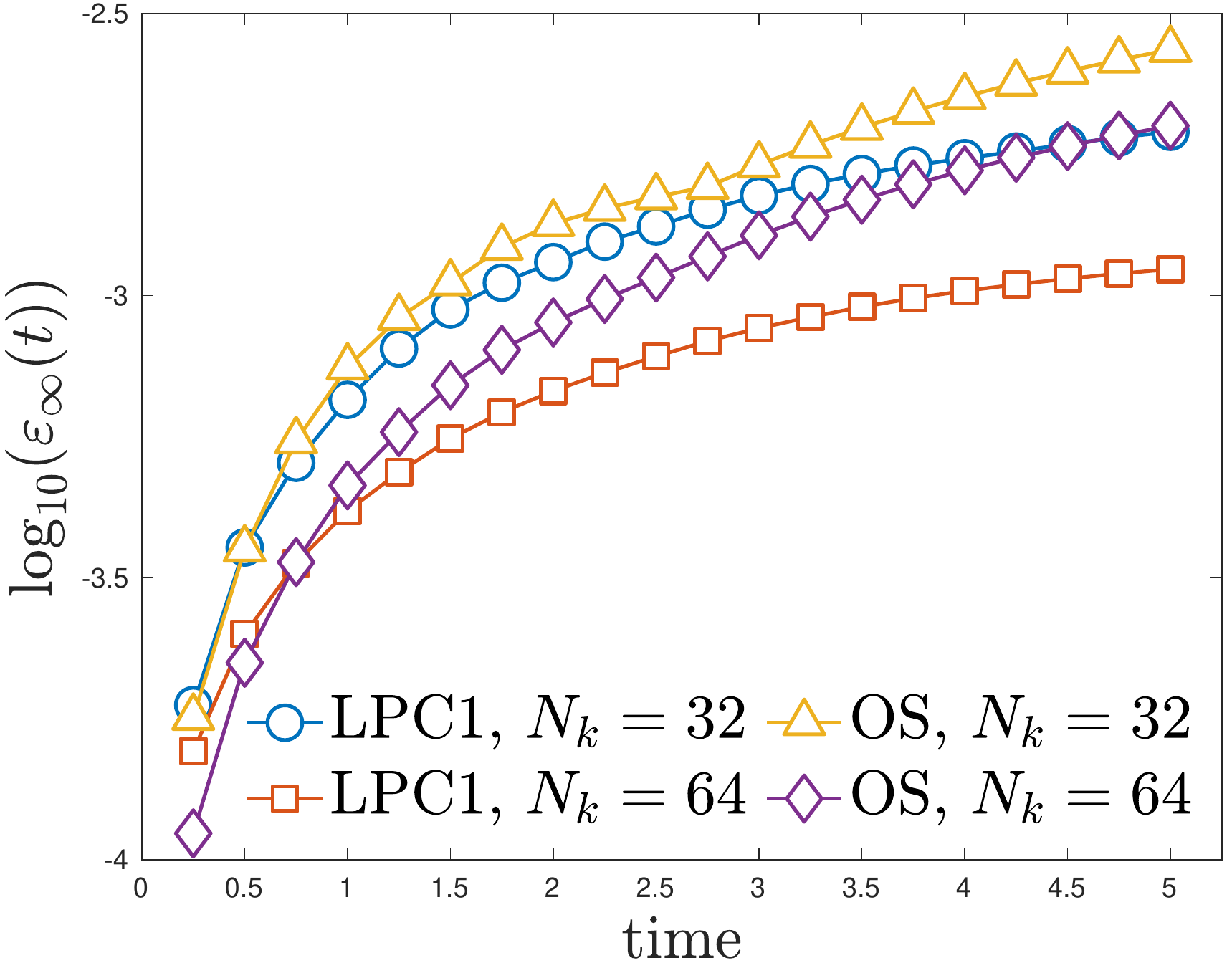}}}
\subfigure[Time evolution of  $\varepsilon_{2}$. ]{
{\includegraphics[width=0.48\textwidth,height=0.27\textwidth]{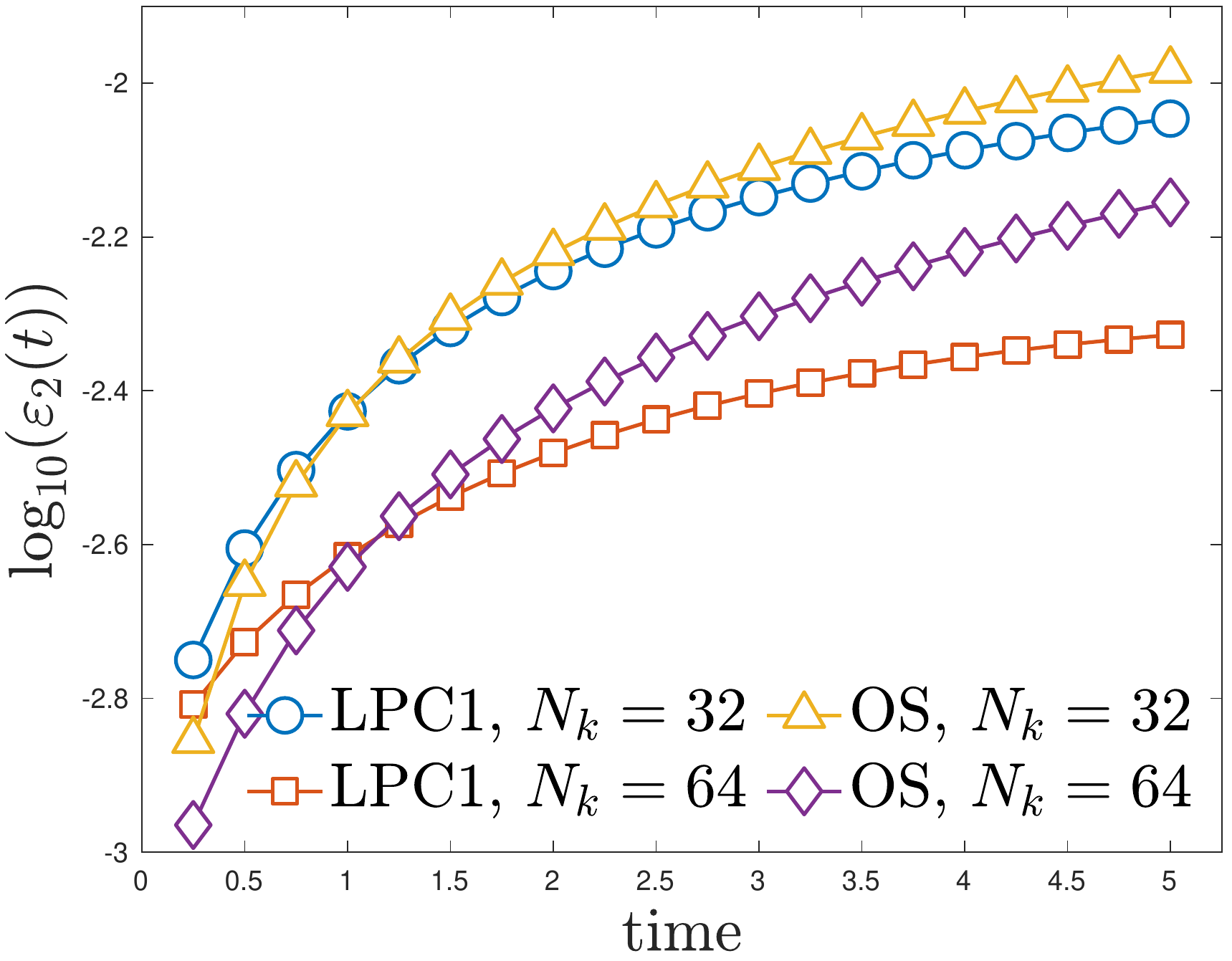}}}
\caption{\small The Hydrogen 1s Wigner function:  The non-splitting scheme outperforms the Strang splitting in accuracy. \label{1s_comp_error}}
\end{figure}

\begin{figure}[!h]
\centering
\subfigure[LPC1, $N_k = 32$.]{
{\includegraphics[width=0.48\textwidth,height=0.27\textwidth]{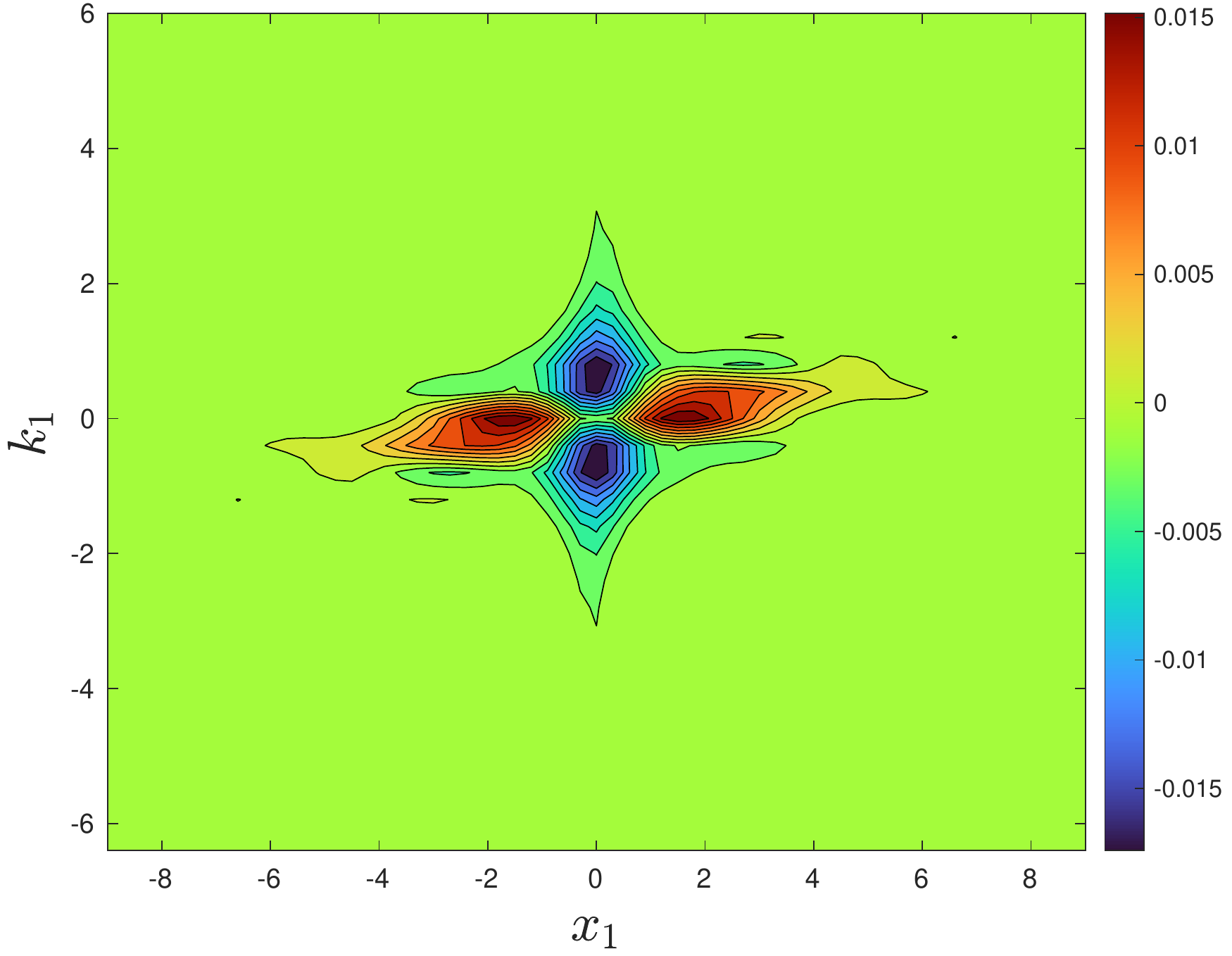}}}
\subfigure[OS, $N_k = 32$. ]{
{\includegraphics[width=0.48\textwidth,height=0.27\textwidth]{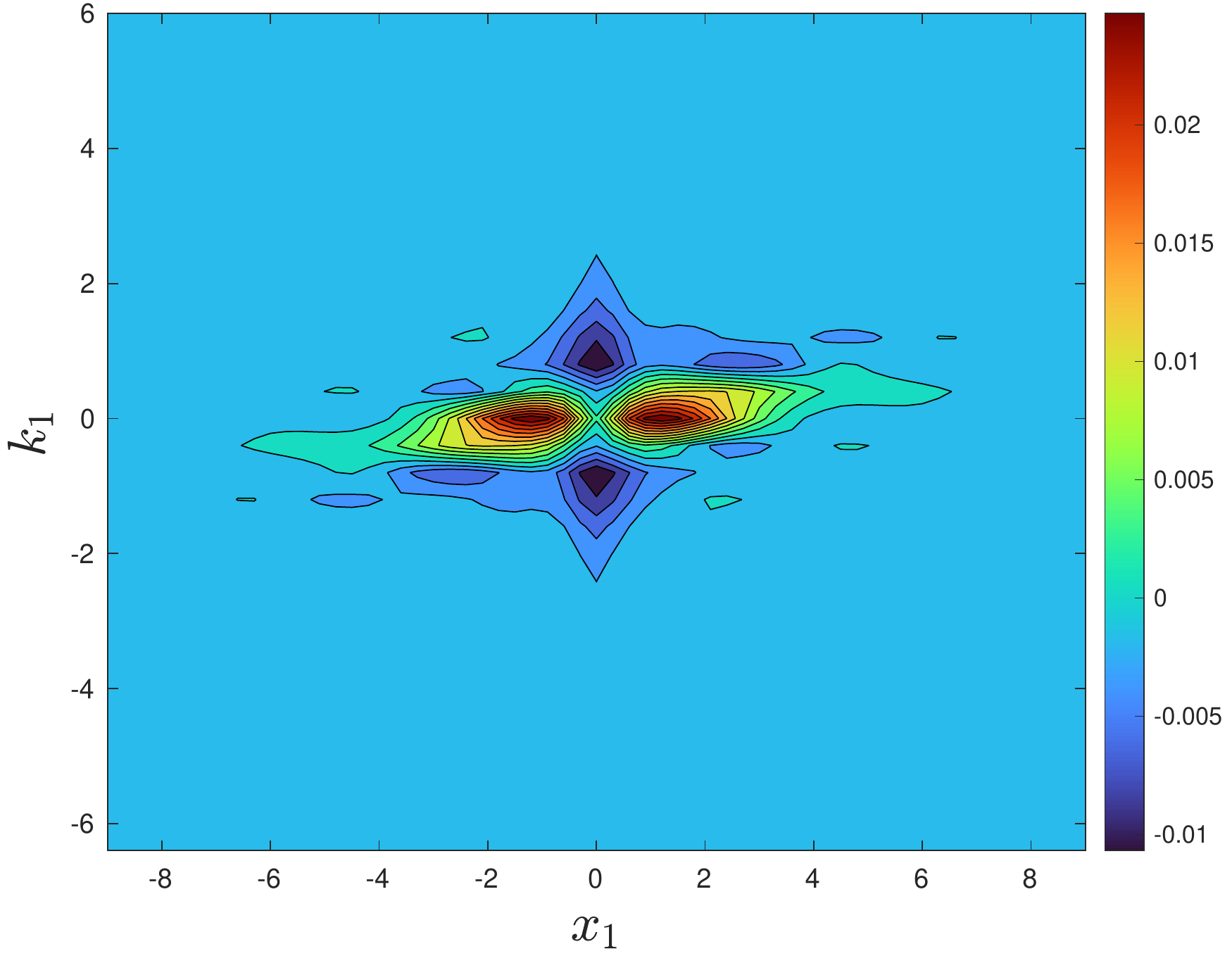}}}
\\
\centering
\subfigure[LPC1, $N_k =64$.\label{1s_error_visual_LPC_N64}]{
{\includegraphics[width=0.48\textwidth,height=0.27\textwidth]{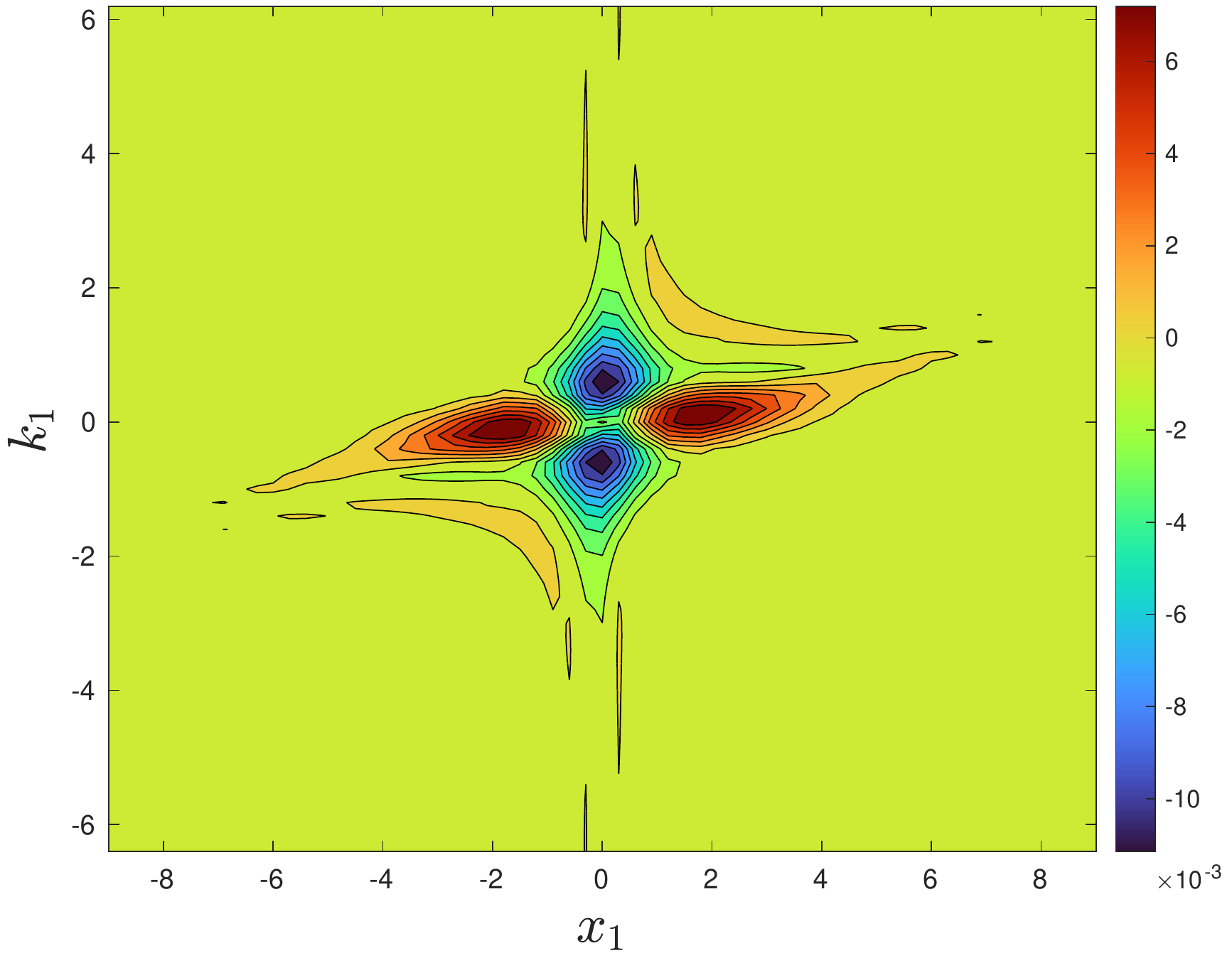}}}
\subfigure[OS, $N_k =64$. \label{1s_error_visual_OS_N64}]{
{\includegraphics[width=0.48\textwidth,height=0.27\textwidth]{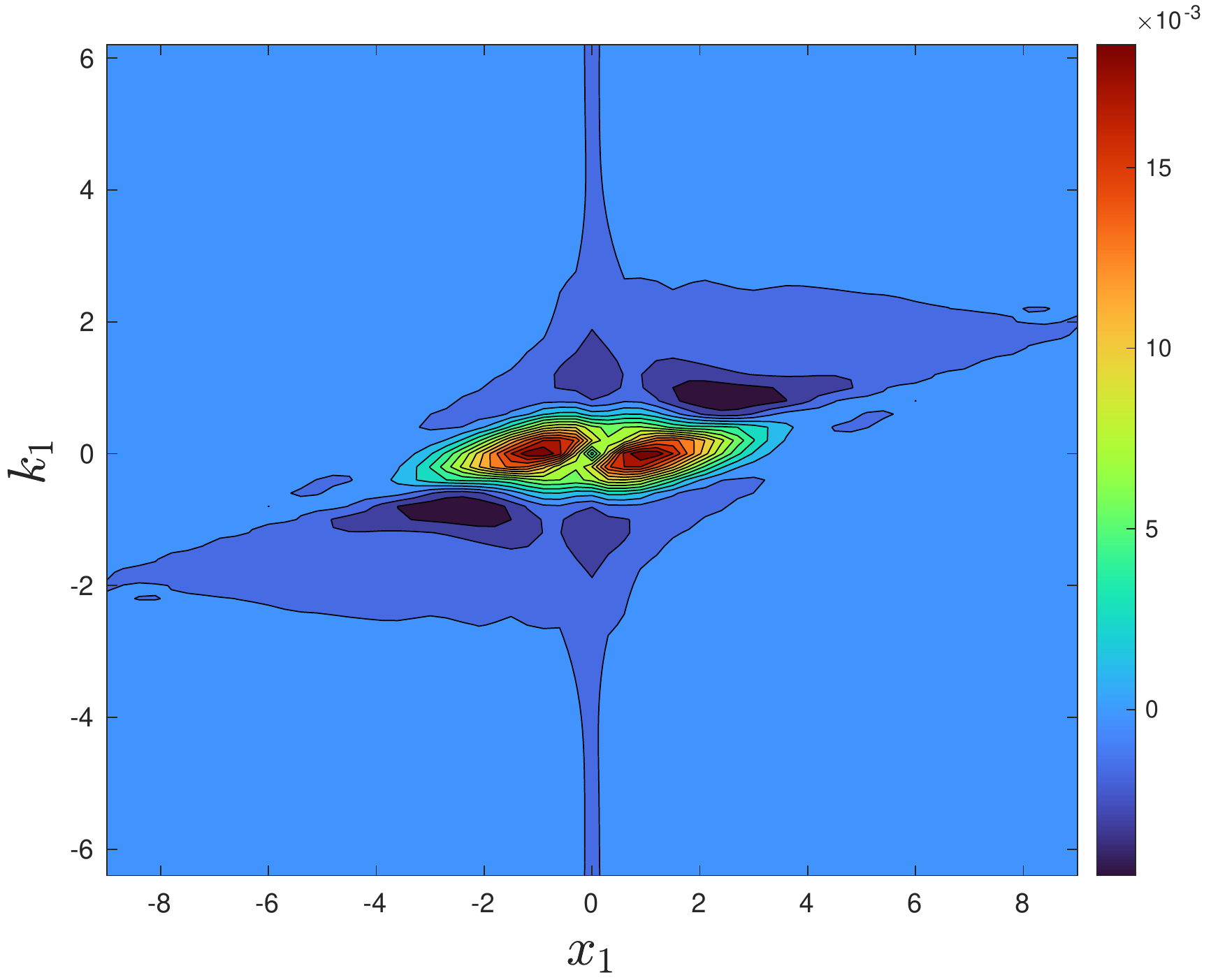}}}
\caption{\small Hydrogen 1s Wigner function:  A visualization of numerical errors $W_1^{\textup{num}}(x, k, t) - W_1^{\textup{ref}}(x, k, t)$ at $t = 5$ induced by truncation of $\bk$-space. Since the 1s Wigner function is not compactly supported in $[-6.4, 6.4]^3$, there are small errors found near the $\bk$-boundary, as well as near $\bx$-boundary due to the natural boundary condition.
\label{1s_error_visualization}}
\end{figure}

Again, a comparison of OS and LPC1 under different $\Delta k$, as well as the convergence in $\bk$-space, is presented in Figures \ref{1s_comp_integrator} and \ref{1s_comp_error}. Numerical errors for reduced Wigner function $W_1^{\textup{num}} - W_1^{\textup{ref}}$ under $N_k = 32$ and $N_k = 64$ are visualized in Figure \ref{1s_error_visualization}. From the results, we can make the following observations.

{\bf Convergence with respect to $\Delta k$}: The convergence of TKM is clearly verified in Figures \ref{1s_convergnce_maxerr} and \ref{1s_convergnce_L2err}, albeit its convergence rate is slower than expectation due to  the mixture of various error sources. Since the initial 1s Wigner function is not compactly supported in $[-6.4, 6.4]^3$ (see Figure \ref{1s_tail}), the overlap with the periodic image may produce small oscillations near the $\bk$-boundary, which is also visualized in Figures  \ref{1s_error_visual_LPC_N64} and \ref{1s_error_visual_OS_N64}.

{\bf Comparison of LPC1 and OS}: Nonetheless, with $61^3 \times 64^3$ uniform grid mesh and LPC1 integrator, CHASM can still achieve relative maximal error about $3.45\%$ and relative $L^2$-error about $7.41\%$ for the reduced Wigner function \eqref{def.reduced_Wigner_function} up to $T=5$, where $\max(|f_{\textup{1s}(\bx, \bk)}| = 1/\pi^3 \approx 0.0323$ and $\sqrt{\iint (|f_{\textup{1s}}(\bx, \bk)|^2 \D \bx \D \bk} \approx 0.0635$. When $N_k = 80$, the relative maximal error and $L^2$-error reduce to $2.93\%$ and $6.33\%$, respectively. By contrast, when the Strang splitting is adopted under the mesh size $61^3 \times 64^3$ , the relative maximal error is $6.20\%$ and relative $L^2$-error is about $11.02\%$. It is also clearly seen in Figure \ref{1s_comp_error} that the non-splitting scheme outperforms the splitting scheme for $N_k \ge 32$ in accuracy.

{\bf Mass conservation}: A slight deviation of the total mass is observed in time evolution. From Figures \ref{LPC_mass} and \ref{OS_mass}, one can see that $\varepsilon_{\textup{mass}}$ up to $t=5$ of LPC1 is about $0.66\%$, regardless of $N_k$, while that of OS is about $1.35\%$ when $N_k = 64$ and becomes even larger when $N_k$ goes smaller.

 \section*{Acknowledgement}
This research was supported by the National Natural Science Foundation of China (No.~1210010642), the Projects funded by China Postdoctoral Science Foundation (No.~2020TQ0011, 2021M690227) and the High-performance Computing Platform of Peking University. 
SS is partially supported by Beijing Academy of Artificial Intelligence (BAAI). The authors are sincerely grateful to  Haoyang Liu and Shuyi Zhang at Peking University for their technical supports on computing environment, which have greatly facilitated our numerical simulations.

\bibliographystyle{spmpsci}      

\begin{thebibliography}{10}
\providecommand{\url}[1]{{#1}}
\providecommand{\urlprefix}{URL }
\expandafter\ifx\csname urlstyle\endcsname\relax
  \providecommand{\doi}[1]{DOI~\discretionary{}{}{}#1}\else
  \providecommand{\doi}{DOI~\discretionary{}{}{}\begingroup
  \urlstyle{rm}\Url}\fi

\bibitem{ChenShaoCai2019}
Chen, Z., Shao, S., Cai, W.: {A high order efficient numerical method for 4-D
  Wigner equation of quantum double-slit interferences}.
\newblock J. Comput. Phys. \textbf{396}, 54--71 (2019).

\bibitem{CrouseillesEinkemmerMassot2020}
Crouseilles, N., Einkemmer, L., Massot, J.: Exponential methods for solving
  hyperbolic problems with application to collisionless kinetic equations.
\newblock J. Comput. Phys. \textbf{420}, 109688 (2020).

\bibitem{CrouseillesLatuSonnendrucker2006}
Crouseilles, N., Latu, G., Sonnendr{\"u}cker, E.: Hermite spline interpolation
  on patches for a parallel solving of the {V}lasov-{P}oisson equation  (2006).
\newblock RR-5926, INRIA.

\bibitem{CrouseillesLatuSonnendrucker2009}
Crouseilles, N., Latu, G., Sonnendrücker, E.: A parallel {V}lasov solver based
  on local cubic spline interpolation on patches.
\newblock J. Comput. Phys. \textbf{228}, 1429--1446 (2009).

\bibitem{DimarcoLoubereNarskiRey2018}
Dimarco, G., Loub$\grave{\textup{e}}$re, R., Narski, J., Rey, T.: An efficient
  numerical method for solving the {B}oltzmann equation in multidimensions.
\newblock J. Comput. Phys. \textbf{353}, 46--81 (2018).

\bibitem{Goudon2002}
Goudon, T.: Analysis of a semidiscrete version of the {Wigner} equation.
\newblock SIAM J. Numer. Anal. \textbf{40}, 2007--2025 (2002).

\bibitem{GreengardJiangZhang2018}
Greengard, L., Jiang, S., Zhang, Y.: The anisotropic truncated kernel method
  for convolution with free-space {G}reen's functions.
\newblock SIAM J. Sci. Comput. \textbf{40}, A3733--A3754 (2018).

\bibitem{GuoLiWang2018b}
Guo, X., Li, Y., Wang, H.: A high order finite difference method for tempered
  fractional diffusion equations with applications to the {CGMY} model.
\newblock SIAM J. Sci. Comput. \textbf{40}, A3322--A3343 (2018).

\bibitem{Kormann2015}
Kormann, K.: A semi-{L}agrangian {V}lasov solver in tensor train format.
\newblock SIAM J. Sci. Comput. \textbf{37}, B613--B632 (2015).

\bibitem{KormannReuterRampp2019}
Kormann, K., Reuter, K., Rampp, M.: A massively parallel semi-{L}agrangian
  solver for the six-dimensional {Vlasov–Poisson} equation.
\newblock Int. J. High Perform. C. \textbf{33}, 924--947 (2019).

\bibitem{MalevskyThomas1997}
Malevsky, A.V., Thomas, S.J.: Parallel algorithms for semi-{L}agrangian
  advection.
\newblock Int. J. Numer. Meth. Fl. \textbf{25}, 455--473 (1997).

\bibitem{PraxmeyerMostowskiWodkiewicz2005}
Praxmeyer, L., Mostowski, J., Wódkiewicz, K.: Hydrogen atom in phase space:
  The {W}igner representation.
\newblock J. Phys. A: Math. Gen. \textbf{39}, 14143--14151 (2005).

\bibitem{Ringhofer1990}
Ringhofer, C.: A spectral method for the numerical simulation of quantum
  tunneling phenomena.
\newblock SIAM J. Numer. Anal. \textbf{27}, 32--50 (1990).

\bibitem{VicoGreengardFerrando2016}
Vico, F., Greengard, L., Ferrando, M.: Fast convolution with free-space
  {G}reen's functions,.
\newblock J. Comput. Phys. \textbf{323}, 191--203 (2016).

\bibitem{XiongChenShao2016}
Xiong, Y., Chen, Z., Shao, S.: An advective-spectral-mixed method for
  time-dependent many-body {Wigner} simulations.
\newblock SIAM J. Sci. Comput. \textbf{38}, B491--B520 (2016).

\end{thebibliography}

\end{document}